\documentclass[pdflatex,sn-mathphys-num]{sn-jnl}

\usepackage{graphicx}%
\usepackage{multirow}%
\usepackage{amsmath,amssymb,amsfonts}%
\usepackage{amsthm}%
\usepackage{mathrsfs}%
\usepackage[title]{appendix}%
\usepackage{xcolor}%
\usepackage{textcomp}%
\usepackage{manyfoot}%
\usepackage{booktabs}%
\usepackage{algorithm}%
\usepackage{algorithmicx}%
\usepackage{algpseudocode}%
\usepackage{listings}%
\usepackage{verbatim}%
\usepackage{natbib}%

\theoremstyle{thmstyleone}%
\newtheorem{theorem}{Theorem}
\newtheorem{proposition}[theorem]{Proposition}%

\theoremstyle{thmstyletwo}%

\theoremstyle{thmstylethree}%

\definecolor{gre}{rgb}{0,.75,0}

\definecolor{red}{rgb}{1,0,0}
\def\red#1{{\color{red}#1}}
\def\norm#1{\|#1\|}
\def\argmin{\mathop{\rm arg\,min}}

\def\R{\mathbb{R}}

\begin{document}

\title[Greedy Newton]{Greedy Newton: Newton's Method with Exact Line Search}


\author[1]{\fnm{Betty} \sur{Shea}}\email{sheaws@cs.ubc.ca}

\author[1,2]{\fnm{Mark} \sur{Schmidt}}\email{schmidtm@cs.ubc.ca}

\affil[1]{\orgname{The University of British Columbia}, \country{Canada}}

\affil[2]{\orgname{Canada CIFAR AI Chair (Amii)}}

\abstract{A defining characteristic of Newton's method is local superlinear convergence within a neighbourhood of a strict local minimum. However, outside this neighborhood Newton's method can converge slowly or even diverge. A common approach to dealing with non-convergence is using a step size that is set by an Armijo backtracking line search. With suitable initialization the line-search preserves local superlinear  convergence, but may give sub-optimal progress when not near a solution. In this work we consider Newton's method under an exact line search, which we call ``greedy Newton'' (GN). We show that this leads to an improved global convergence rate, while retaining a local superlinear convergence rate. We empirically show that GN may work better than backtracking Newton by allowing significantly larger step sizes.}


\maketitle

\section{Introduction}

For minimizing a twice-differentiable  function $f:\R^n\rightarrow\R$, 
the pure Newton iteration $x_{k+1}^N$ starting from some vector $x_k$ is given by 
\begin{equation}
	\label{eq:newton_update}
	x_{k+1}^N=x_{k}-\nabla^2f(x_{k})^{-1}\nabla f(x_{k}).
\end{equation}
This method dates back to work by Newton and Raphson in the 1600s for finding roots of polynomials~\cite{deuflhard2012short}. Since then, the method has evolved to work in a variety of settings. In optimization, Newton's method is a powerful tool for minimizing non-linear objectives mainly due to its remarkable property of superlinear (quadratic) convergence in a neighborhood of a strict local minimum under appropriate conditions~\citep[see][Section~9.5]{Boyd2004}. 

However, Newton's method also has known weaknesses. For example, the method is not guaranteed to converge in general or even decrease $f$. One of the standard fixes to the non-convergence is to introduce a step size $\alpha_k$,
\begin{equation}
	x_{k+1} =x_{k}-\alpha_k\nabla^2f(x_{k})^{-1}\nabla f(x_{k}).
	\label{eq:GN_update}
\end{equation}
The step size is typically set by first considering $\alpha_k=1$ and dividing the step size by a fixed constant (``backtracking'') until the Armijo condition is satisfied~\citep[][Section~9.5]{Boyd2004}. Provided that we eventually get close enough to a strict local minimum, Armijo backtracking preserves the superlinear convergence of the pure Newton method. Many variations of Newton's method exist such as modifications for cases where $\nabla^2 f(x_k)$ is not invertible. Other variations include those based on trust-region methods instead of line searches~\cite[see][]{Nocedal2006}, but the superlinear convergence proofs in the literature that we are aware of for~\eqref{eq:GN_update} assume we first test $\alpha_k=1$ and accept this step size if it satisfies a variant of the Armijo condition.

While $\alpha_k=1$ becomes asymptotically optimal in the neighbourhood of a strict local minimizer, it may not be the optimal step size even when close to a minimizer. Further, when far away from a local minimizer using $\alpha_k=1$ or the smaller values obtained by backtracking may converge slowly.
This paper instead investigates Newton's method where the step size is set to minimize the function value,
\begin{equation}
	\label{eq:newton_exactLS}
	\alpha_k\in\argmin_{\alpha}f(x_k-\alpha\nabla^2f(x_k)^{-1}\nabla f(x_k)).
\end{equation}
We call using this exact line search within Newton's method the ``greedy Newton'' (GN) method. We first address two mis-conceptions the reader may have about this method:
\begin{itemize}
	\item It typically {\bf does not significantly increase the cost of Newton's method} to find a local minimizer of~\eqref{eq:newton_exactLS}. For most problems just computing the Hessian is $n$-times more expensive than evaluating the function or directional derivatives. Thus, you can evaluate the objective in~\eqref{eq:newton_exactLS} and its derivative several times without changing the overall cost of the method. As an example, consider logistic regression with $m$ training examples $\{a_i,b_i\}$ with dense features $a_i\in\R^n$ and binary labels $b_i\in\{-1,1\}$,
	\begin{equation}
		f(x) = \sum_{i=1}^m\log(1+\exp(-b_ix^Ta_i)).
		\label{eq:logreg}
	\end{equation}
	The cost of computing the Hessian for this problem is $O(mn^2)$, but the cost of evaluating $f$ or a directional derivative of it is only $O(mn)$. With bisection we can solve the one-dimensional problem~\eqref{eq:newton_exactLS} to $\epsilon$ accuracy over a bounded domain in $O(\log(1/\epsilon))$ iterations, so the cost of a naive black-box numerical method is only $O(mn\log(1/\epsilon))$. 
	
	Further, if we exploit the linear composition structure of ~\eqref{eq:logreg} the cost of bisection can be reduced to $O(mn + m\log(1/\epsilon))$. It is also possible to use faster one-dimensional minimizers like the secant method. Indeed,  low-cost line searches are possible for a wide variety of problems including linear models, matrix factorization models, and certain neural networks~\citep{Narkiss2005,Narkiss2005_SVM,Zibulevsky2008,Zibulevsky2010,sorber2016exact,Shea2023}.
	\item The {\bf exact line search can yield a significantly smaller function value} than Armijo backtracking. The optimal step size~\eqref{eq:newton_exactLS} may be significantly larger than the maximum step size of $1$ considered in standard implementations. While we could backtrack from a step size larger than $1$, the Armijo condition itself can exclude the optimal step size. Indeed, for non-quadratic functions the \emph{maximum step size allowed by the Armijo condition can be arbitrarily smaller} than the optimal step size.
\end{itemize}

By modifying standard arguments for the convergence of Newton's method, we show:
\begin{enumerate}
	\item For strongly-convex functions with a Lipschitz-continuous gradient, GN slightly improves the global convergence rate compared to Armijo backtracking (Section~2.1).
	\item Under the additional assumption that the Hessian is Lipschitz-continuous, superlinear convergence is achieved by any method (like GN) that decreases the function value by at least as much as the pure Newton iteration~\eqref{eq:newton_update} (Section~2.2).
	\item Under the same assumptions, we analyze the local convergence rate of Newton's method with arbitrary step sizes (Section~2.3).
	\item The global convergence rate of Netwon's method can be further improved using hybrids with gradient descent  (Section~2.4).
\end{enumerate}
We are not aware of the superlinear convergence of GN appearing previously in the literature, although recent work bounds the GN step size for self-concordant functions~\cite{ivanova2023optimal} and various progress measures  give superlinear convergence for solving non-linear equations~\cite{burdakov1980some}. We also note that asymptotic superlinear convergence of quasi-Newton methods with an exact line search is known classically~\citep{powell1971convergence} and concurrent with this work explicit superlinear rates of greedy quasi-Newton methods have been derived~\citep{jin2024non}.
In Section \ref{sec:GN_experiments} we experiment with the GN method for logistic regression.
Our findings suggest that GN consistently works better than using Armijo backtracking, and substantially better for certain problems where the optimal step sizes can be much larger than 1. 

\section{Convergence of Greedy Newton Methods}
\label{sec:GN_analysis}

All our results assume that $f$ is twice-differentiable and that the eigenvalues of the Hessian $\nabla^2 f$ are bounded between positive constants $\mu$ and $L$ for all $x$,
\begin{equation}
	\mu I \preceq \nabla^2 f(x) \preceq L I.
	\label{eq:Lmu}
\end{equation}
These assumptions are equivalent to assuming that $\nabla f$ is $L$-Lipschitz continuous and that $f$ is $\mu$-strongly convex. We note that without these assumptions GN may not converge~\cite{mascarenhas2007divergence,jarre2016simple}, but various Hessian modifications guarantee convergence~\cite[see][Section~3.4]{Nocedal2006}.
The local superlinear convergence results also require that the Hessian is $M$-Lipschitz continuous,
\begin{equation}
	\norm{\nabla^2 f(x) - \nabla^2f(y)} \leq M\norm{x-y},
	\label{eq:M}
\end{equation}
where the matrix norm $\norm{\cdot}$ on the left is the spectral norm while the vector norm on the right is the Euclidean norm. We give proofs of the results in this section in Appendix~\ref{sec:app_analysis}.

\subsection{Global Convergence of Greedy Newton}
\label{sec:global}

We first give a global rate of convergence for the GN method.
\begin{proposition}
	Let a twice-differentiable $f$ be $\mu$-strongly convex with an $L$-Lipschitz continuous gradient~\eqref{eq:Lmu}. Then the iterations of Newton's method~\eqref{eq:GN_update} with the greedy step size~\eqref{eq:newton_exactLS} satisfy
	\[
	f(x_k) - f(x_*) \leq \left(1 - \frac{\mu^2}{L^2}\right)^k[f(x_0) - f(x_*)].
	\]
\end{proposition}
This result implies that in order for the sub-optimality $d_k = (f(x_k) - f(x_*))$ to be less than $\epsilon$, we require at most $(L^2/\mu^2)\log(d_0/\epsilon)$ iterations.
If we instead set the step size by starting from a sufficiently large guess for $\alpha_k$ and halving it until the Armijo condition is satisfied, then with a sufficient decrease factor of $\sigma=1/2$ we have a slower rate of
\[
f(x_k) - f(x_*) \leq \left(1 - \frac{\mu^2}{\red{2}L^2}\right)^k[f(x_0) - f(x_*)].
\]
This requires
$\red{2}(L^2/\mu^2)\log(d_0/\epsilon)$ iterations to guarantee that we reach an accuracy of $\epsilon$. Thus, GN halves the worst-case number of steps compared to this standard approach. If we use an Armijo sufficient decrease factor of $\sigma < 1/2$ and multiply the step size by $\beta < 1$ instead of $1/2$ when we backtrack, we require $(L^2/\red{2\sigma\beta}\mu^2)\log(d_0/\epsilon)$~\cite[see][Section~9.5]{Boyd2004} (this again assumes the initial guess for $\alpha_k$ is sufficiently large, and note that it may need to be larger than 1). Note that the extra factor of $1/2\sigma\beta$ is greater than 1 since $\sigma\beta < 1/2$.
Thus, GN performs as well as backtracking with an aribtrarily large initial $\alpha_k$, and an arbitrarily small backtracking and sufficient decrease factor.

\subsection{Local Convergence of ``As Fast as Newton'' Methods}
\label{sec:local}

We next consider a local rate for any method that decreases the function as much as the pure Newton method.
\begin{proposition}
	Let a twice-differentiable $f$ be $\mu$-strongly convex with an $L$-Lipschitz continuous gradient~\eqref{eq:Lmu}, and an $M$-Lipschitz continuous Hessian~\eqref{eq:M}. Consider a method that is guaranteed to decrease the function as much as the pure Newton step~\eqref{eq:newton_update}, $f(x_{k+1}) \leq f(x_k - \nabla^2 f(x_k)^{-1}\nabla f(x_k))$. The iterations of such methods satisfy
	\[
	\norm{x_{k+1}-x_*}\leq\red{\sqrt{\frac{L}{\mu}}}\frac{M}{2\mu}\norm{x_k-x_*}^2.
	\]
\end{proposition}
This result implies superlinear (quadratic) convergence beginning at the first iteration where we have
\begin{equation}
	\label{eq:easyGN_fastRegionLS}
	\norm{x_k-x_*}<\red{\sqrt{\frac{\mu}{L}}}\frac{2\mu}{M}.
\end{equation}
Note that this radius of fast convergence is smaller than the radius for the pure Newton method by a factor of $\sqrt{\mu/L}$~\cite{Nesterov2006,Sun2021}, and thus we must be closer to the solution in order to guarantee superlinear convergence. Note that this result applies not only to the GN method but a variety of other possible methods.

\subsection{Local Convergence of Newton with Arbitrary Step Size}
\label{sec:non-unit}

We next consider a similar result, but for Newton's method with arbitrary step sizes.
\begin{proposition}
	Let a twice-differentiable $f$ be $\mu$-strongly convex with an $L$-Lipschitz continuous gradient~\eqref{eq:Lmu} and an $M$-Lipschitz continuous Hessian~\eqref{eq:M}. Then Newton's method with a step size of $\alpha_k$~\eqref{eq:GN_update} satisfies
	\[
	\norm{x_{k+1}-x_*}\leq\red{|\alpha_k|}\frac{M}{2\mu}\norm{x_k-x_*}^2+\red{\left|\alpha_k-1\right|\frac{L}{\mu}\norm{x_k-x_*}}
	\]
\end{proposition}
Note that if we assume \red{$|\alpha_k-1| \leq \norm{x_k-x_*}$} then we have
\[
\norm{x_{k+1}-x_*}\leq\red{|\alpha_k|}\frac{M}{2\mu}\norm{x_k-x_*}^2+\red{\frac{L}{\mu}\norm{x_k-x_*}^2} = \frac{\red{|\alpha_k|}M + \red{ 2L}}{2\mu}\norm{x_k-x_*}^2.
\]
Thus we have superliner (quadratic) convergence if for all $k$ large enough we have
\[
\red{|\alpha_k - 1| \leq \norm{x_k-x_*}} \quad \text{and} \quad
\norm{x_k-x_*} <
\frac{2\mu}{\red{|\alpha_k|}M + \red{2L}}.
\]
Thus, if $L$ is similar to $M$ and if $\alpha_k$ converges to 1 at least as fast as $\norm{x_k-x_*}$ converges to zero, then Newton's method with non-zero step sizes has a similar radius of fast convergence to the pure Newton method. In the specific case of GN we have that $\alpha_k$ converges to 1 asymptotically as the quadratic approximation in the pure Newton method becomes exact. But the rate that $\alpha_k$ converges to 1 is less clear. 

\subsection{Global Convergence of Hybrid Gradient-Newton Methods}

In Section~\ref{sec:global} we review how GN improves on the linear convergence rate of Newton's method with backtracking from $(1-2\beta\sigma\mu^2/L^2)$ to $(1-\mu^2/L^2)$. However, under the same assumptions gradient descent with an exact line search achieves a rate of $(1-\mu/L)$ while with backtracking gradient descent achieves a rate $(1-2\beta\sigma\mu/L)$~\cite[see][Section 9.3]{Boyd2004}. Fortunately, it is possible to use the result of Section~\ref{sec:local} to design methods that have these faster global linear convergence rates while maintaining a local superlinear convergence rate.

Perhaps the simplest hybrid method is the following:
\begin{itemize}
	\item Let $x_{k+1}^N$ be the pure Newton step~\eqref{eq:newton_update} and $x_{k+1}^G$ be the gradient descent step with exact line search,
	\[
	x_{k+1}^G = x_k - \alpha_k^G\nabla f(x_k), \quad \alpha_k^G \in \argmin_{\alpha}\left\{f(x_k - \alpha\nabla f(x_k))\right\}
	\]
	\item If $f(x_{k+1}^G) < f(x_{k+1}^N)$ take the gradient step, otherwise take the pure Newton step.
\end{itemize}
This approach guarantees the $(1-\mu/L)$ linear rate is achieved at all iterations, while the result of Section~\ref{sec:local} guarantees that this approach has a superlinear convergence rate. However, in our experiments this hybrid approach tended to perform worse than GN.

Other hybrid methods are possible, such as ones based on backtracking for either the gradient or Newton step. Another option discussed informally in the literature~\citep{Conn1994,zibulevsky2013speeding} is to use a step size on both the gradient and Newton step,
\[
x_{k+1} = x_k - \alpha_k^a\nabla f(x_k) - \alpha_k^b\nabla^2f(x_k)^{-1}\nabla f(x_k),
\]
and optimize the step sizes $\alpha_k^a$ and $\alpha_k^b$. This ``plane search'' approach to setting two step sizes is efficient for many problems arising in machine learning~\cite[see][]{Shea2023}. However, in our experiments we found that this approach only gave small gains over the basic GN method (with $\alpha_k^a$ consistently being chosen close to zero).

\section{Experiments}
\label{sec:GN_experiments}

Our first experiment considers logistic regression~\eqref{eq:logreg} with the synthetic data included in the minFunc package~\cite{Schmidt2005}. This generates $m=500$ examples where the elements of $a_i$ are sampled from a standard normal, a true $\tilde{x}$ is sampled from a standard normal, and we set $b_i$ to be the sign of $(a_i^T\tilde{x} + \delta_i)$ with $\delta_i$ sampled from a standard normal. We generated 4 versions: one with $n=20$ yielding a strongly-convex problem, one with $n=20$ where 10 of the features are repeated yielding a convex problem, one with $n=200$ yielding a strictly convex problem, and one with $n=2000$ yielding a convex problem. In the latter two cases the data is linearly separable, and thus no minimizing solution $x^*$ exists. In these separable cases the function converges to its minimal value of 0 for vectors of the form $\kappa x$ as $\kappa$ goes to $\infty$ for any vector $x$ that separates the data. For the convex cases where the Hessian is not positive-definite we used $\nabla^2f(x_k)+10^{-12}I$ in place of the Hessian. We also considered L2-regularized variants of these problems with a regularization strength of $\lambda=1$ (this makes all the problems strongly-convex with unique finite minimizers). We initialized all methods with the zero vector. For methods using the Armijo condition we used the standard choices of a sufficient decrease factor of $\sigma=10^{-4}$ and a backtracking factor of $\beta=1/2$. Although faster methods are possible, our experiments approximated an exact line search with the following simple bisection procedure for a descent direction $d_k$:
\begin{enumerate}
	\item We initialize with $\alpha=1$.
	\item While the directional derivative $\nabla f(x_k + \alpha d_k)^Td_k$ is negative, we double $\alpha$. Once we have an $\alpha$ large enough that the directional derivative is positive, this gives an interval $(0,\alpha)$ containing a minimizer.
	\item Compute the mid-point of the interval, and cut the interval in half based on the directional derivative at the mid-point.
	\item Repeat the previous step until the interval length is less than $10^{-8}$.
\end{enumerate}
We give the precise code used to implement GN for logistic regression in Appendix~\ref{app:code}.

\begin{figure}
	\includegraphics[width=.24\textwidth]{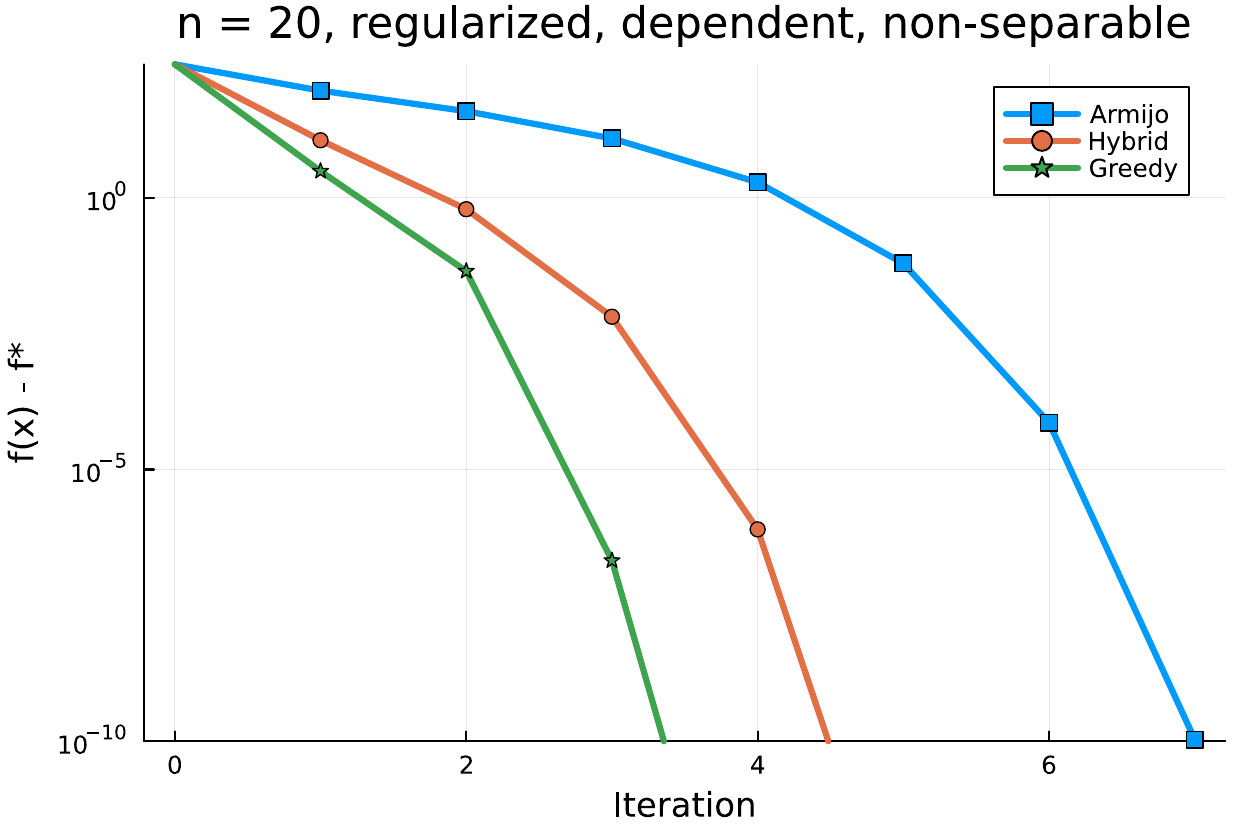}
	\includegraphics[width=.24\textwidth]{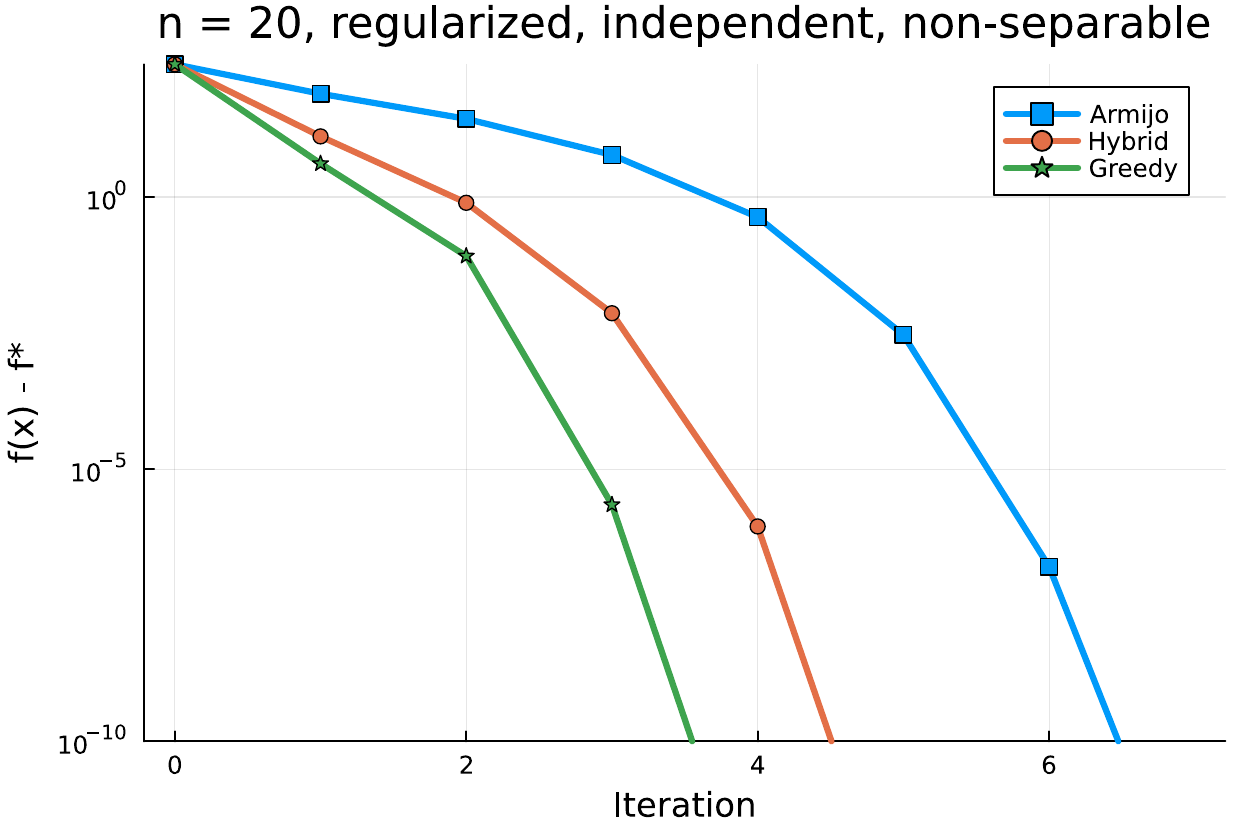}
	\includegraphics[width=.24\textwidth]{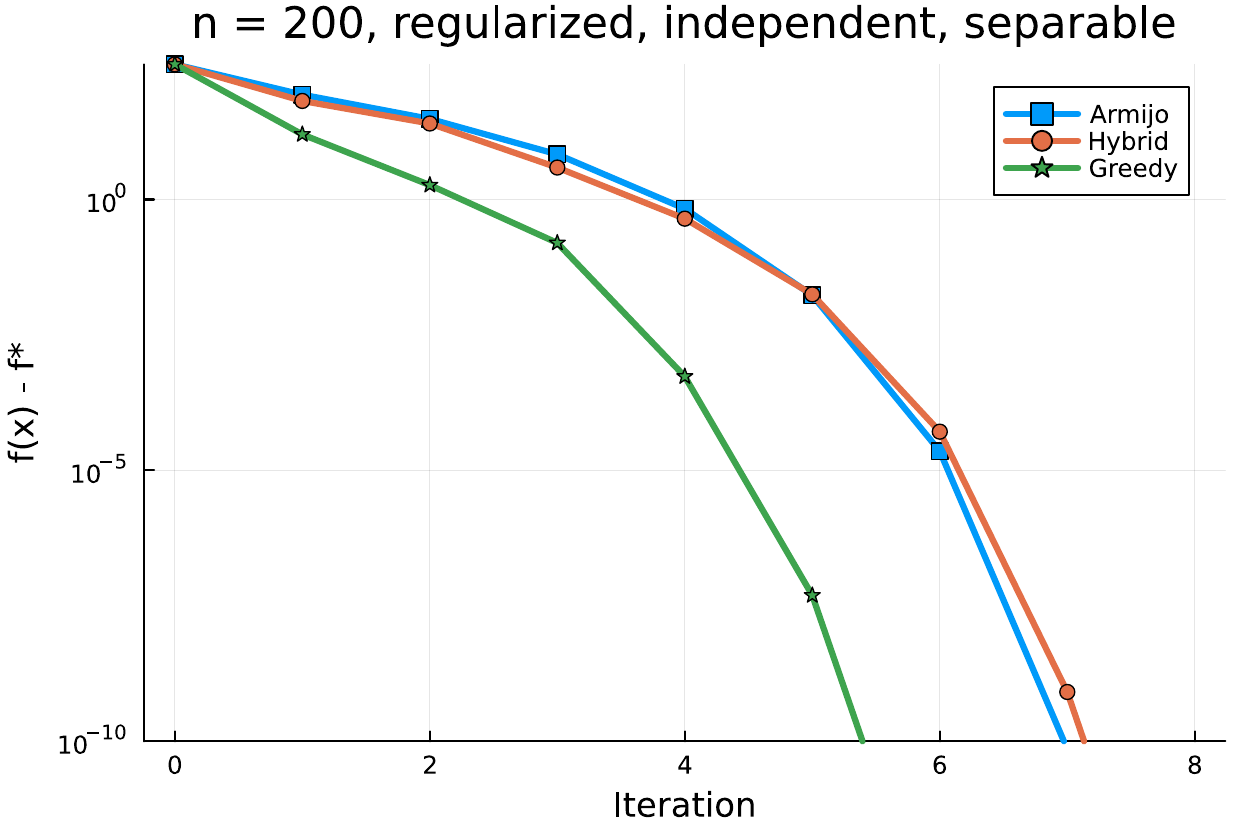}
	\includegraphics[width=.24\textwidth]{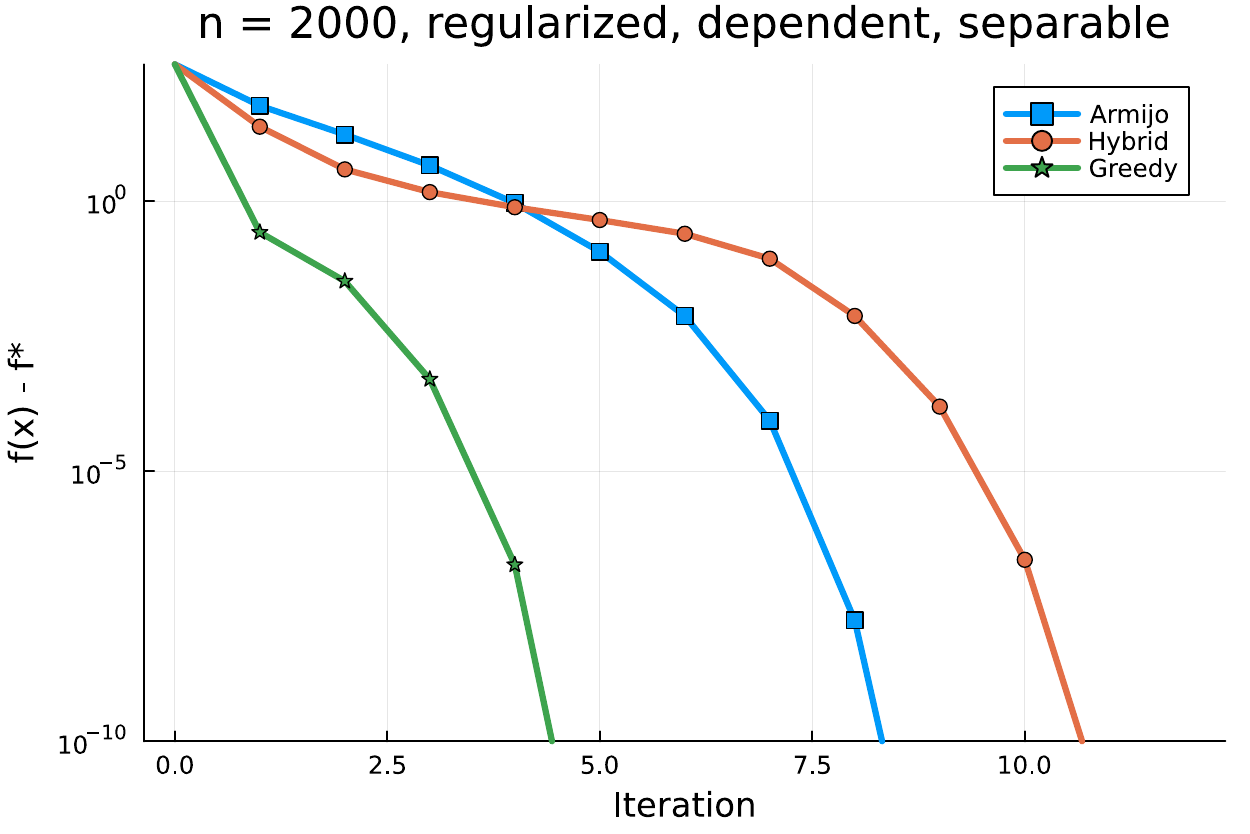}
	\includegraphics[width=.24\textwidth]{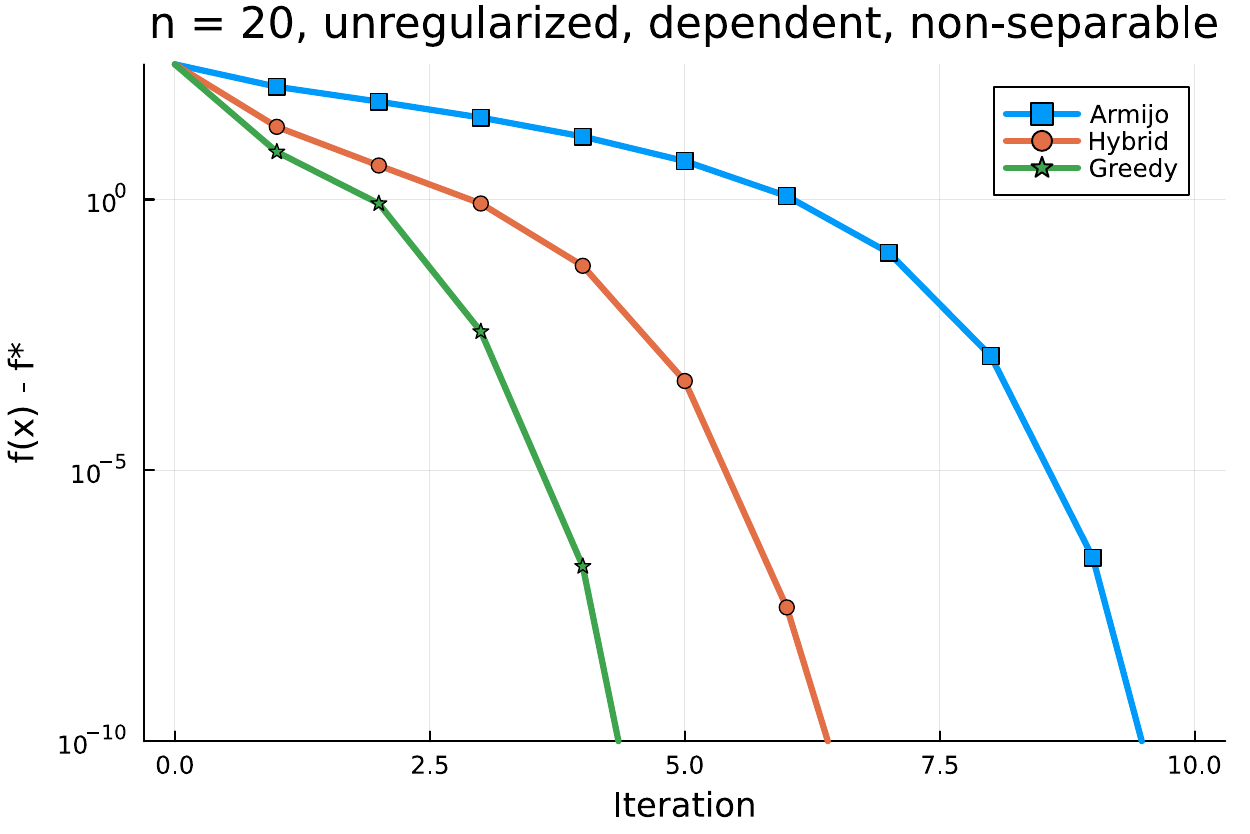}
	\includegraphics[width=.24\textwidth]{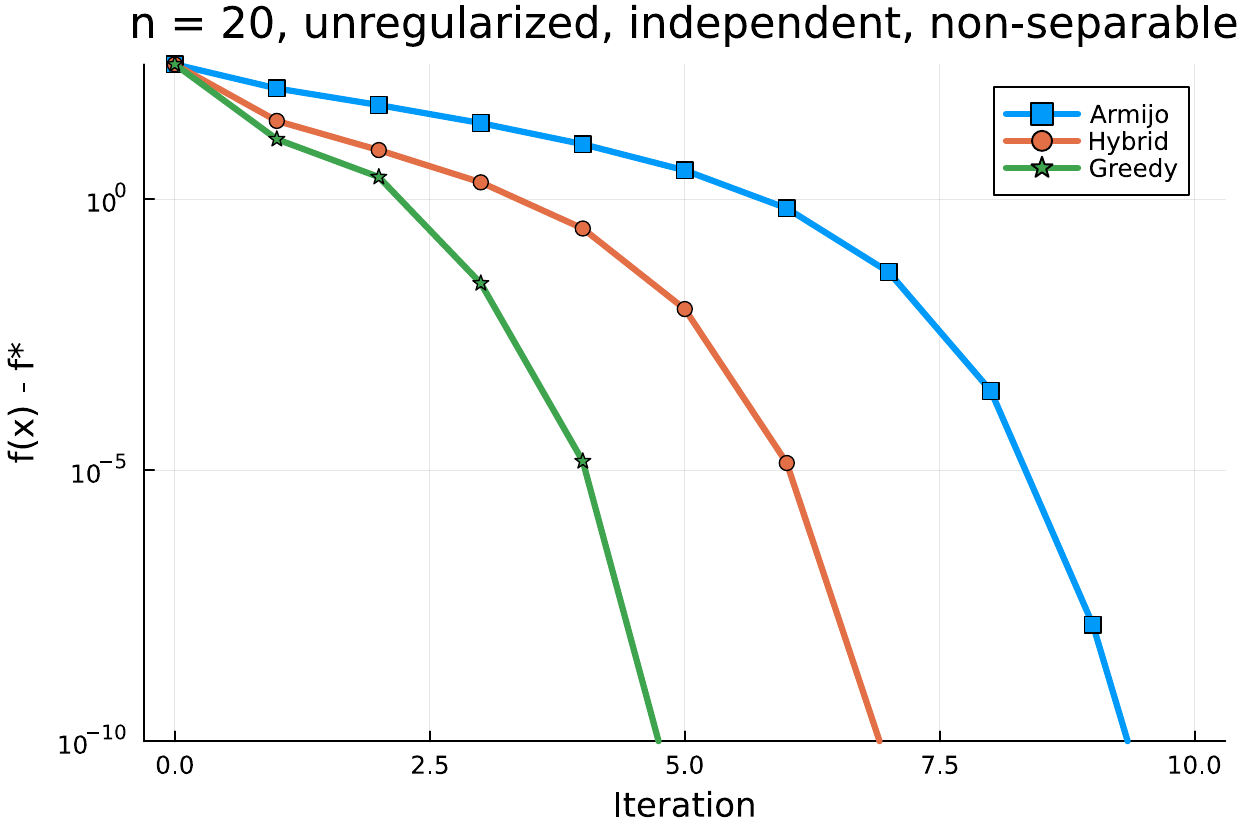}
	\includegraphics[width=.24\textwidth]{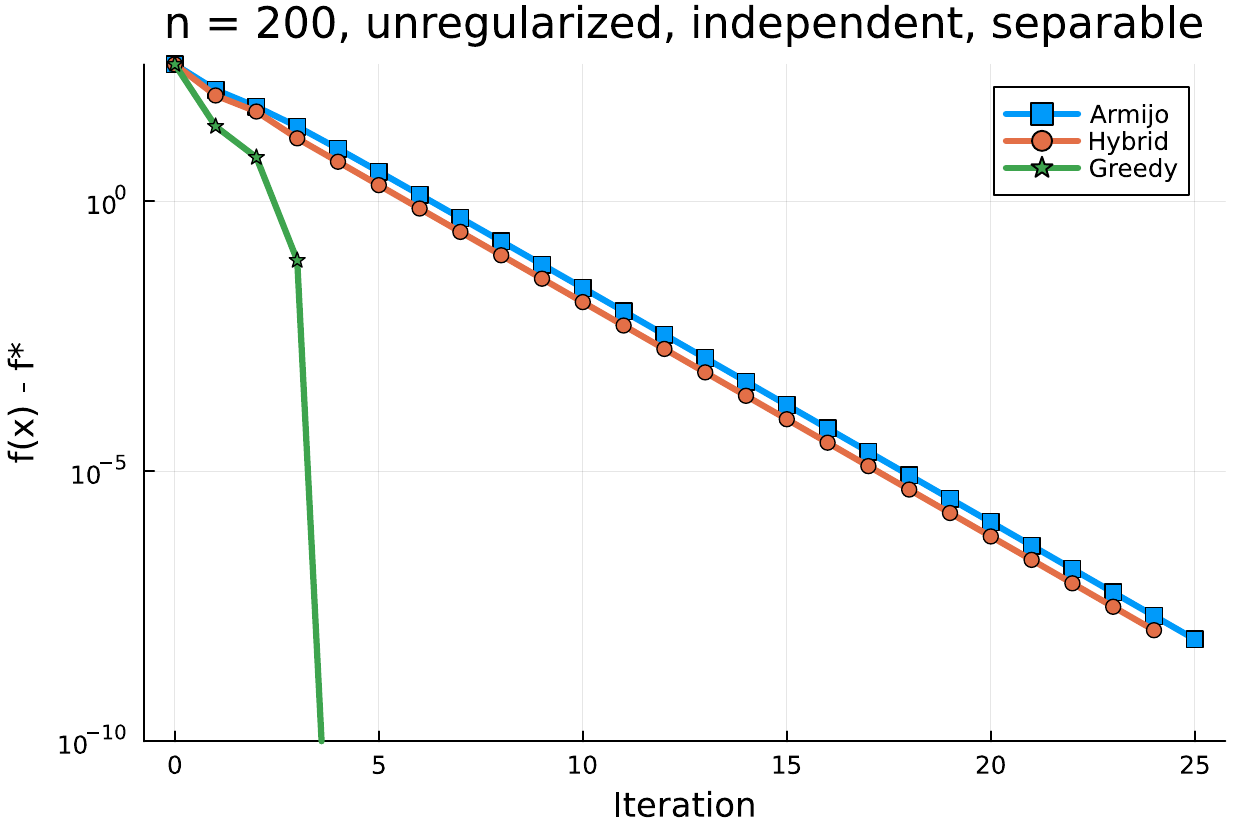}
	\includegraphics[width=.24\textwidth]{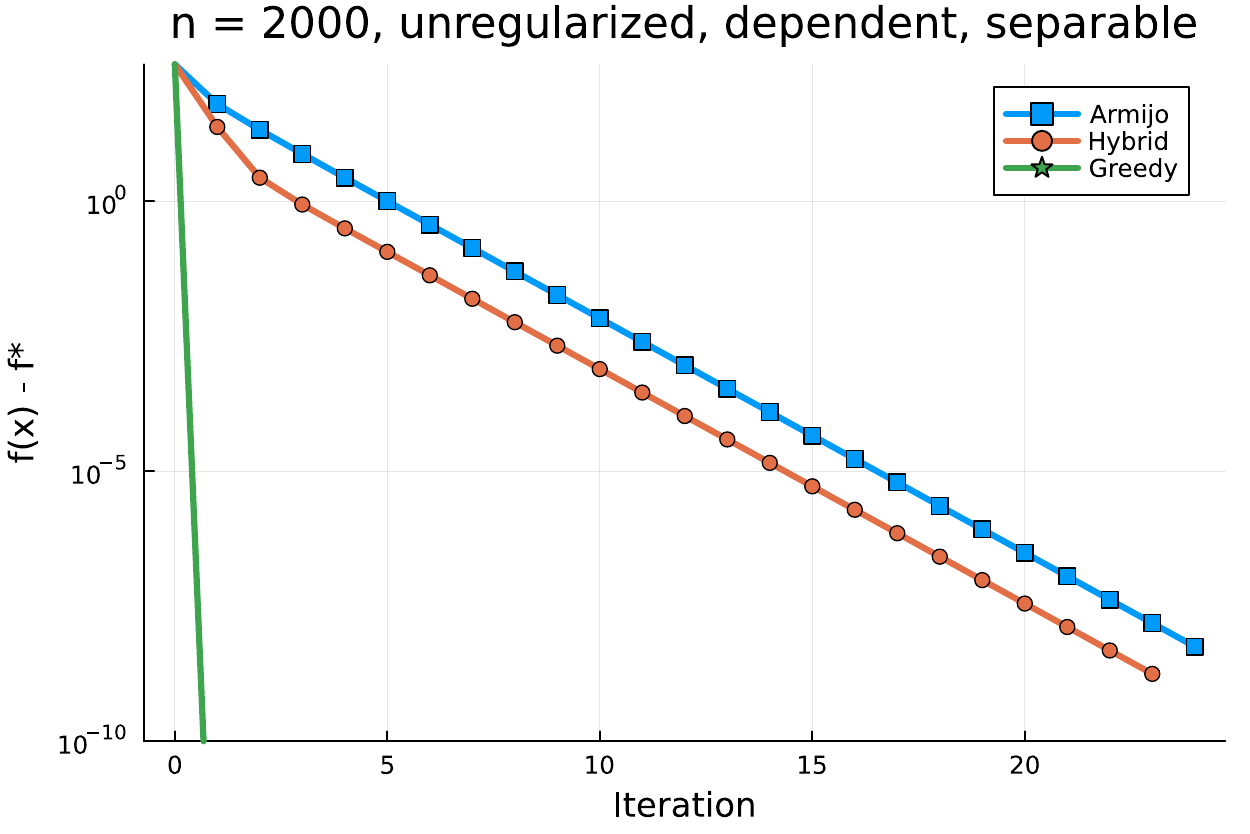}
	\caption{Comparison of Newton with Armijo backtracking, hybrid gradient-Newton, and greedy Newton on logistic regression problems with regularization (top row) and without regularization (bottom row). From left to right, in the bottom row  the problems are convex, strongly-convex, strictly-convex, and convex (all problems in the top row are strongly-convex). The right two datasets are linearly separable and the left two datasets are non-separable. 
	}
	\label{fig:logregf}
\end{figure}

In Figure~\ref{fig:logregf} we compare GN with Newton's method with Armijo backtracking from an initial step size of 1, and the hybrid of greedy gradient descent and pure Newton method discussed in Section 2.4. 
We see that GN outperformed the other two methods in all settings. In Appendix~\ref{app:exp} we report results on real data where the GN continued to perform the same or better than these other methods across 40 datasets.
In Appendix~\ref{app:time} we repeat the synthetic data experiment but measure performance using runtime. We found that the GN method continued to outperform the other methods using the crude runtime measure of performance.

\begin{figure}
	\includegraphics[width=.24\textwidth]{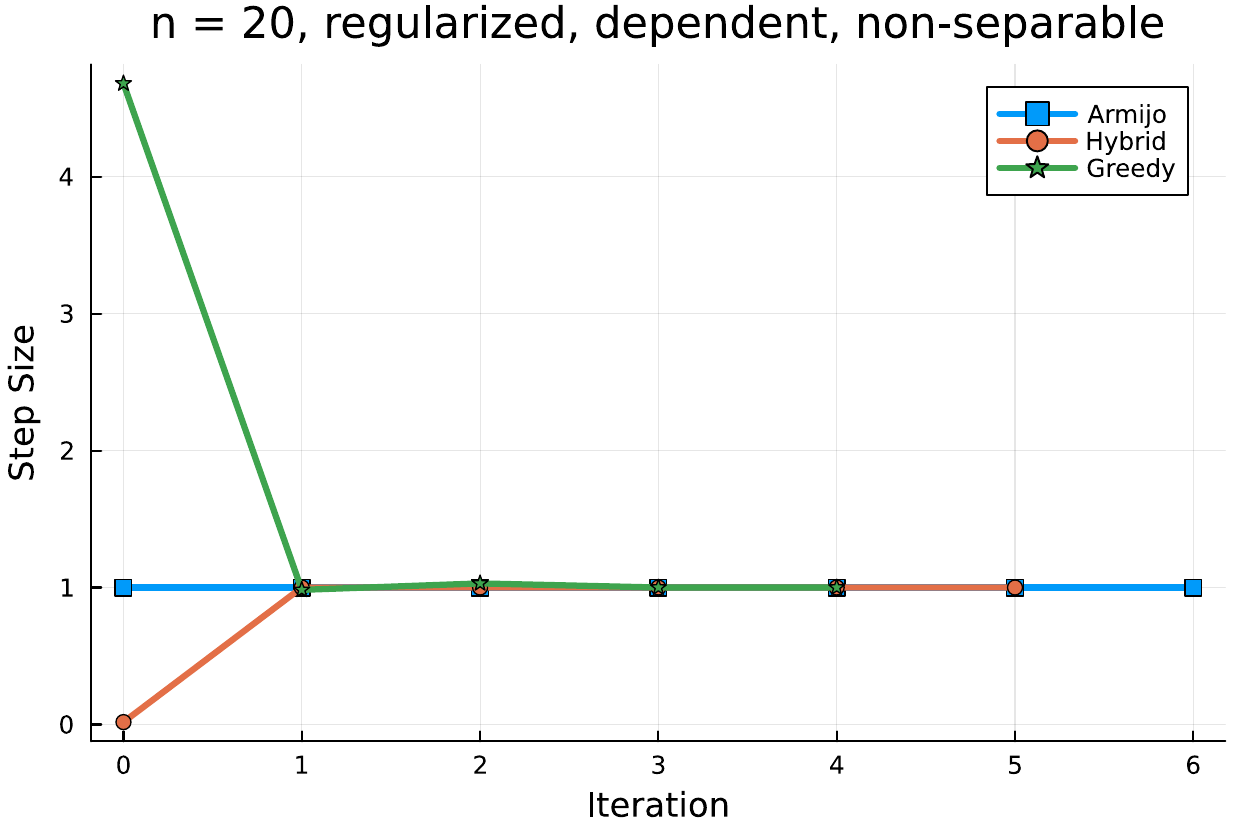}
	\includegraphics[width=.24\textwidth]{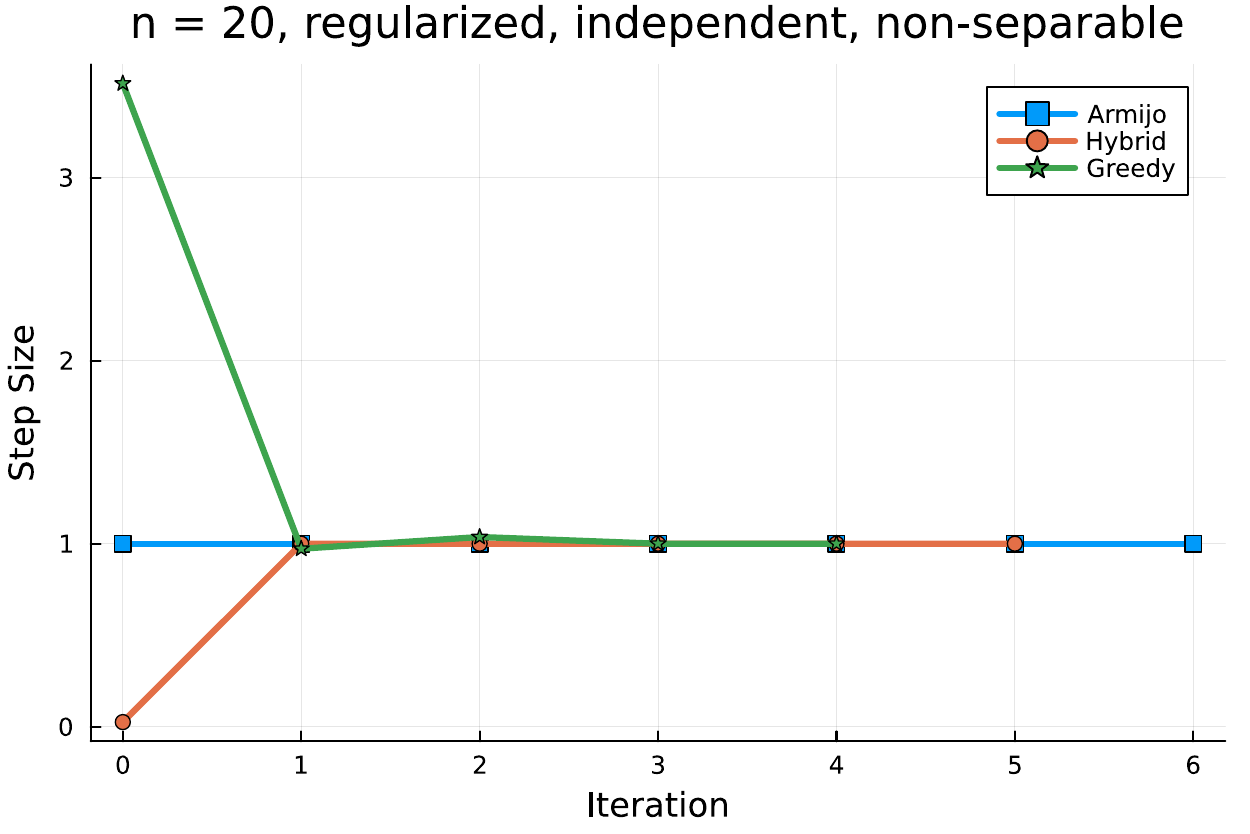}
	\includegraphics[width=.24\textwidth]{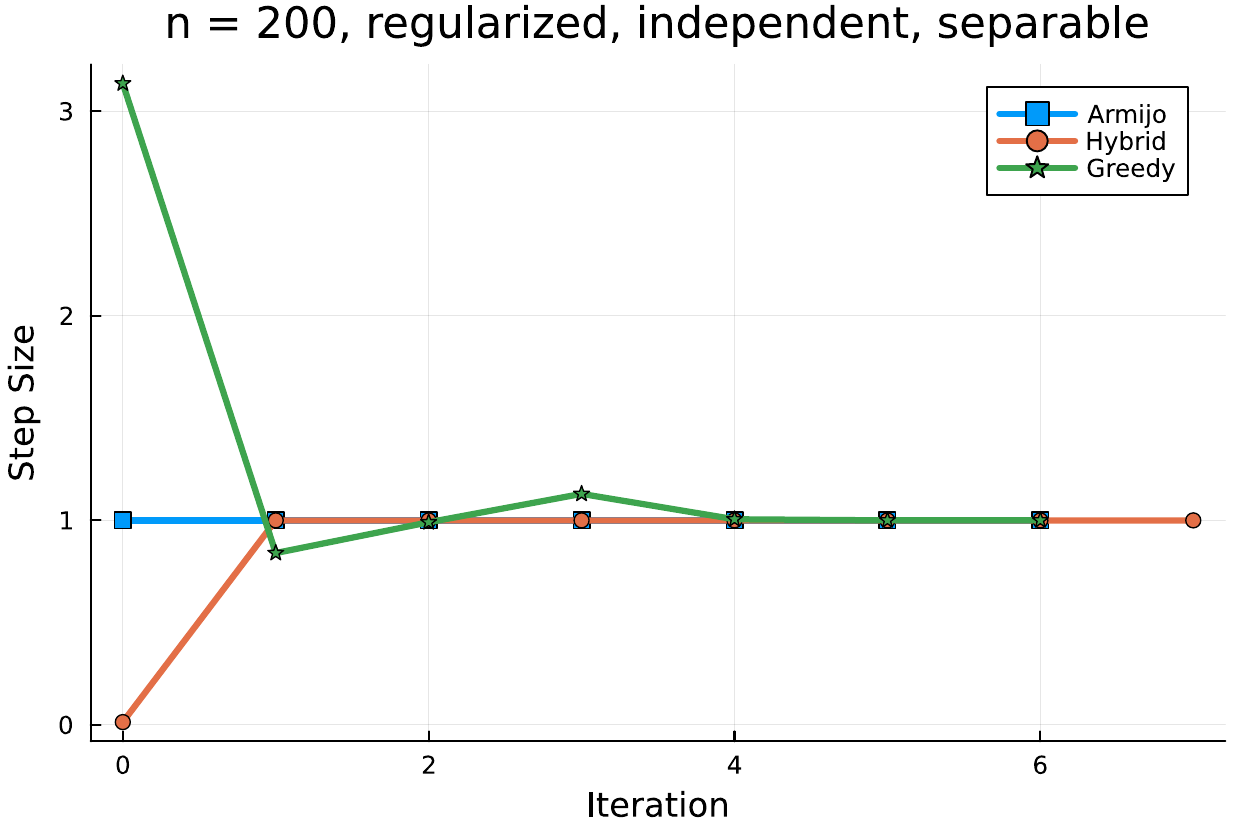}
	\includegraphics[width=.24\textwidth]{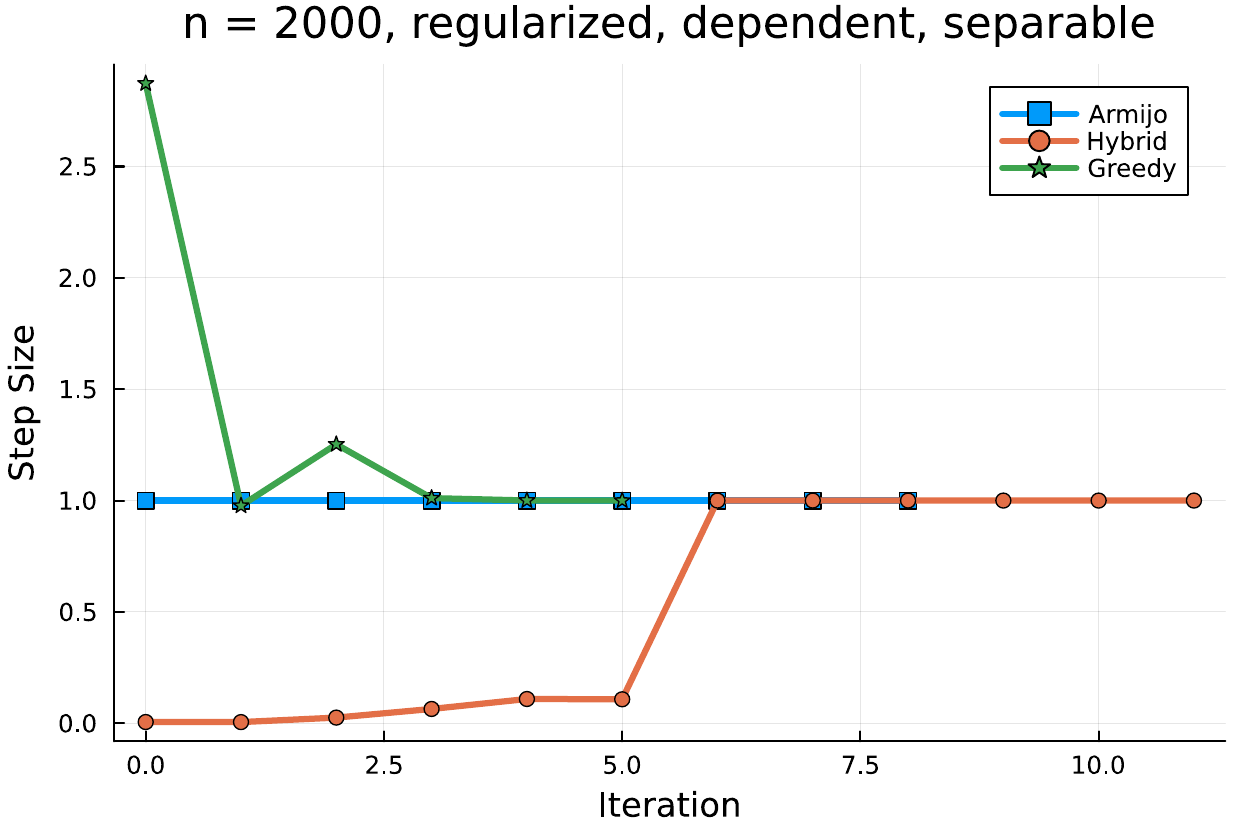}
	\includegraphics[width=.24\textwidth]{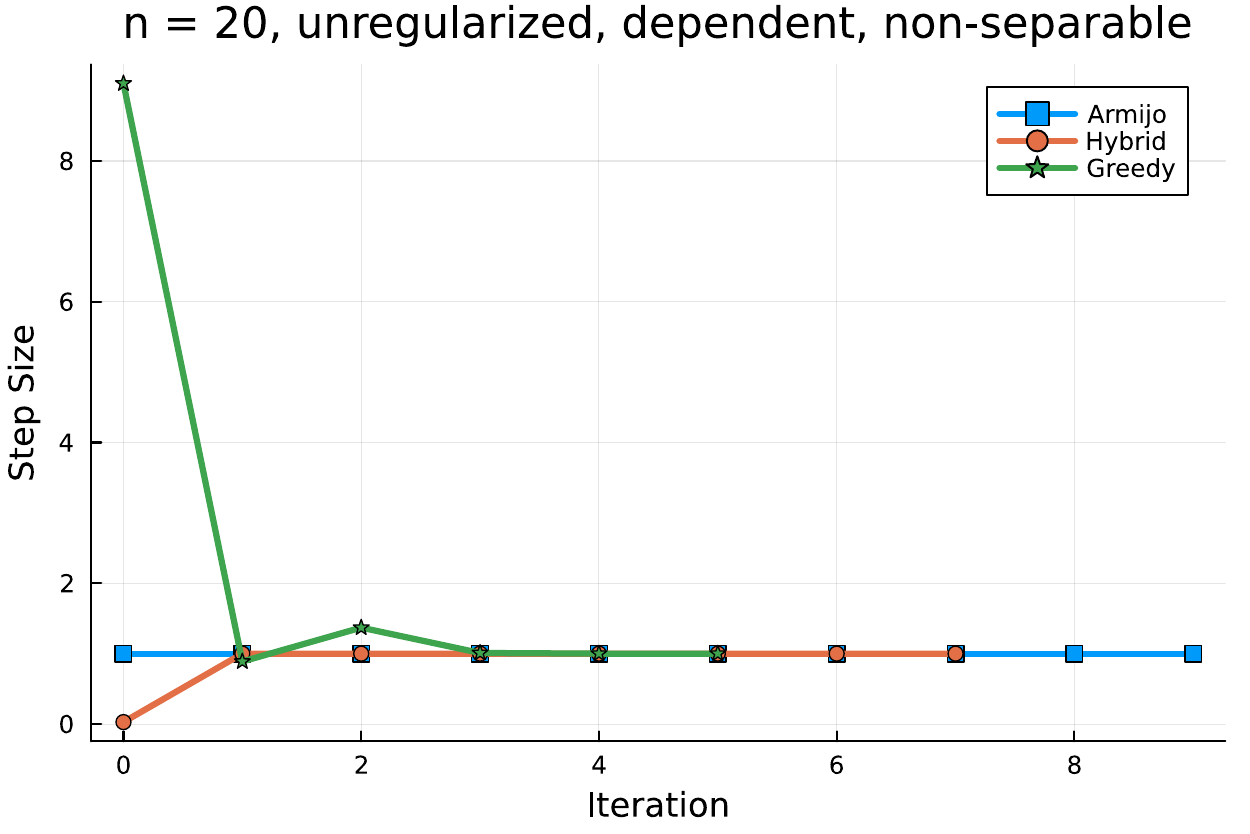}
	\includegraphics[width=.24\textwidth]{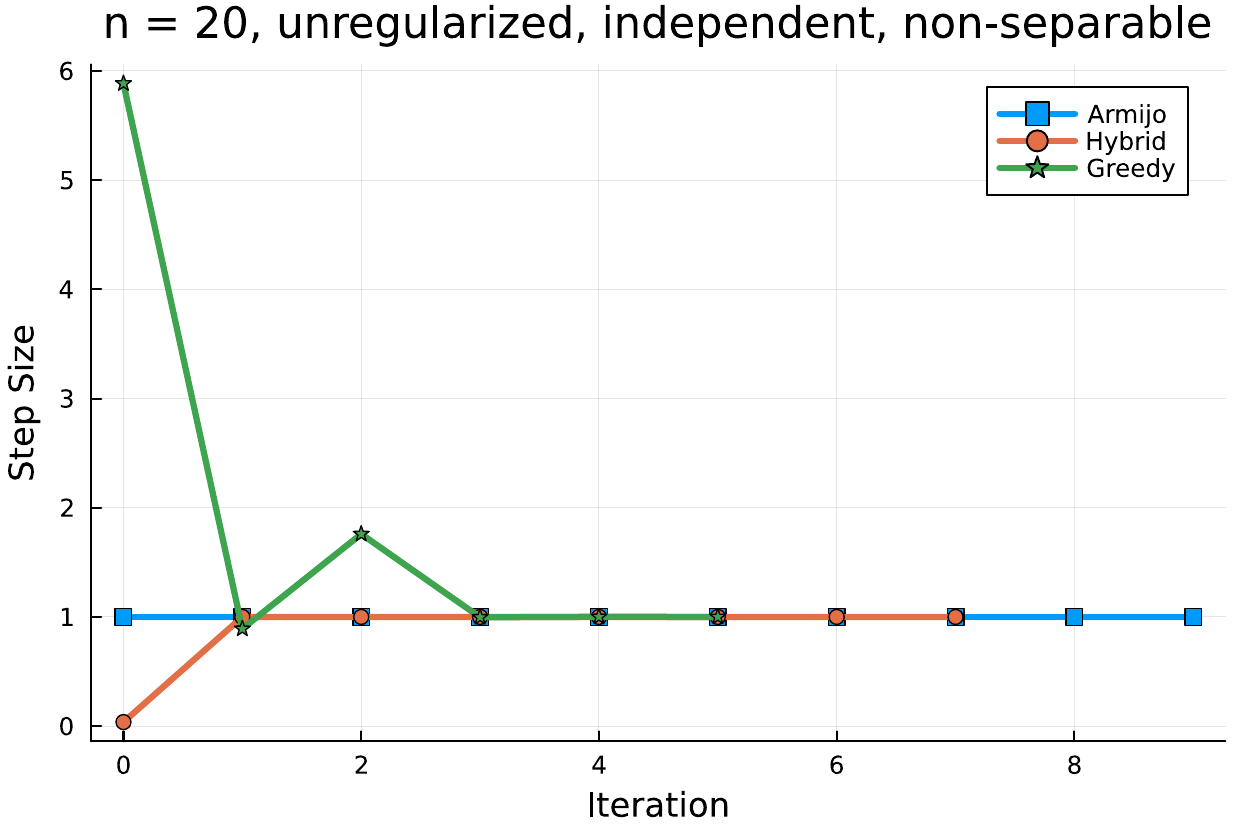}
	\includegraphics[width=.24\textwidth]{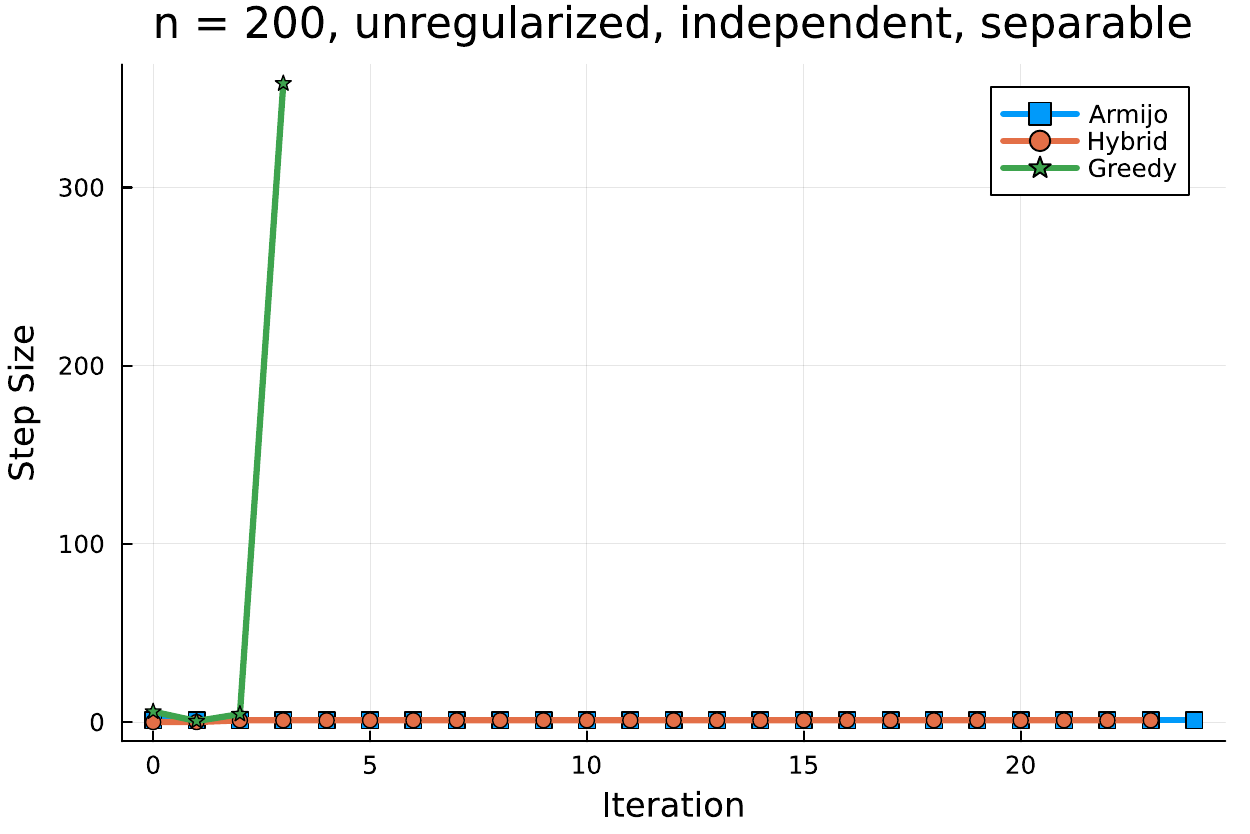}
	\includegraphics[width=.24\textwidth]{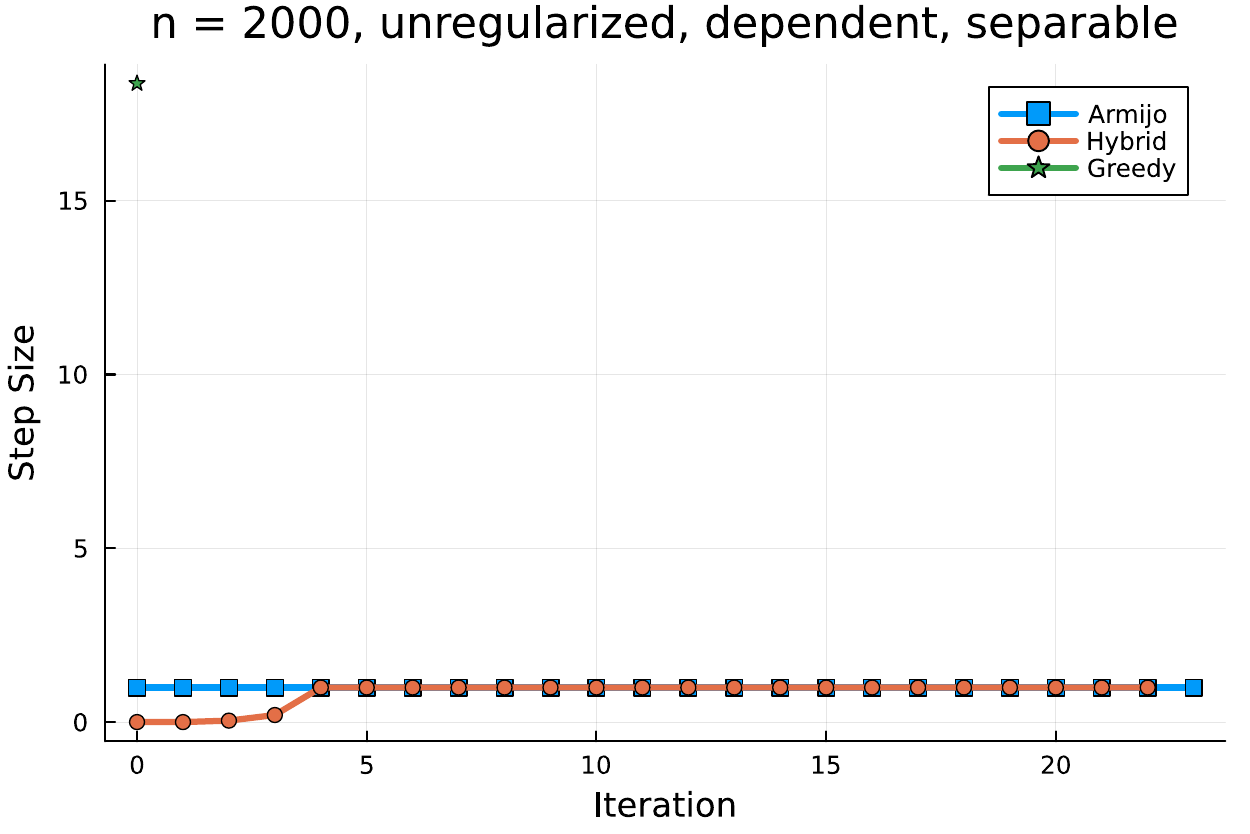}
	\caption{
		Step size chosen on each iteration for the plots in Figure~\ref{fig:logregf}. For iterations where the hybrid method took a gradient step, the gradient descent step size is shown (when a pure Newton step was taken the step size is 1).
	}
	\label{fig:logregt}
\end{figure}

In Figure~\ref{fig:logregt} we plot the step sizes used in the methods in Figure~\ref{fig:logregf}. For the hybrid method, a step size is 1 is reported when the pure Newton step is taken and otherwise the gradient descent step size is shown. Observe that GN used an initial step of at least 2 on every dataset, and the largest step size used was greater than 300.

The performance gain for GN was particularly large in the unregularized case for the two separable datasets (where the function is minimized by taking an $x$ that separates the data and making its magnitude arbitrarily large); in these cases the Armijo and hybrid methods performed poorly while GN achieved numerical accuracy extremely quickly (in 4 iterations and 1 iteration respectively). For these datasets GN used step sizes much larger than 1. 

In our experiments on non-separable datasets or with regularization, GN initially takes a large step then immediately takes steps that become close to 1. In contrast, 
 the Armijo backtracking method always accepted a step size of 1. Thus, GN's performance gain is largely due to using larger step sizes on early iterations. In light of this, one might consider modifying the Armijo line search to use a larger initialization. We explore this possibility in Appendix~\ref{app:armijo}. We found that in some cases initializing the Armijo backtracking procedure with a larger step could narrow or remove the initial gap between backtracking and GN. However, the Armijo backtracking procedure with a larger initialization often led to poor performance on later iterations (accepting large step sizes with worse performance than using a step size of 1).

\section{Open Problems}
\label{sec:GN_conclusion}

Our experiments show that we can use Newton's method more advantageously when we do not restrict the step size to be less than 1. However, our theory does not reflect the large performance increases we saw in practice. 
Below we list some open problems:
\begin{enumerate}
	\item Section 2.1: can we prove that step sizes bigger than 1 improve the global rate?
	\item Section 2.2: is the additional $\sqrt{L/\mu}$ term in the superlinear rate necessary?
	\item Section 2.3: can we analyze the rate at which $\alpha_k$ converges to 1?
	\item Section 2.4: can we justify why GN outperforms the theoretically-faster hybrid method?
	\item Section 3: can we prove a faster rate for GN on separable logistic regression problems?
\end{enumerate}
We close by discussing using an exact line search in other Newton-like algorithms. First, we note that a precise step size search could also be added to Newton's method with cubic regularization and that this does not change the radius of superlinear convergence of that method (Appendix~\ref{app:cubic}). For practical quasi-Newton and Hessian-free Newton methods that avoid computing the Hessian explicitly~\citep[see][]{Nocedal2006}, the cost of a precise step size search may significantly increase the iteration cost of the method. However, many problem structures arising in practice including logistic regression still allow an efficient line search even when function and gradient calculations are the bottleneck~\citep[see][]{Shea2023}.

\begin{appendices}

\section{Analysis of Greedy Newton}
\label{sec:app_analysis}

In this section, we prove the results in Sections~2.1-2.3. Our analyses are modifications of existing convergence analyses of the pure and backtracking Newton method~\citep{Boyd2004,Nesterov2018,Sun2021} to use an exact line search. 

\subsection{Global Convergence of Greedy Newton}

In this section we give the proof of Proposition~2.1.
Our assumption is that $\nabla^2 f(x_k)$ is positive definite with eigenvalues in $[\mu,L]$. This implies that $\nabla^2 f(x_k)^{-1}$ is symmetric and positive definite with eigenvalues in $[1/L,1/\mu]$. 
Using these facts in a Taylor expansion gives
{\scriptsize
	\begin{align*}
		f(x_{k+1}) & = f(x_k) + \nabla f(x_k)^T(x_{k+1}-x_k) + \frac 1 2(x_{k+1}-x_k)^T\nabla^2 f(z)(x_{k+1}-x_k) & \text{(for $z$ between $x_{k+1}$ and $x_k$)}\\
		&  \leq f(x_k) + \nabla f(x_k)^T(x_{k+1}-x_k) + \frac{L}{2}\norm{x_{k+1}-x_k}^2 &\text{($\nabla^2f(z) \preceq LI$)}\\
		& =  f(x_k)  - \alpha_k \nabla f(x_k)^T \nabla^2f(x_k)^{-1}\nabla f(x_k)+ \frac{L\alpha_k^2}{2}\norm{\nabla^2f(x_k)^{-1}\nabla f(x_k)}^2& \text{($x_{k+1}$ is Newton step~\eqref{eq:GN_update})}\\
		& =  f(x_k)  - \alpha_k \nabla f(x_k)^T \nabla^2f(x_k)^{-1}\nabla f(x_k)+ \frac{L\alpha_k^2}{2}\nabla f(x_k)^T\nabla^2 f(x_k)^{-2}\nabla f(x_k)\\
		& \leq f(x_k) - \alpha_k \nabla f(x_k)^T \nabla^2f(x_k)^{-1}\nabla f(x_k) + \frac{L\alpha_k^2}{2\mu} \nabla f(x_k)^T \nabla^2 f(x_k)^{-1}\nabla f(x_k) & \text{($\nabla^2 f(x_k)^{-1} \preceq (1/\mu)I)$}\\
		& = f(x_k) - \alpha_k\left(1 - \frac{\alpha_k}{2}\frac{L}{\mu}\right)\nabla f(x_k)^T\nabla^2f(x_k)^{-1}\nabla f(x_k)
	\end{align*}
}
With an exact line search~\eqref{eq:newton_exactLS}, we decrease the function by at least as much as choosing $\alpha_k=\mu/L$, so we have
\begin{align*}
	f(x_{k+1}) 
	& \leq f(x_k) - \frac{\mu}{2L}\nabla f(x_k)^T \nabla^2f(x_k)^{-1}\nabla f(x_k)\\
	& \leq f(x_k) - \frac{\mu}{2L^2}\norm{\nabla f(x_k)}^2,
\end{align*}
using that $\nabla^2 f(x_k)^{-1}\succeq (1/L)I$.
By strong convexity we have $\frac{1}{2\mu}\norm{\nabla f(x_k)}^2 \geq f(x_k) - f(x_*)$, using this and subtracting $f(x_*)$ from both sides gives
\begin{align*}
	f(x_{k+1}) - f^* & \leq f(x_k) - f(x_*) - \frac{\mu^2}{L^2}[f(x_k) - f(x_*)]\\
	& = \left(1-\frac{\mu^2}{L^2}\right)[f(x_k) - f^*].
\end{align*}
Applying this recursively gives the result.

\subsection{Local Convergence of ``As Fast as Newton" Methods}

In this section we give the proof of Proposition~2.2.
If $f$ is $\mu$-strongly convex with an $M$-Lipschitz Hessian, then the pure Newton update $x_{k+1}^N$~\eqref{eq:newton_update} satisfies~\cite{Bertsekas2016,Sun2021}
\begin{equation}
	\label{eq:easyGN_globalNewtonProgress}
	\norm{x_{k+1}^N-x_*}\leq\frac{M}{2\mu}\norm{x_k-x_*}^2
\end{equation}
If $\nabla f$ is $L$-Lipschitz then we have
\[f(x)\leq f(x_*)+\nabla f(x_*)^\intercal (x-x_*)+\frac{L}{2}\norm{x-x_*}^2,\]
which using $\nabla f(x_*)=0$ implies
\begin{equation}
	\label{eq:easyGN_Lsmooth}
	f(x)-f(x_*)\leq\frac{L}{2}\norm{x-x_*}^2
\end{equation}
Similarly, $\mu$-strongly convexity of $f$ implies that
\begin{equation}
	\label{eq:easyGN_strongConvexity}
	f(x)-f(x_*)\geq\frac{\mu}{2}\norm{x-x_*}^2
\end{equation}
Thus, if an algorithm has $f(x_{k+1}) \leq f(x_{k+1}^N)$ then we have
\begin{align*}
	\norm{x_{k+1} - x_*}^2 & \leq
	\frac{2}{\mu}[f(x_{k+1}) - f(x_*)]\\
	& \leq \frac{2}{\mu}[f(x_{k+1}^N) - f(x_*)]\\
	& \leq \frac{L}{\mu}\norm{x_{k+1}^N - x_*}^2.
\end{align*}
Combined with the progress of Newton's method given by Equation \eqref{eq:easyGN_globalNewtonProgress}, we get
\begin{equation}
	\label{eq:easyGN_bound}
	\norm{x_{k+1}-x_*}\leq\sqrt{\frac{L}{\mu}}\frac{M}{2\mu}\norm{x_k-x_*}^2
\end{equation}

\subsection{Local Convergence of Newton with Arbitrary Step Size}

In this section we give the proof of Proposition~2.3.
The result holds if $\alpha_k=0$ since $L/\mu \geq 1$. Thus, we focus on the case of $\alpha_k\neq 0$. Subtracting $x_*$ from both sides of the Newton update with a step size of $\alpha_k\neq 0$~\eqref{eq:GN_update} gives
\[x_{k+1}-x_* = x_k-\alpha_k\nabla^2f(x_k)^{-1}\nabla f(x_k)-x_*\]
Since $\nabla f(x_*)=0$, we can add the quantity $\alpha_k\nabla^2f(x_k)^{-1}\nabla f(x_*)$ on the right hand side and rearrange
\begin{align*}
	x_{k+1}-x_* &= (x_k-x_*)-\alpha_k\nabla^2f(x_k)^{-1}\left(\nabla f(x_k)-\nabla f(x_*)\right)\\
	&=\frac{\alpha_k}{\alpha_k}\nabla^2f(x_k)^{-1}\nabla^2f(x_k)(x_k-x_*)-\alpha_k\nabla^2f(x_k)^{-1}\left(\nabla f(x_k)-\nabla f(x_*)\right)
\end{align*}
to get 
\begin{equation}
	\label{eq:GNstep_optGap}
	x_{k+1}-x_*=\alpha_k\nabla^2f(x_k)^{-1}\left[\frac{1}{\alpha_k}\nabla^2f(x_k)(x_k-x_*)-\left(\nabla f(x_k)-\nabla f(x_*)\right)\right]
\end{equation}
From Taylor's theorem we have
\[\nabla f(x_*) = \nabla f(x_k+(x_*-x_k))=\nabla f(x_k)+\int_0^1\nabla^2f(x_k+t(x_*-x_k))(x_*-x_k)dt\]
\[\implies \nabla f(x_k)-\nabla f(x_*)=\int_0^1\nabla^2f(x_k+t(x_*-x_k))(x_k-x_*)dt\]
Substituting this into \eqref{eq:GNstep_optGap} gives
\begin{align*}
	x_{k+1}-x_*&=\alpha_k\nabla^2f(x_k)^{-1}\left[\frac{1}{\alpha_k}\nabla^2f(x_k)(x_k-x_*)-\int_0^1\nabla^2f(x_k+t(x_*-x_k))(x_k-x_*)dt\right]\\
	&=\alpha_k\nabla^2f(x_k)^{-1}\left[\frac{1}{\alpha_k}\nabla^2f(x_k)-\int_0^1\nabla^2f(x_k+t(x_*-x_k))dt\right](x_k-x_*)\\
	&=\alpha_k\nabla^2f(x_k)^{-1}\left[\int_0^1\left(\frac{1}{\alpha_k}\nabla^2f(x_k)-\nabla^2f(x_k+t(x_*-x_k))\right)dt\right](x_k-x_*).
\end{align*}
Taking norms on both sides and using the Cauchy-Schwartz inequality gives
\begin{align}
	\norm{x_{k+1}-x_*}& \leq|\alpha_k|\norm{\nabla^2f(x_k)^{-1}}\left\|\int_0^1\frac{1}{\alpha_k}\nabla^2f(x_k)-\nabla^2f(x_k+t(x_*-x_k))dt\right\|\norm{x_k-x_*}\label{eq:hessDiffBound}
\end{align}
We bound the factor containing the integral using the triangle inequality, $\norm{\int_a^bg(t)dt}\leq\int_a^b\norm{g(t)}dt$,
\begin{align*}
	\left\|\int_0^1\frac{1}{\alpha_k}\nabla^2f(x_k)-\nabla^2f(x_k+t(x_*-x_k))dt\right\|&\leq\int_0^1\left\|\frac{1}{\alpha_k}\nabla^2f(x_k)-\nabla^2f(x_k+t(x_*-x_k))\right\|dt
\end{align*}
Using Lipschitz continuity of the gradient~\eqref{eq:Lmu} and Hessian~\eqref{eq:M}, this could be rewritten
\begin{align*}
	&\int_0^1\left\|-\nabla^2f(x_k)+\frac{1}{\alpha_k}\nabla^2f(x_k)+\nabla^2f(x_k)-\nabla^2f(x_k+t(x_*-x_k))\right\|dt\\
	&=\int_0^1\left\|\left(-1+\frac{1}{\alpha_k}\right)\nabla^2f(x_k)+\left(\nabla^2f(x_k)-\nabla^2f(x_k+t(x_*-x_k)\right)\right\|dt\\
	&\leq\int_0^1\left|1-\frac{1}{\alpha_k}\right|\left\|\nabla^2f(x_k)\right\|+\left\|\nabla^2f(x_k)-\nabla^2f(x_k+t(x_*-x_k)\right\|dt\\
	&\leq\int_0^1\left|1-\frac{1}{\alpha_k}\right|\left\|\nabla^2f(x_k)\right\|+Mt\norm{x_k-x_*}dt\\
	&\leq\int_0^1\left|1-\frac{1}{\alpha_k}\right|L+Mt\norm{x_k-x_*}dt,
\end{align*}
Substituting this back into~\eqref{eq:hessDiffBound}  gives
\begin{align*}
	\norm{x_{k+1}-x_*}&\leq|\alpha_k|\norm{\nabla^2f(x_k)^{-1}}\int_0^1\left[\left|1-\frac{1}{\alpha_k}\right|L+Mt\norm{x_k-x_*}\right]dt\norm{x_k-x_*}\\
	&=|\alpha_k|\norm{\nabla^2f(x_k)^{-1}}\left[\left|1-\frac{1}{\alpha_k}\right|L+\frac{M}{2}\norm{x_k-x_*}\right]\norm{x_k-x_*}\\
	&\leq\frac{|\alpha_k|}{\mu}\left[\left|1-\frac{1}{\alpha_k}\right|L+\frac{M}{2}\norm{x_k-x_*}\right]\norm{x_k-x_*}  & \text{(by strong convexity)}\\
	& = |\alpha_k-1|\frac{L}{\mu}\norm{x_k-x_*}+|\alpha_k|\frac{M}{2\mu}\norm{x_k-x_*}^2.
\end{align*}

\section{Greedy Newton Implementation}
\label{app:code}

\begin{figure}
\includegraphics[width=1\textwidth]{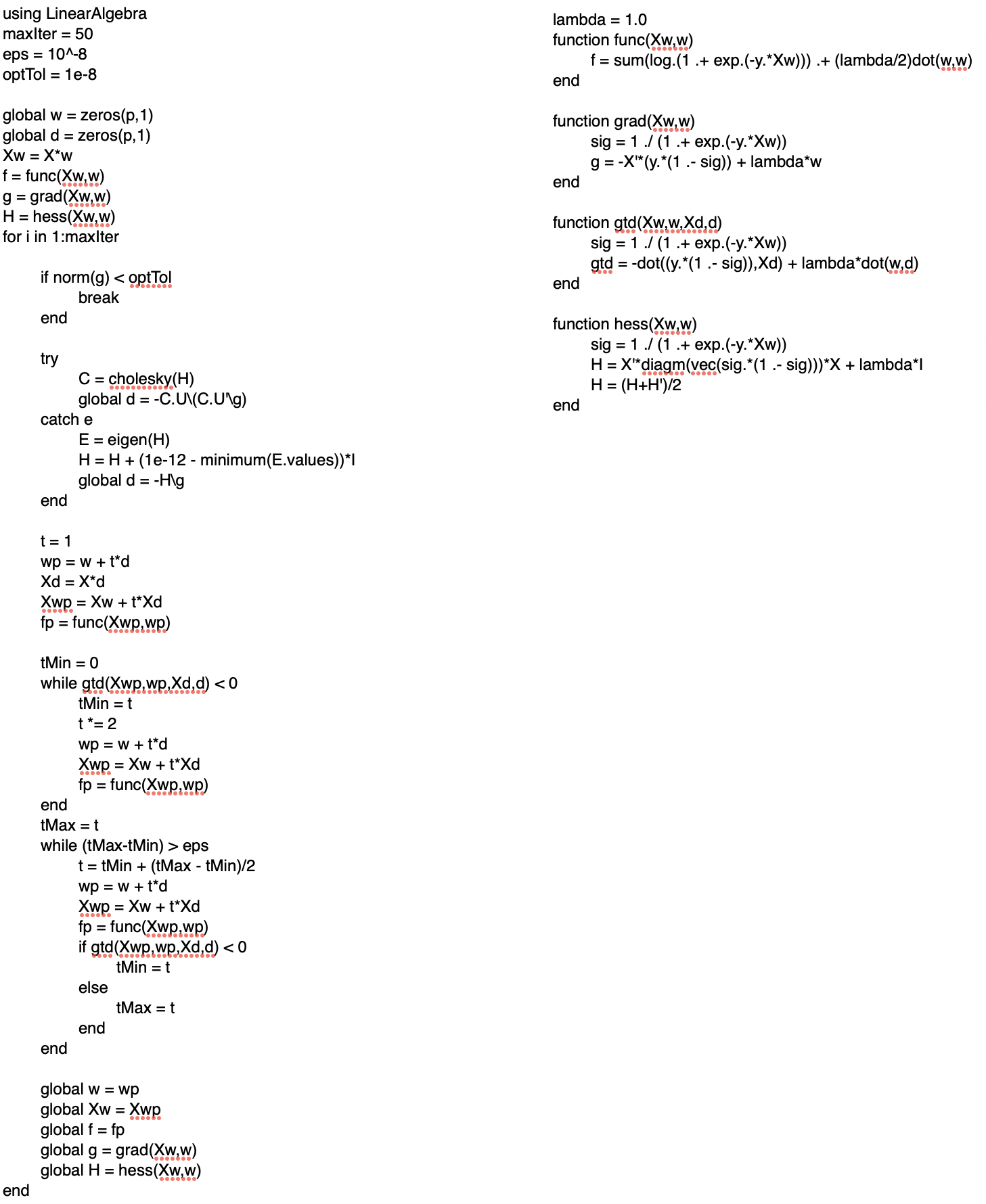}
\caption{On the left is the Julia code used in our experiments implemented the greedy Newton method. On the right is our code for regularized logistic regression implementing the function value, gradient, directional derivative, and Hessian. The unregularized version corresponds to setting lambda to 0 (or equivalently removing the terms involving lambda). Note that the bottlenecks on each iteration of the GN code are the matrix operations involving the Hessian in the try-catch statement, computing the matrix-vector product Xd, and finally calling the gradient and Hessian functions at the end. All other operations are performed on vectors or scalars.}
\label{fig:code}
\end{figure}

In Figure~\ref{fig:code} we give our Julia code implementation of the GN method, as well as the functions required to apply the code to logistic regression.

\section{Real Data Experiments}
\label{app:exp}

\begin{figure}
	\includegraphics[width=.24\textwidth]{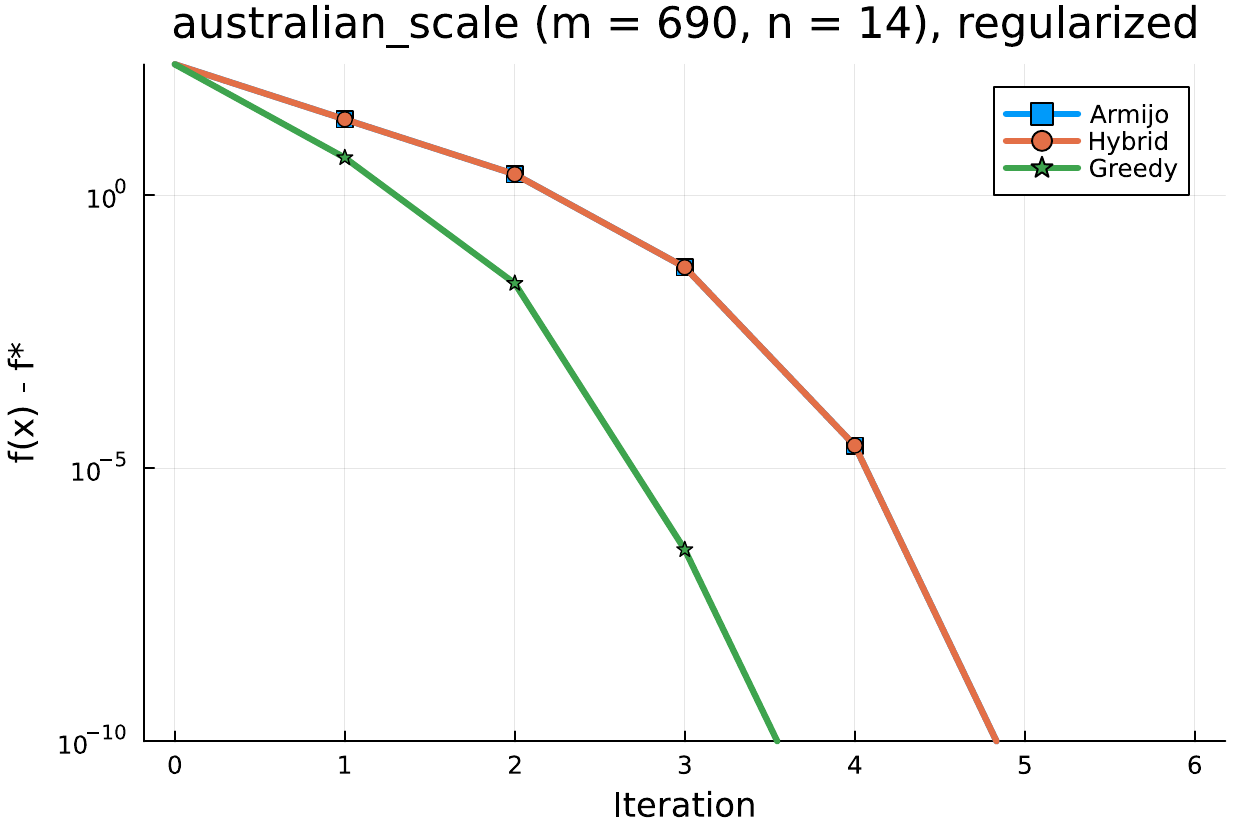}
	\includegraphics[width=.24\textwidth]{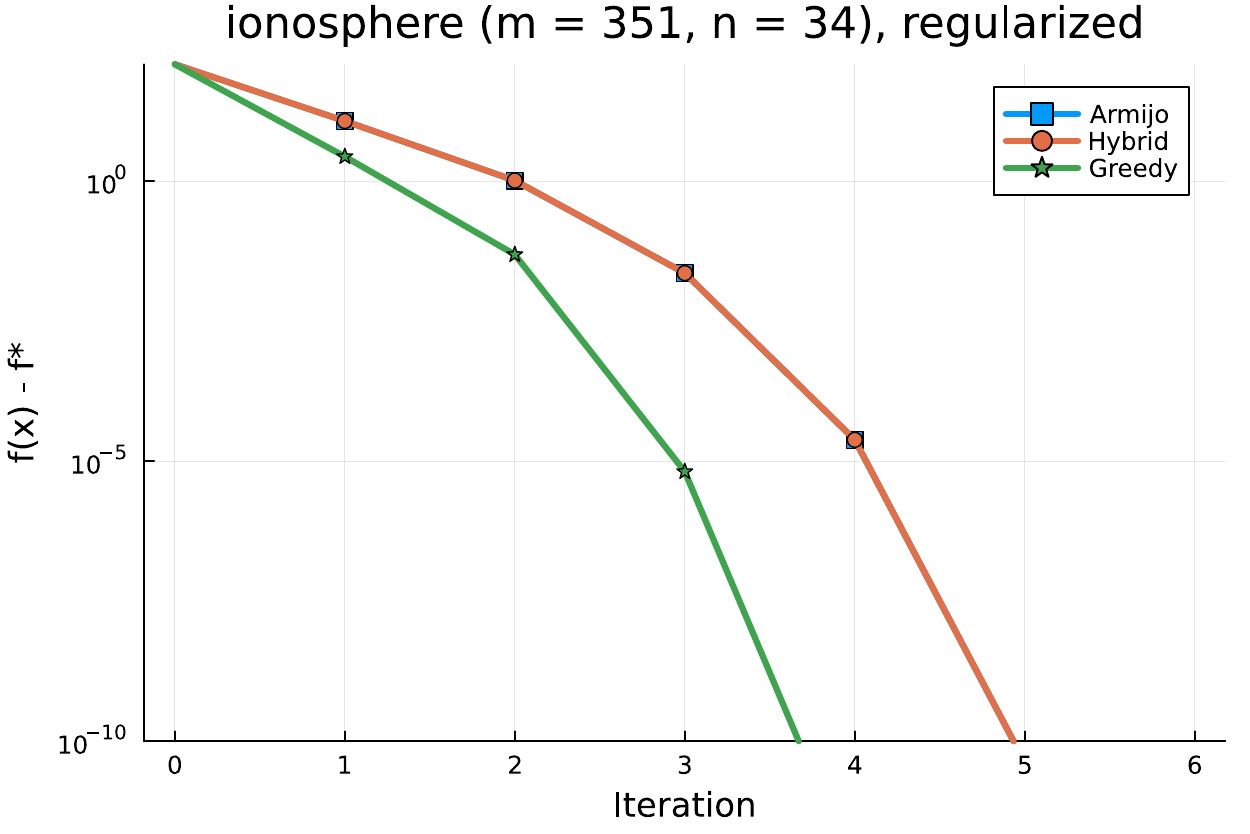}
	\includegraphics[width=.24\textwidth]{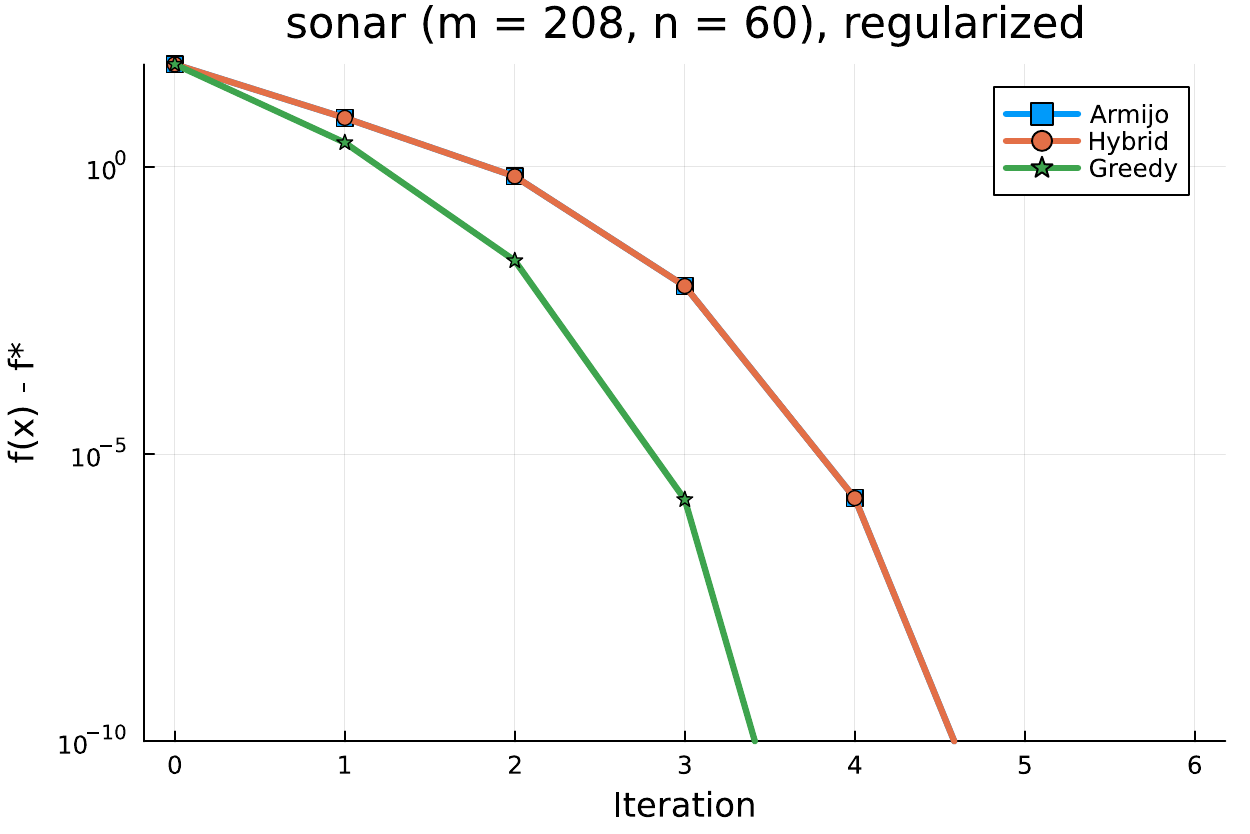}
	\includegraphics[width=.24\textwidth]{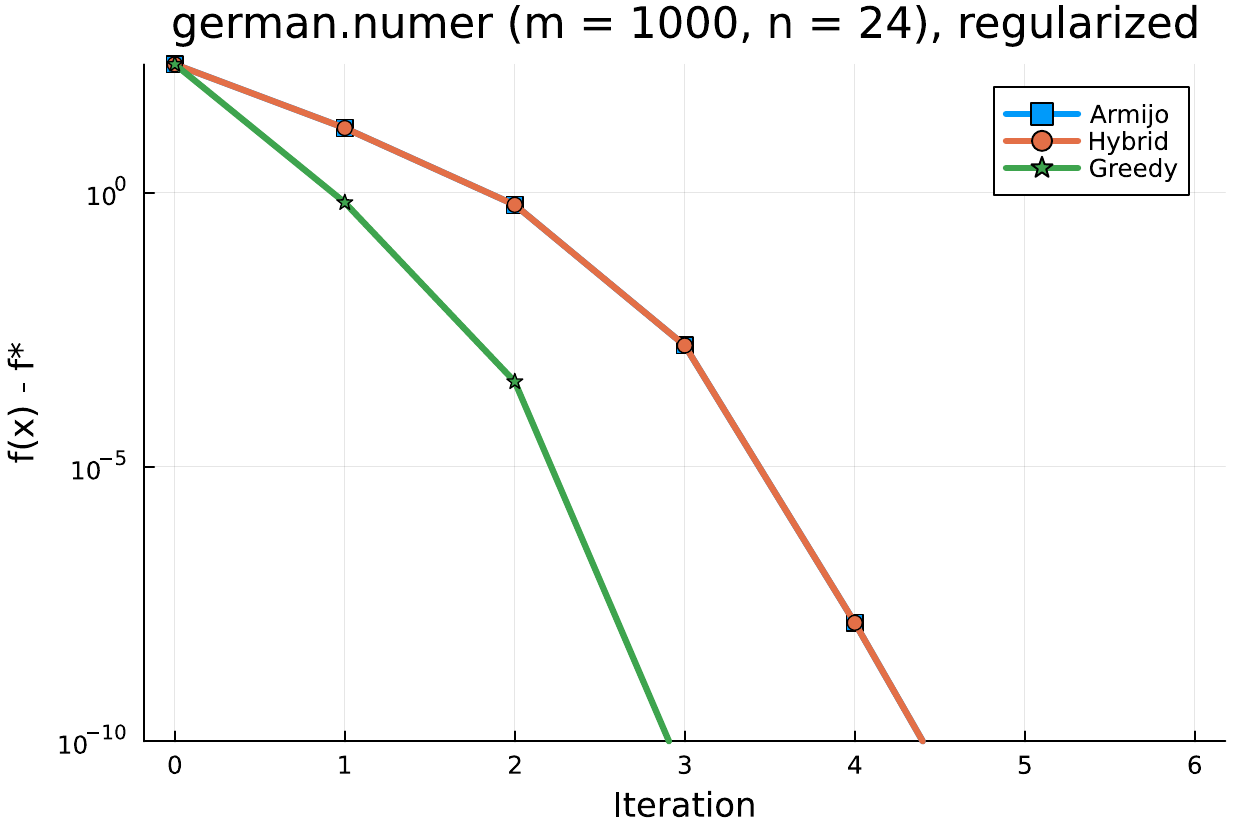}
	\includegraphics[width=.24\textwidth]{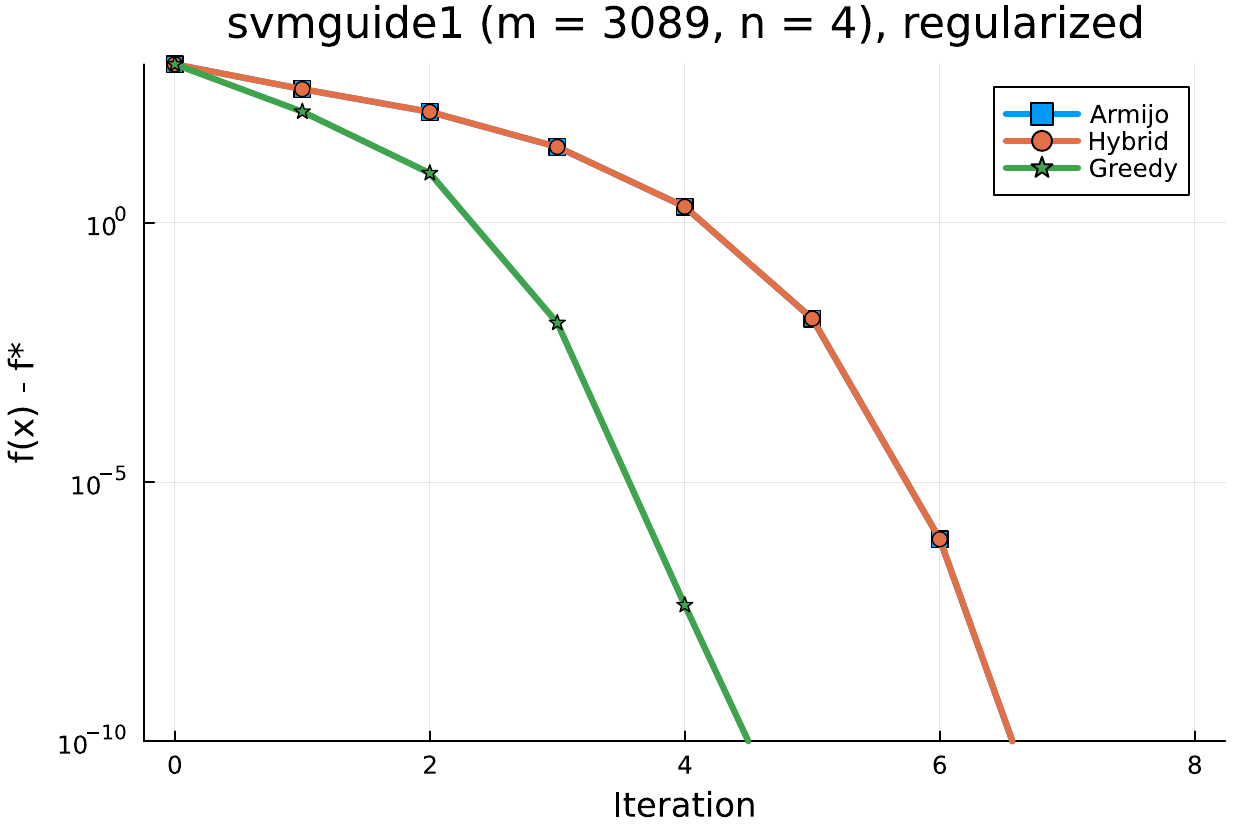}
	\includegraphics[width=.24\textwidth]{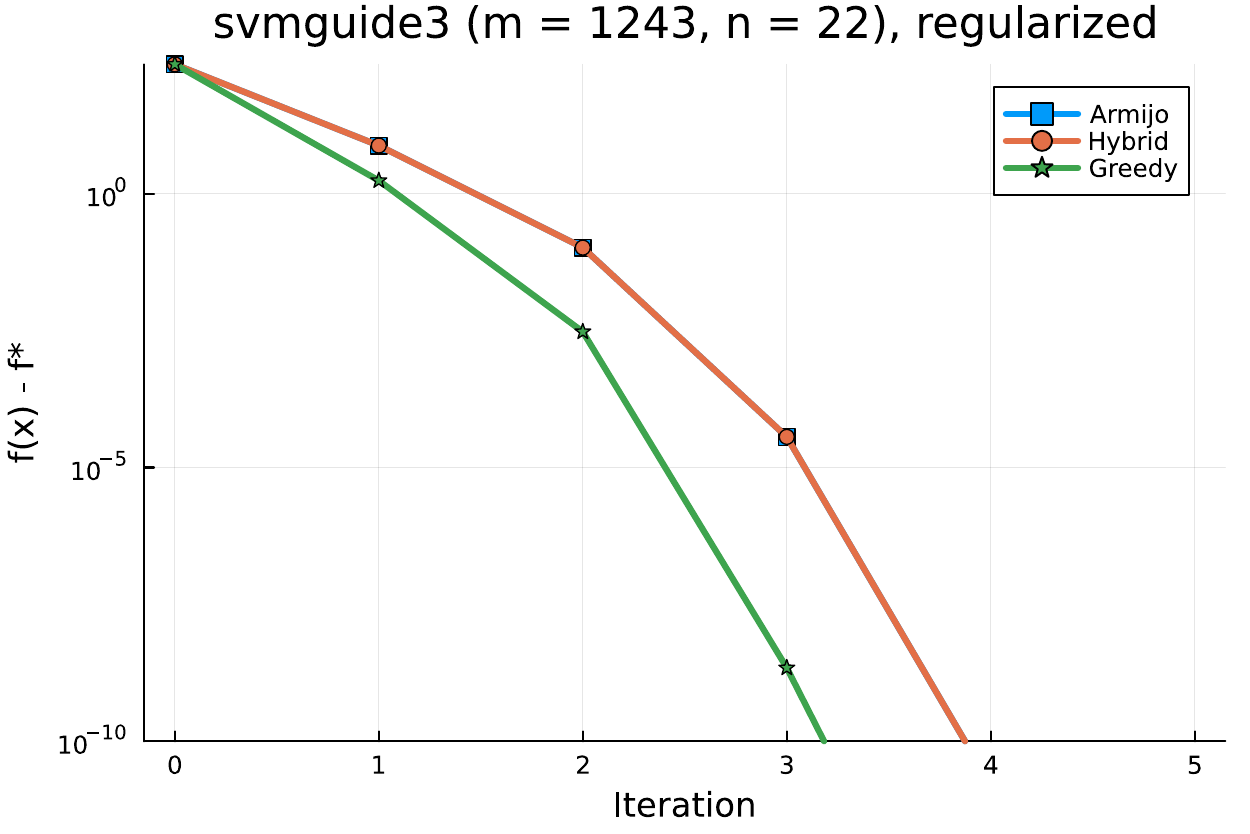}
	\includegraphics[width=.24\textwidth]{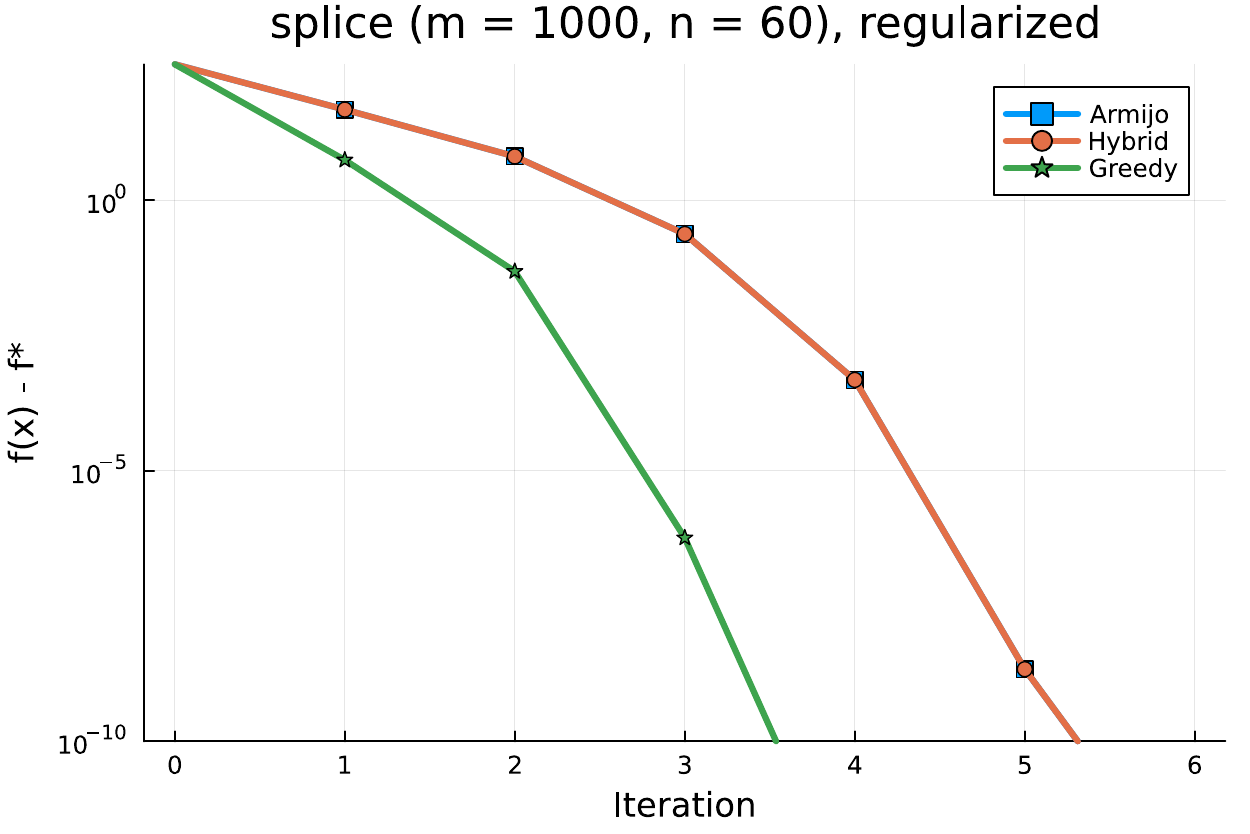}
	\includegraphics[width=.24\textwidth]{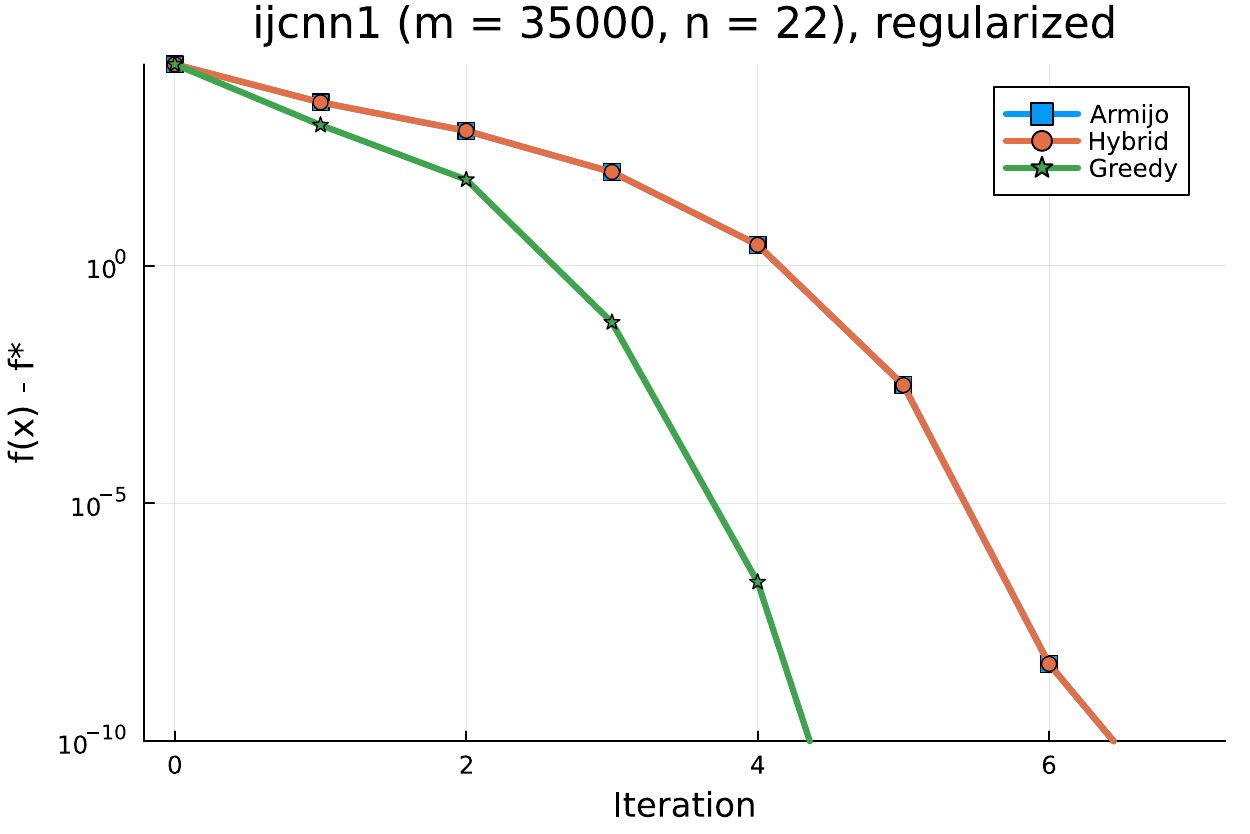}
	\includegraphics[width=.24\textwidth]{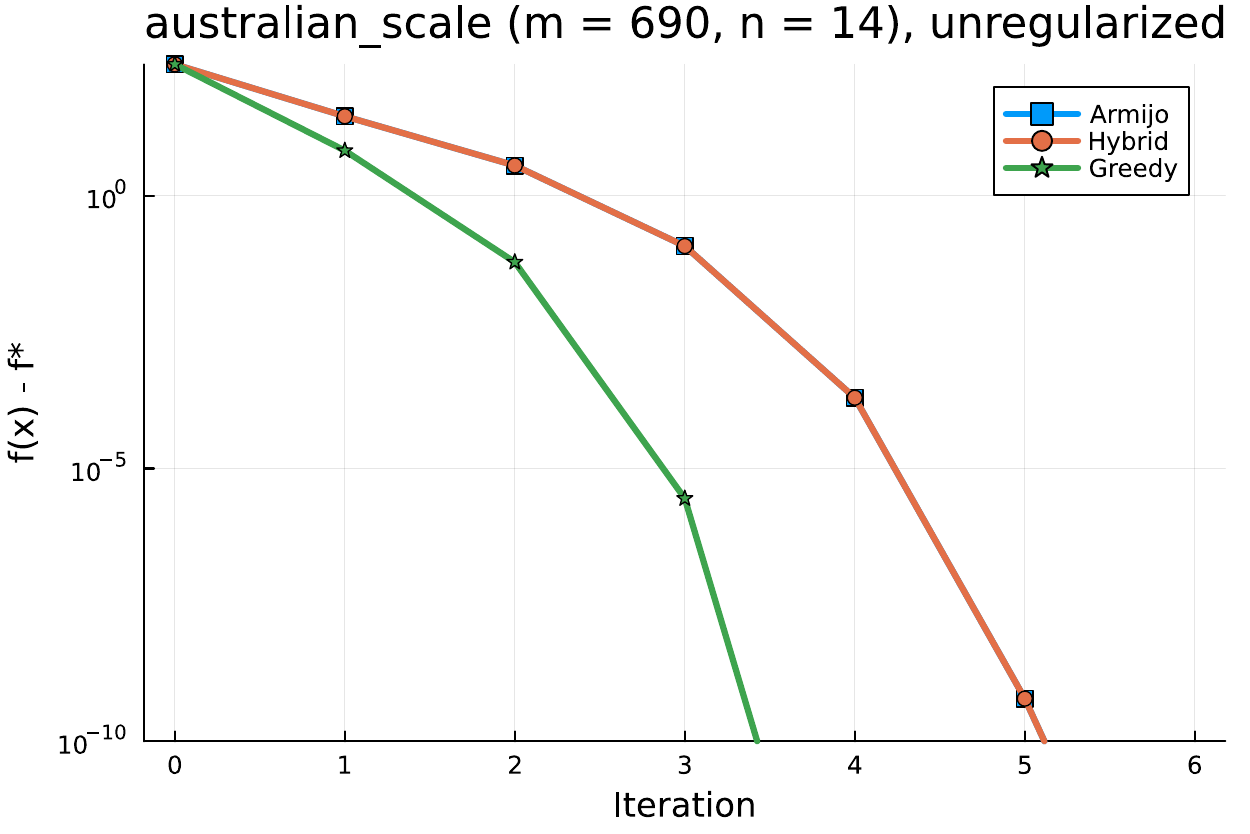}
	\includegraphics[width=.24\textwidth]{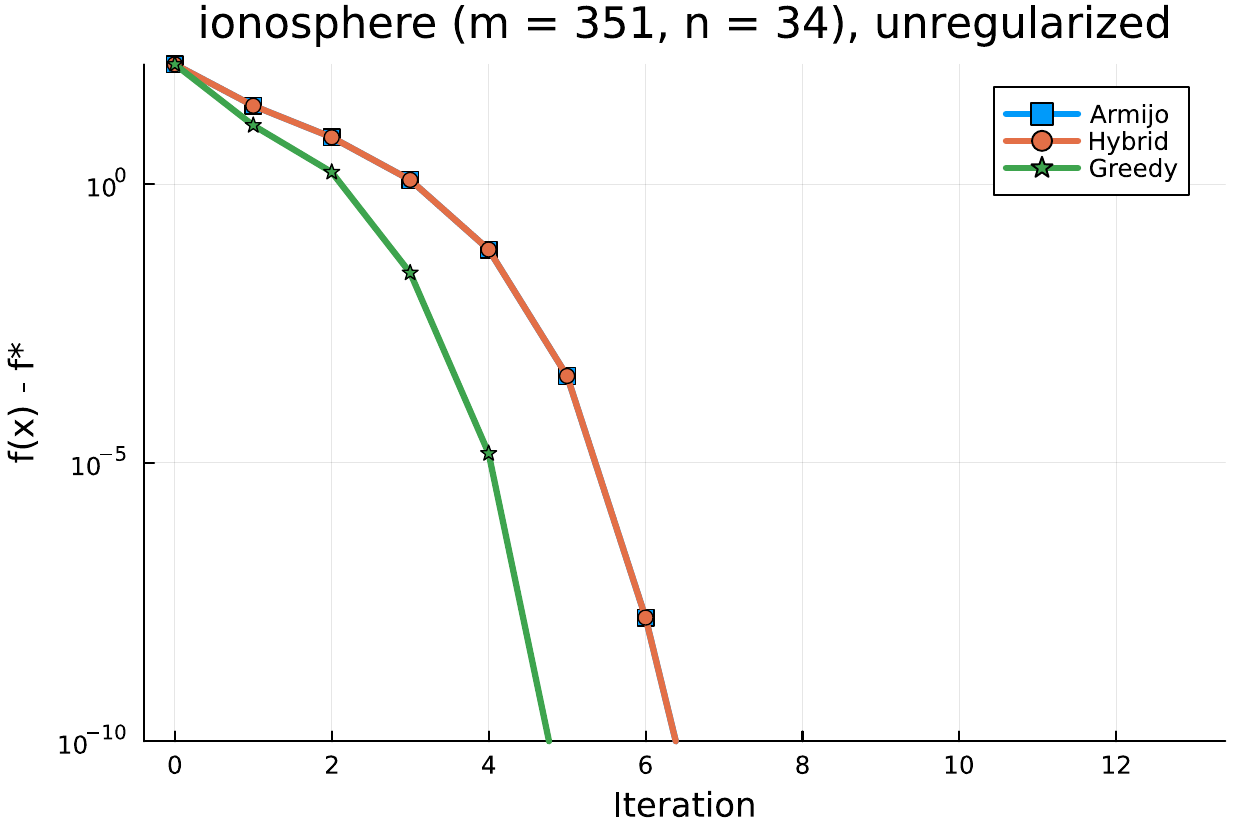}
	\includegraphics[width=.24\textwidth]{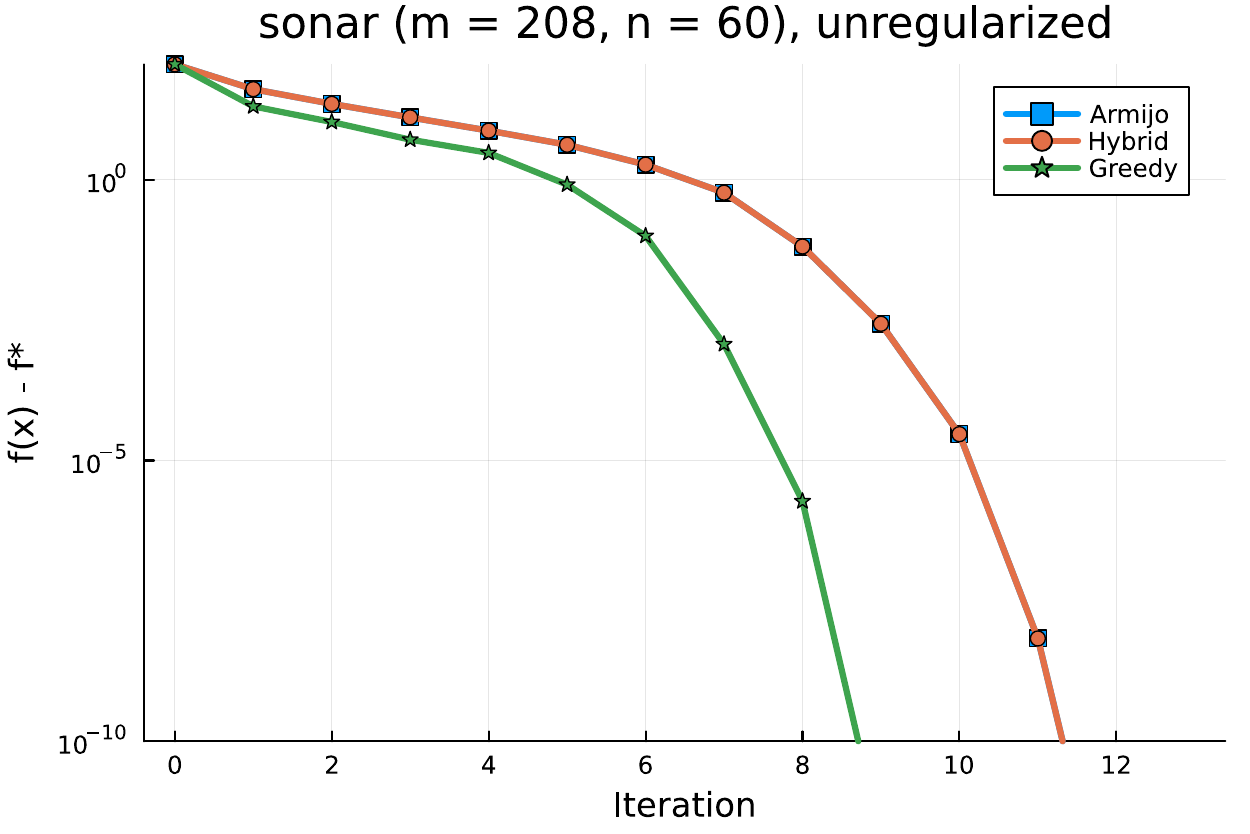}
	\includegraphics[width=.24\textwidth]{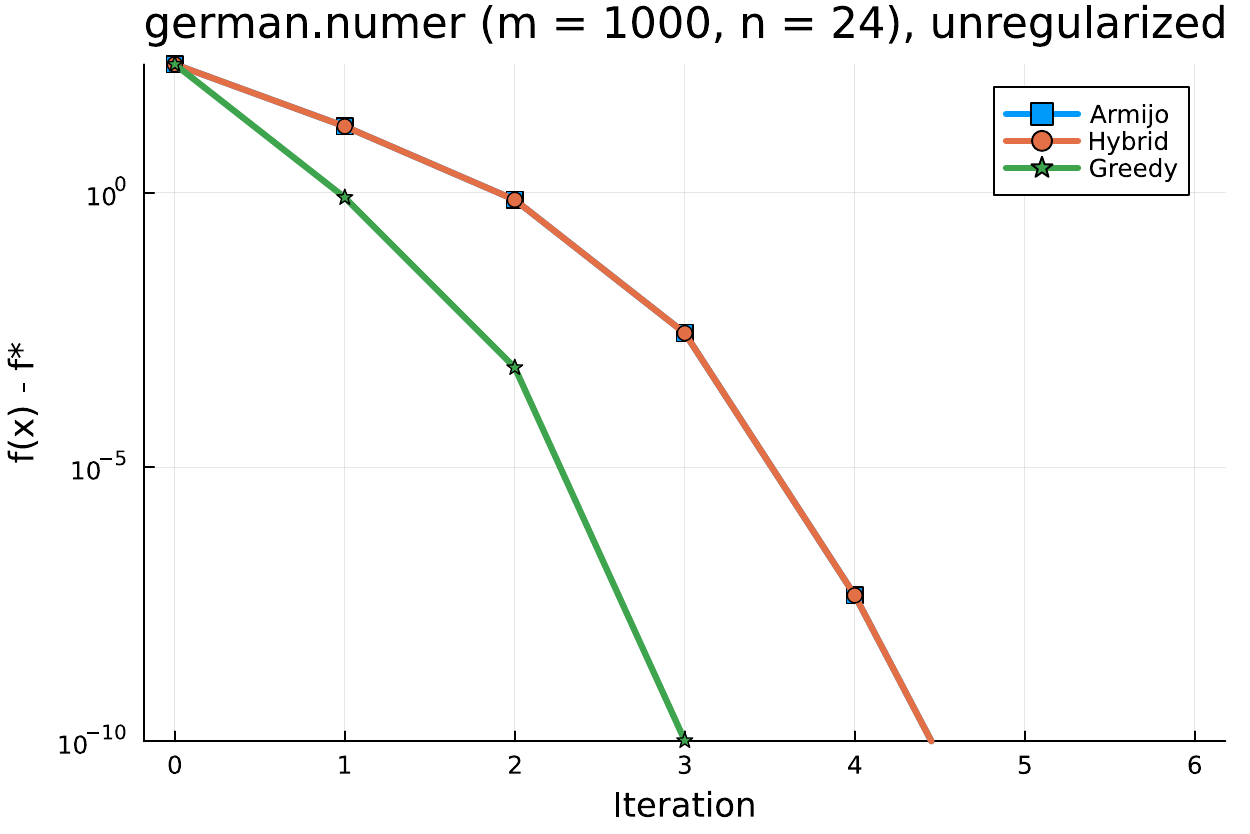}
	\includegraphics[width=.24\textwidth]{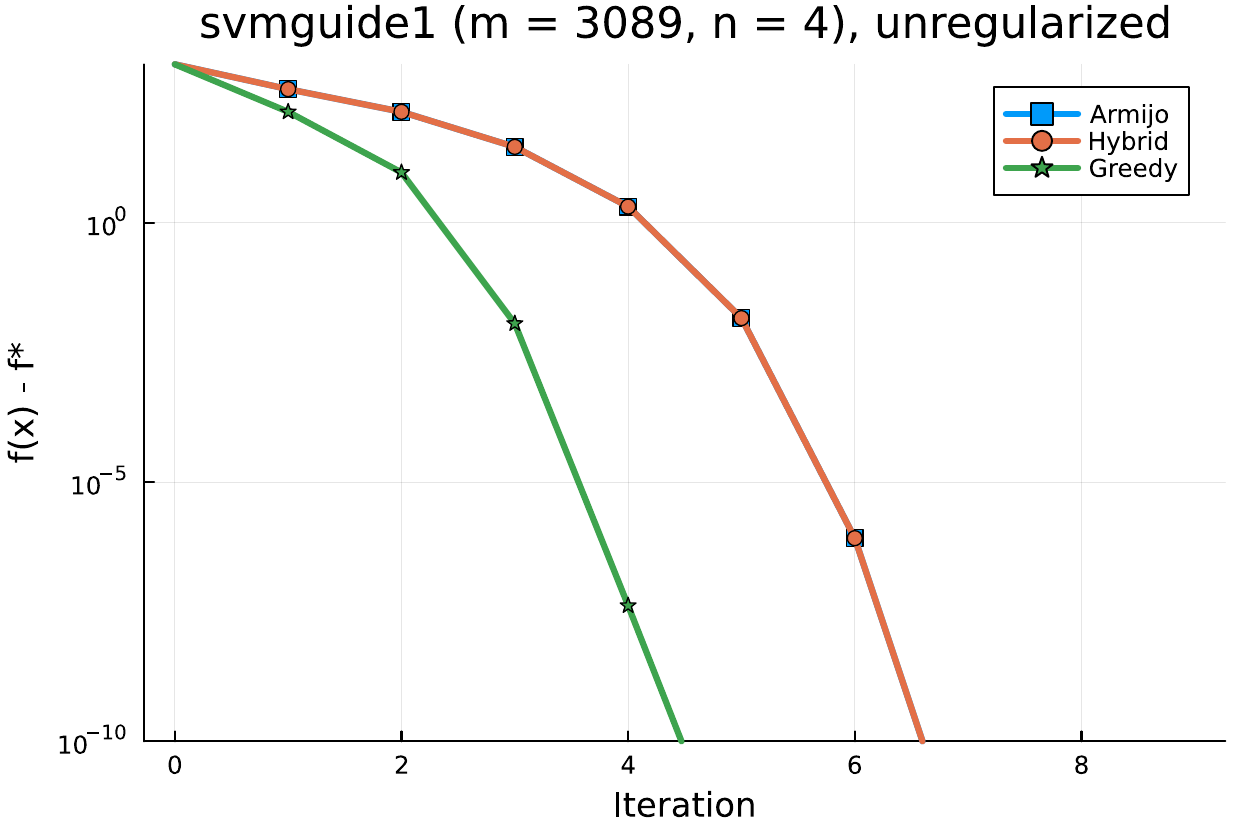}
	\includegraphics[width=.24\textwidth]{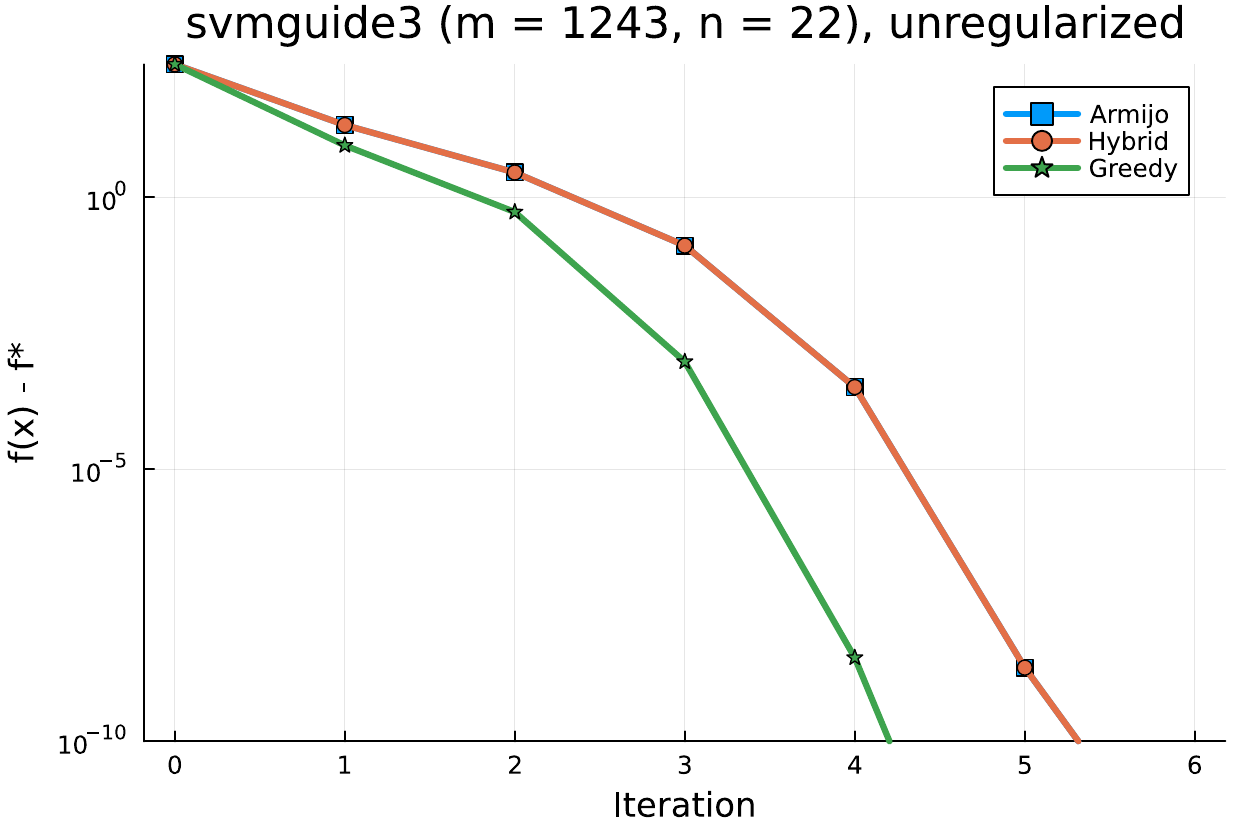}
	\includegraphics[width=.24\textwidth]{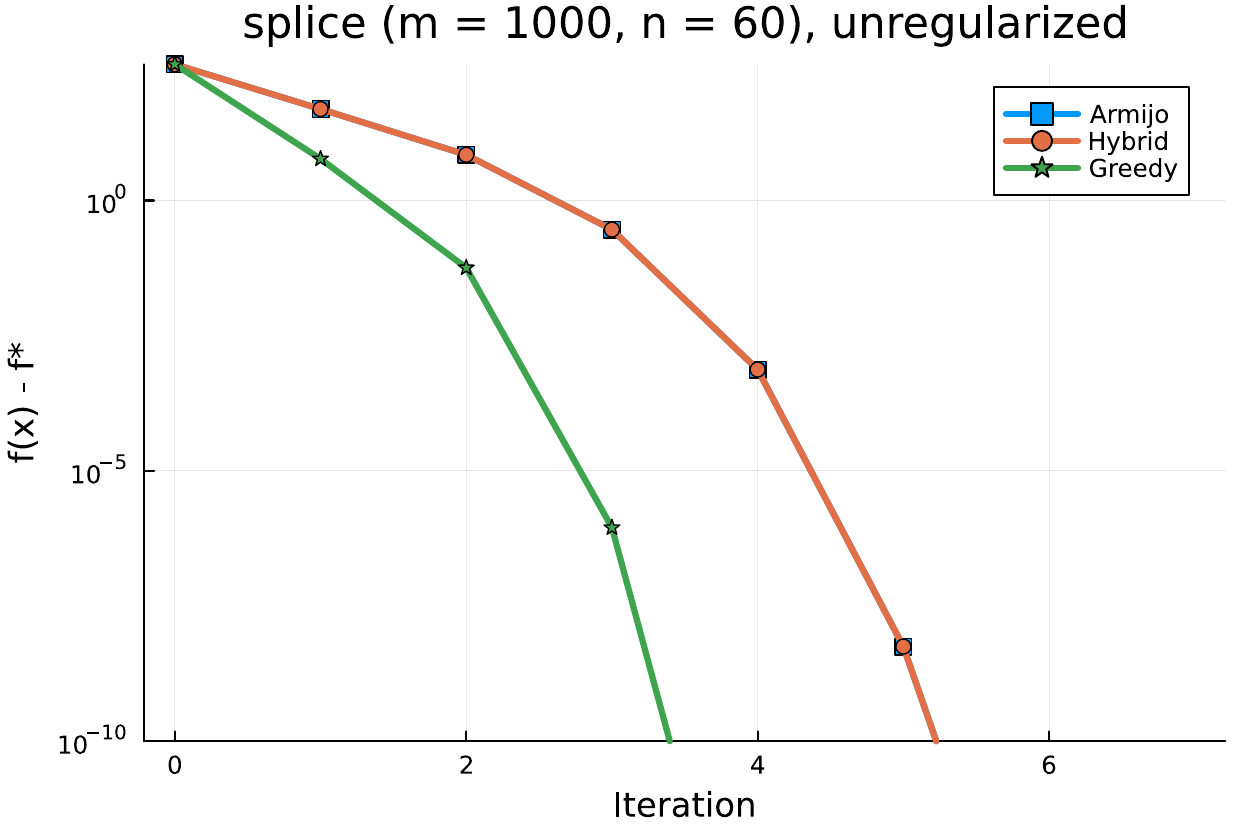}
	\includegraphics[width=.24\textwidth]{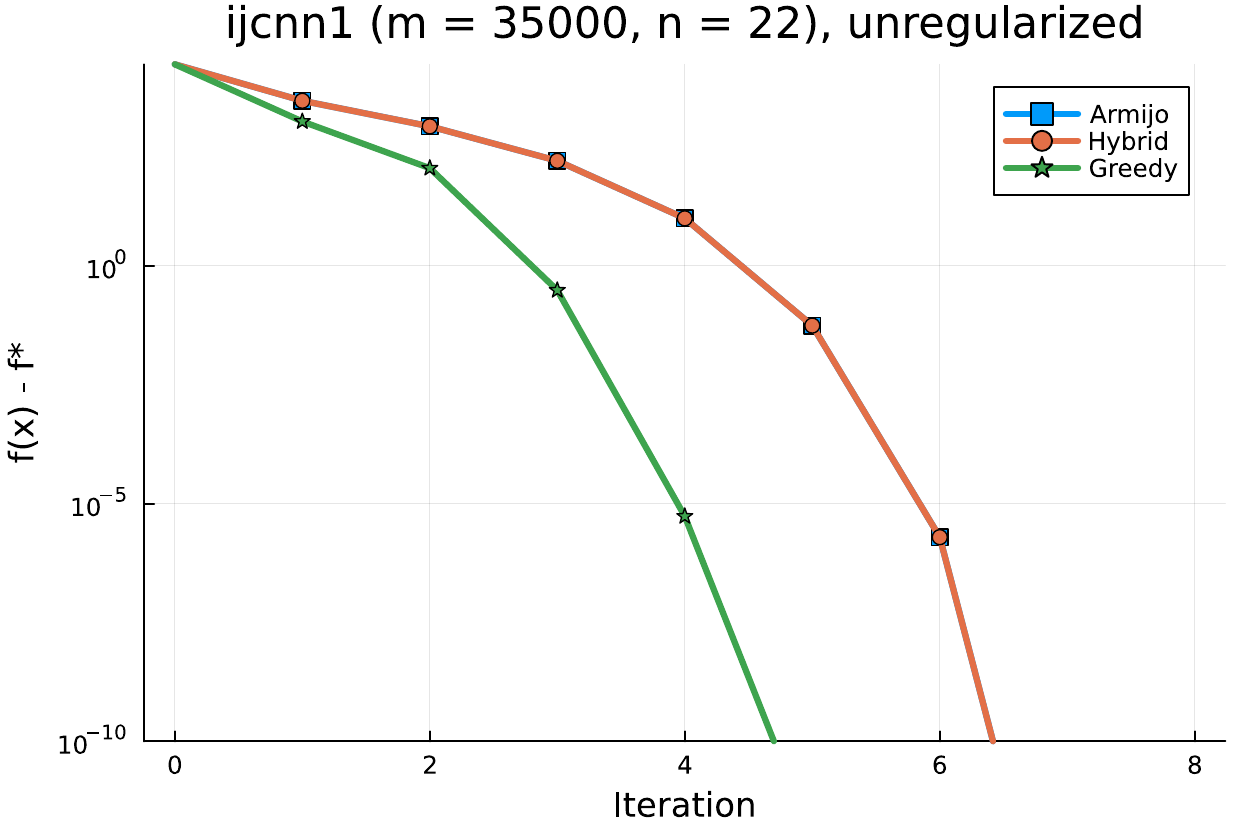}
	\caption{
		Comparison of methods on real logistic regression datasets, in 8 cases where we observed typical performance.
	}
	\label{fig:logregRealf1}
\end{figure}

\begin{figure}
	\includegraphics[width=.24\textwidth]{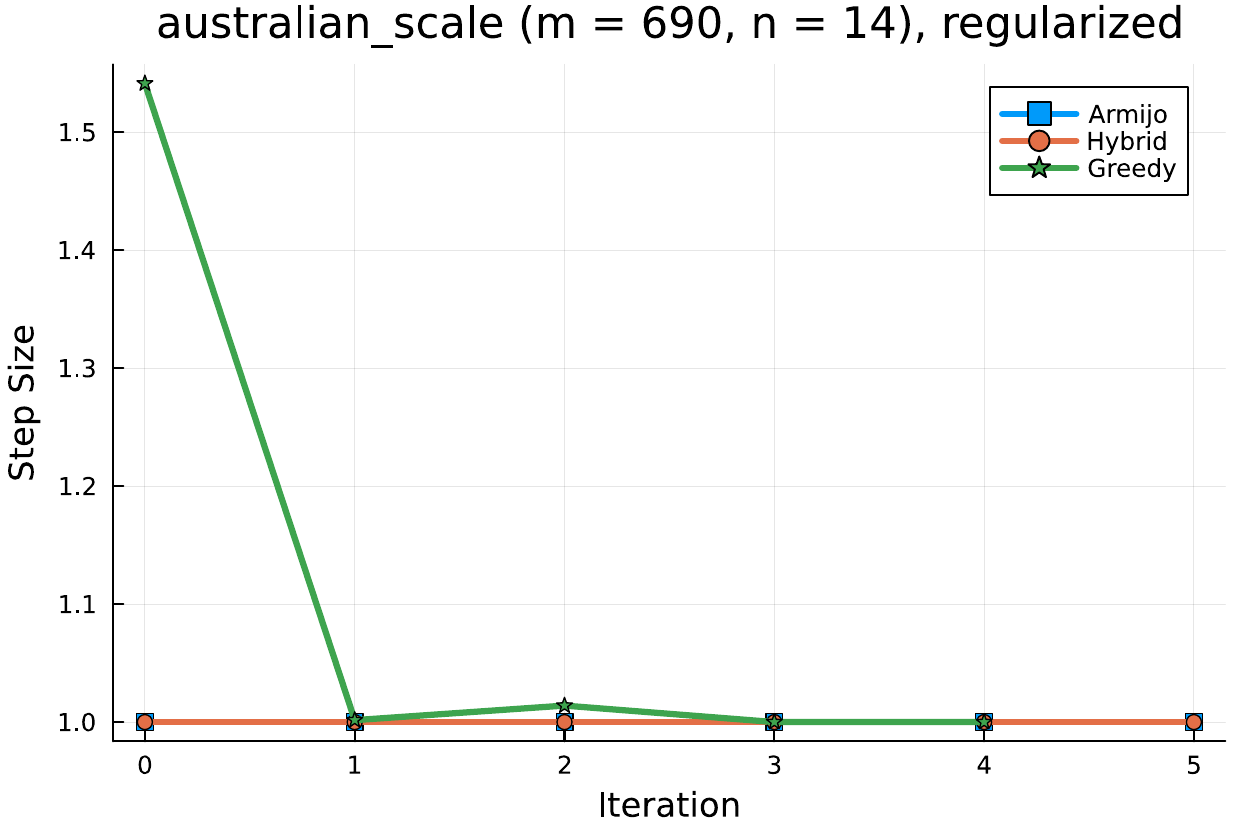}
	\includegraphics[width=.24\textwidth]{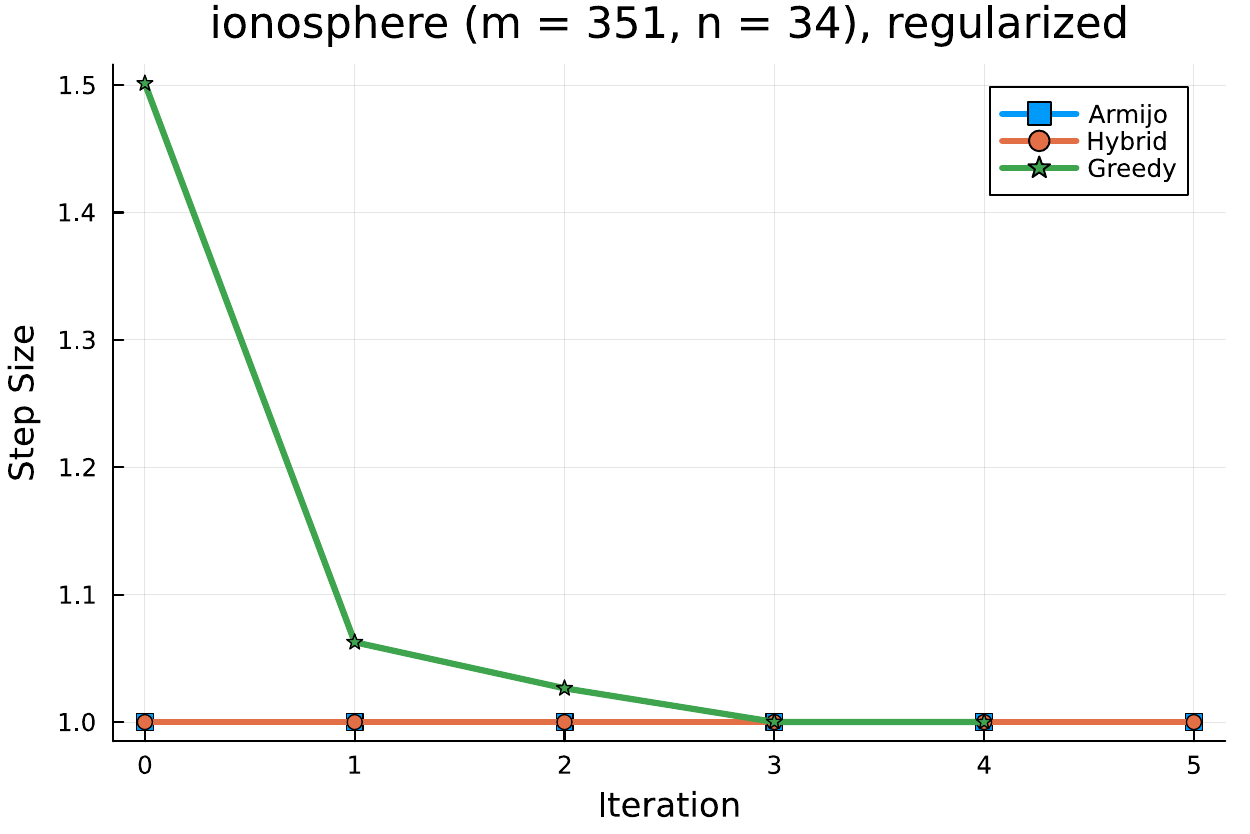}
	\includegraphics[width=.24\textwidth]{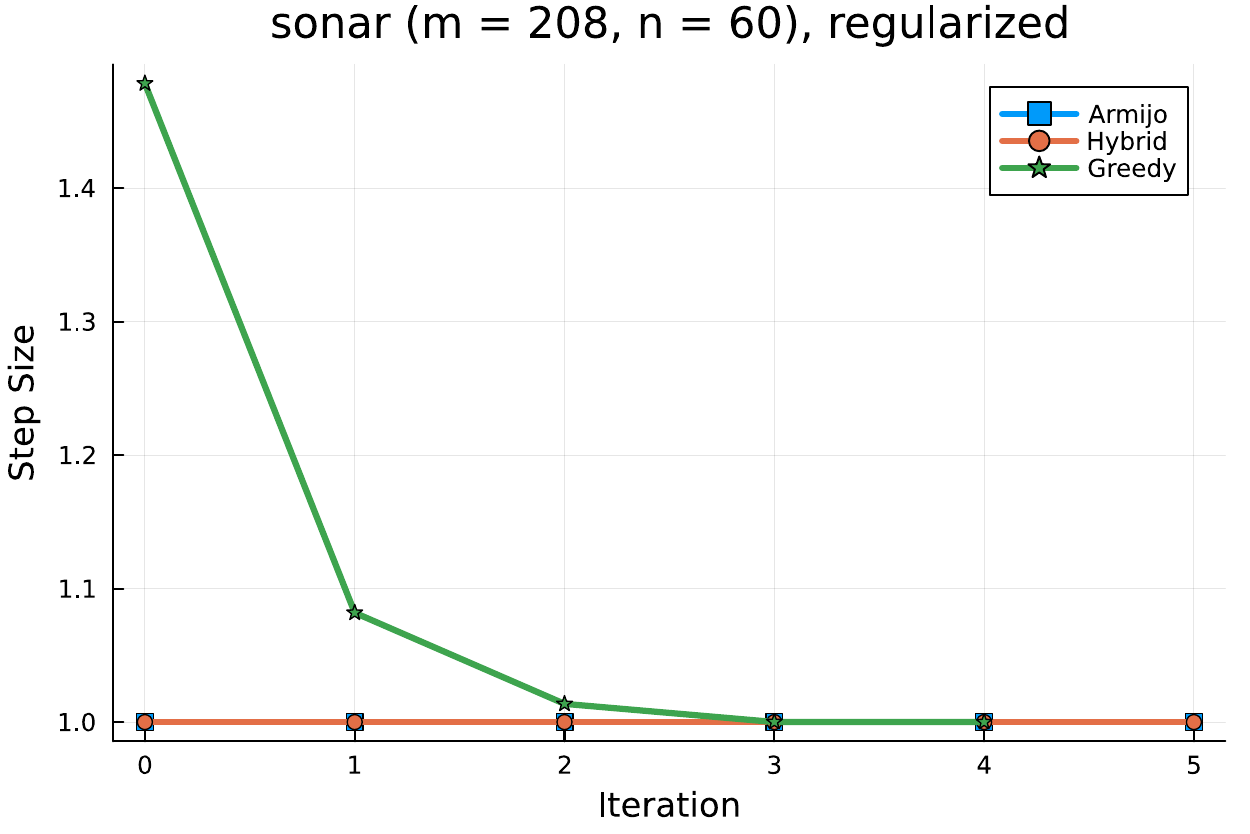}
	\includegraphics[width=.24\textwidth]{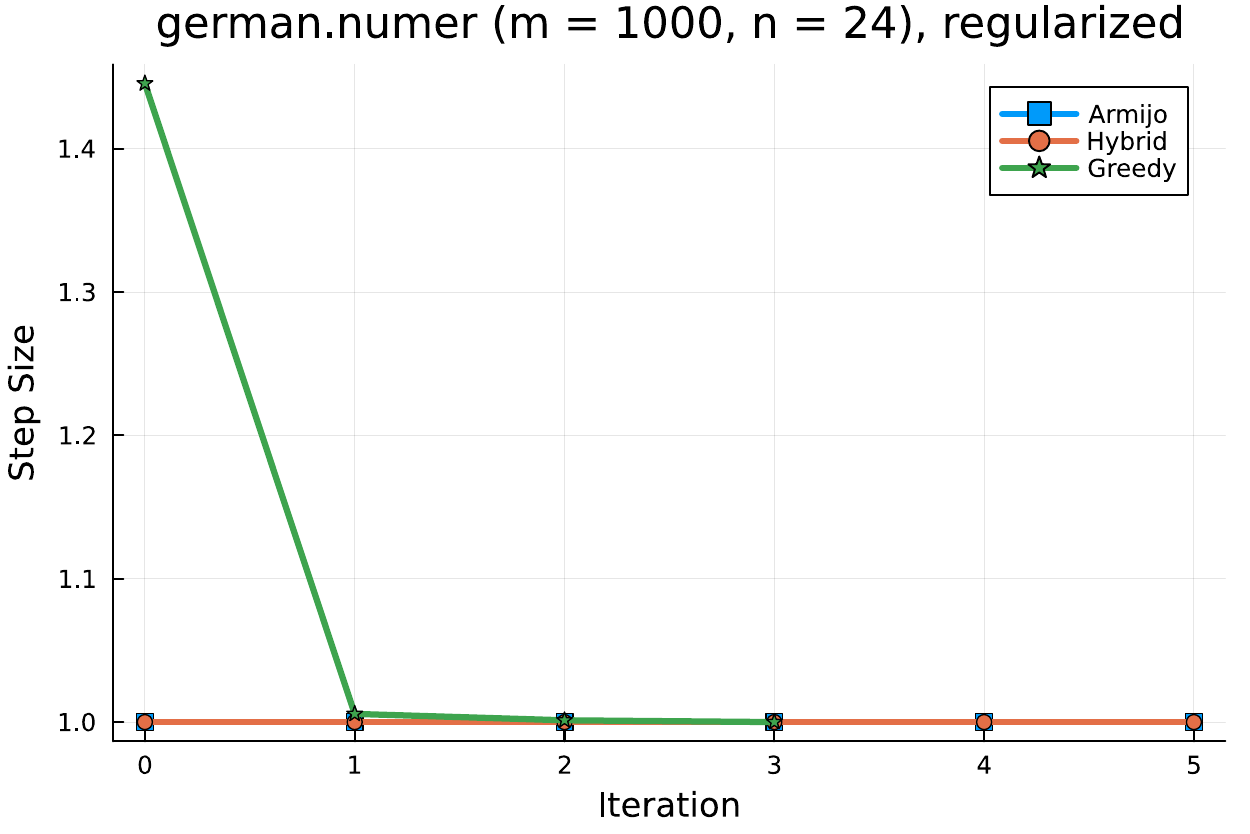}
	\includegraphics[width=.24\textwidth]{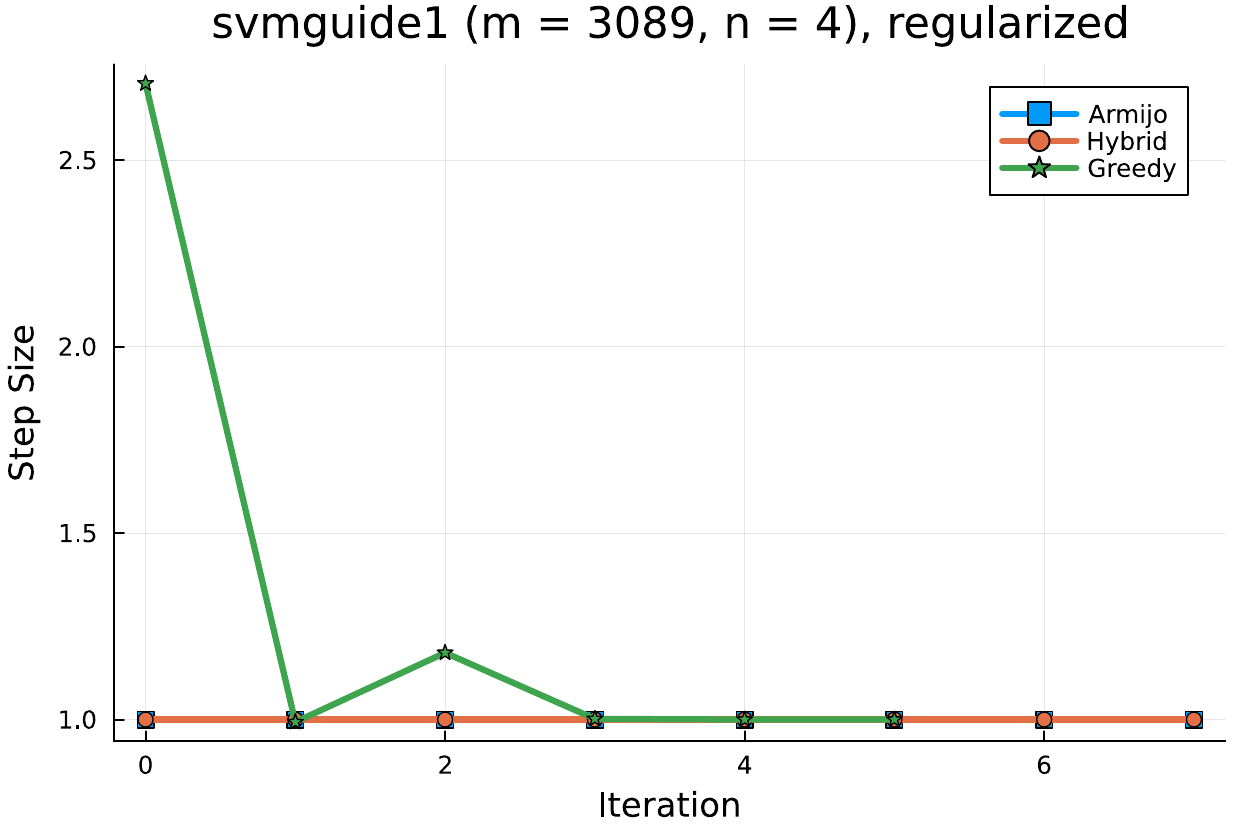}
	\includegraphics[width=.24\textwidth]{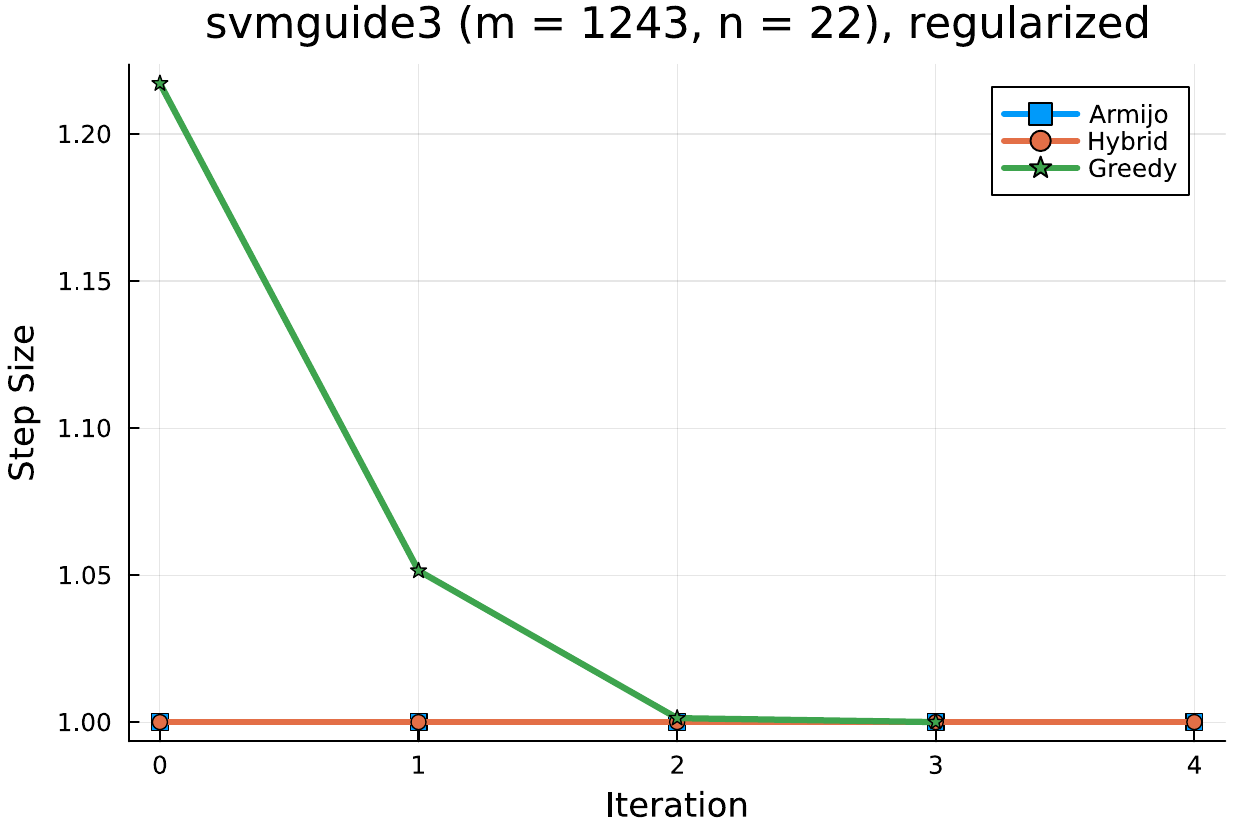}
	\includegraphics[width=.24\textwidth]{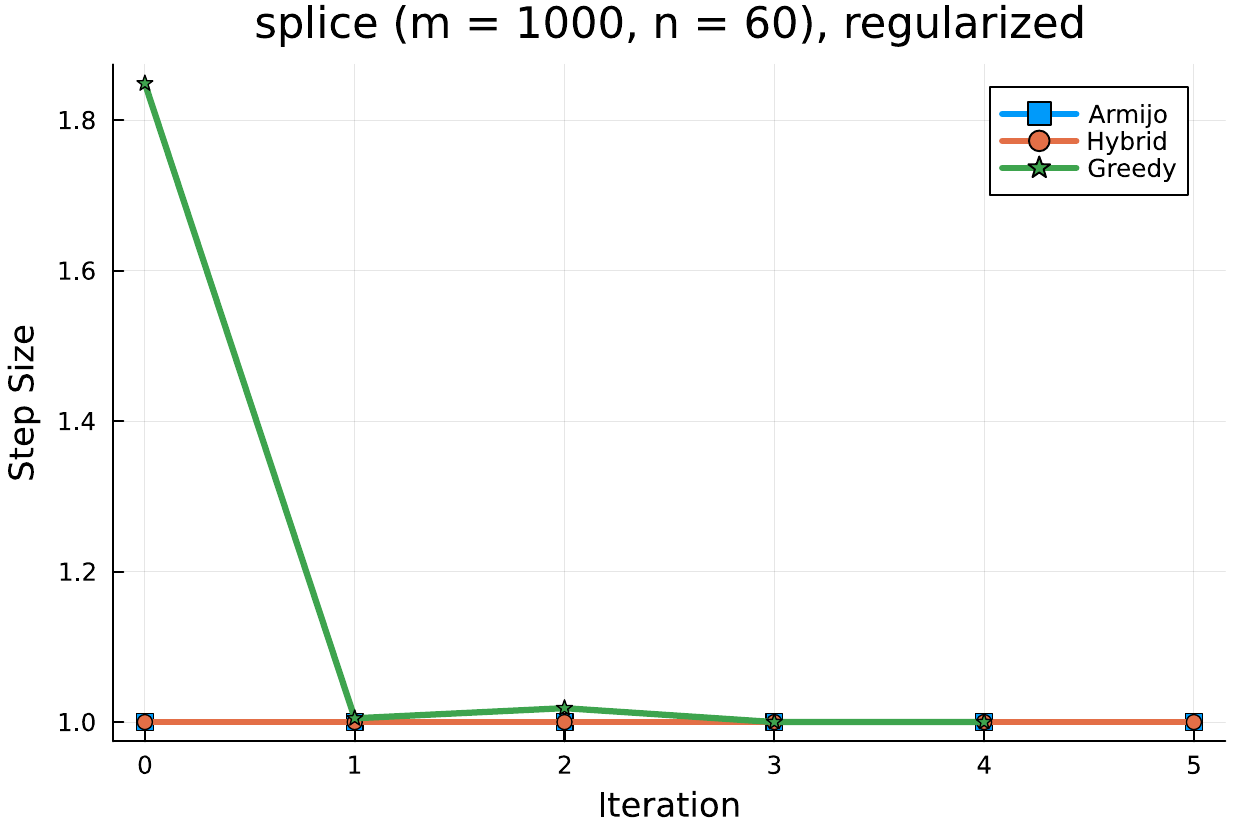}
	\includegraphics[width=.24\textwidth]{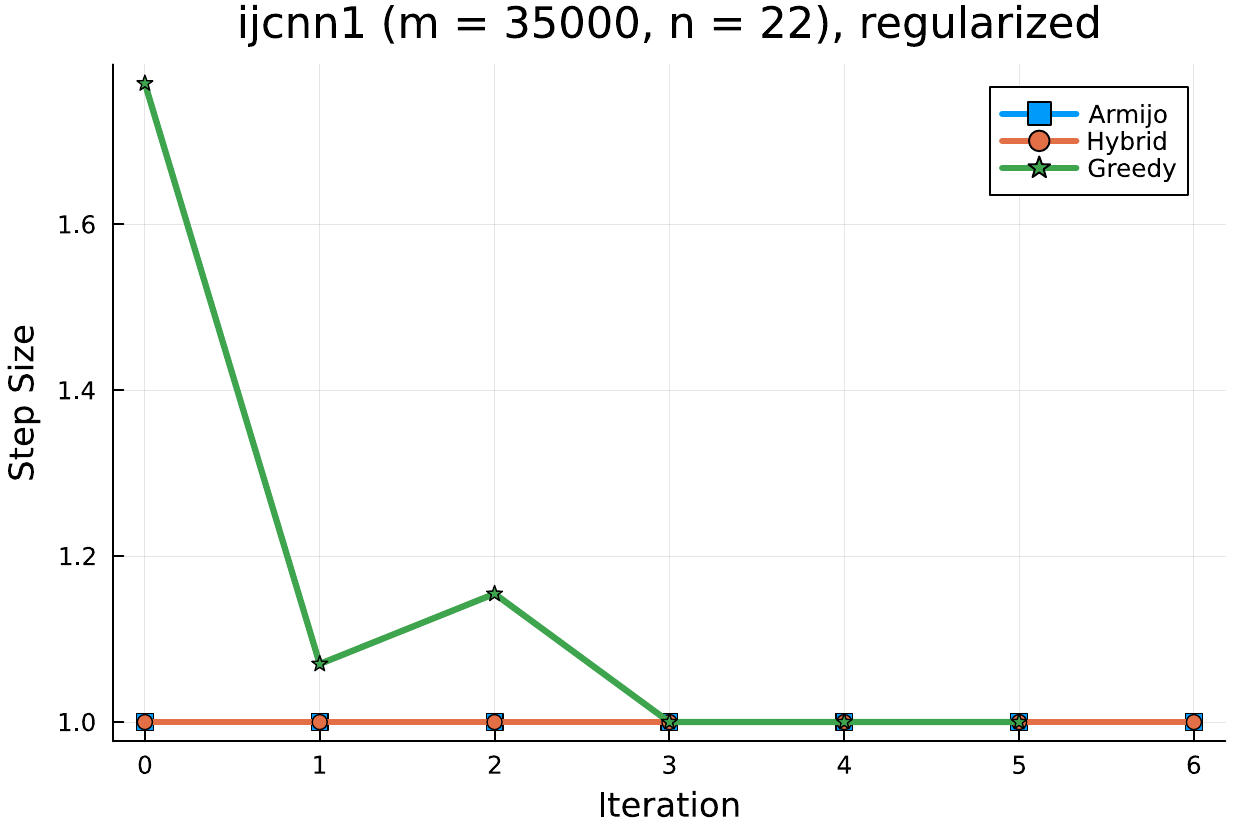}
	\includegraphics[width=.24\textwidth]{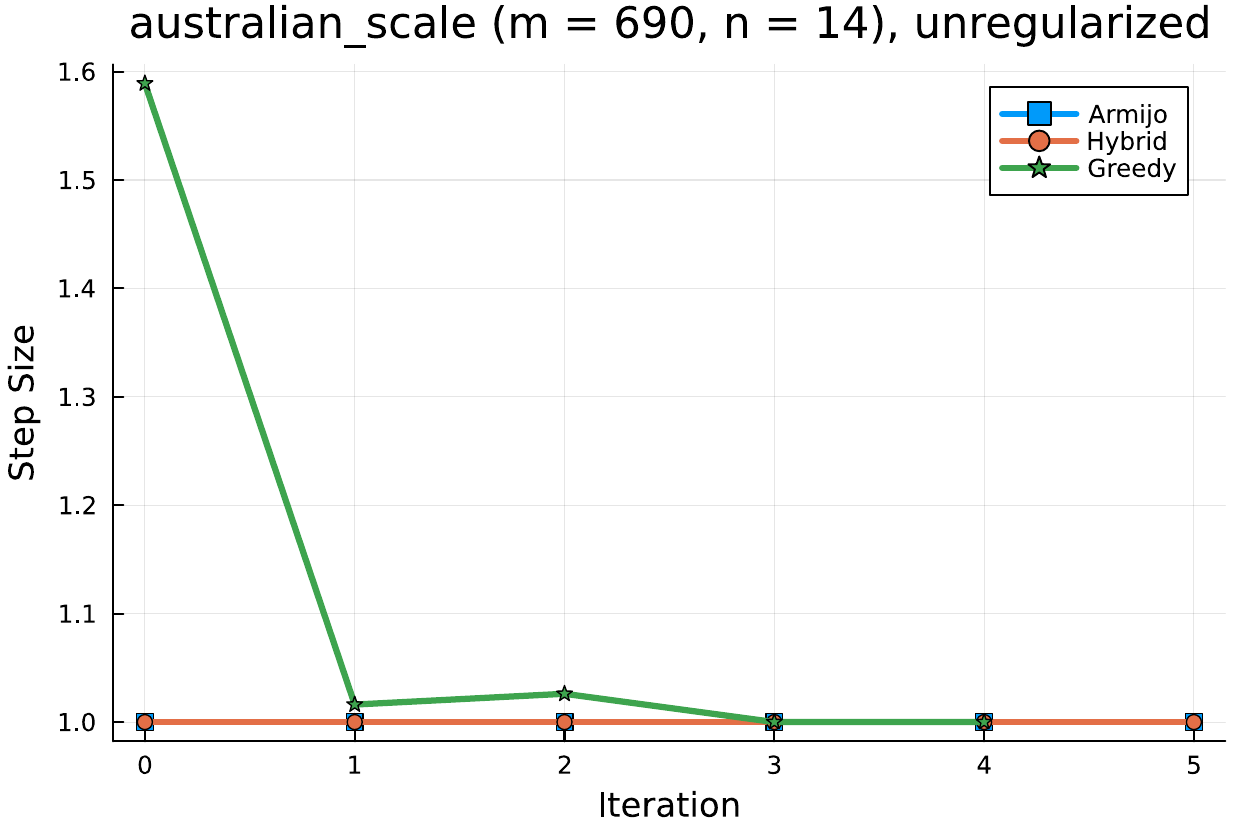}
	\includegraphics[width=.24\textwidth]{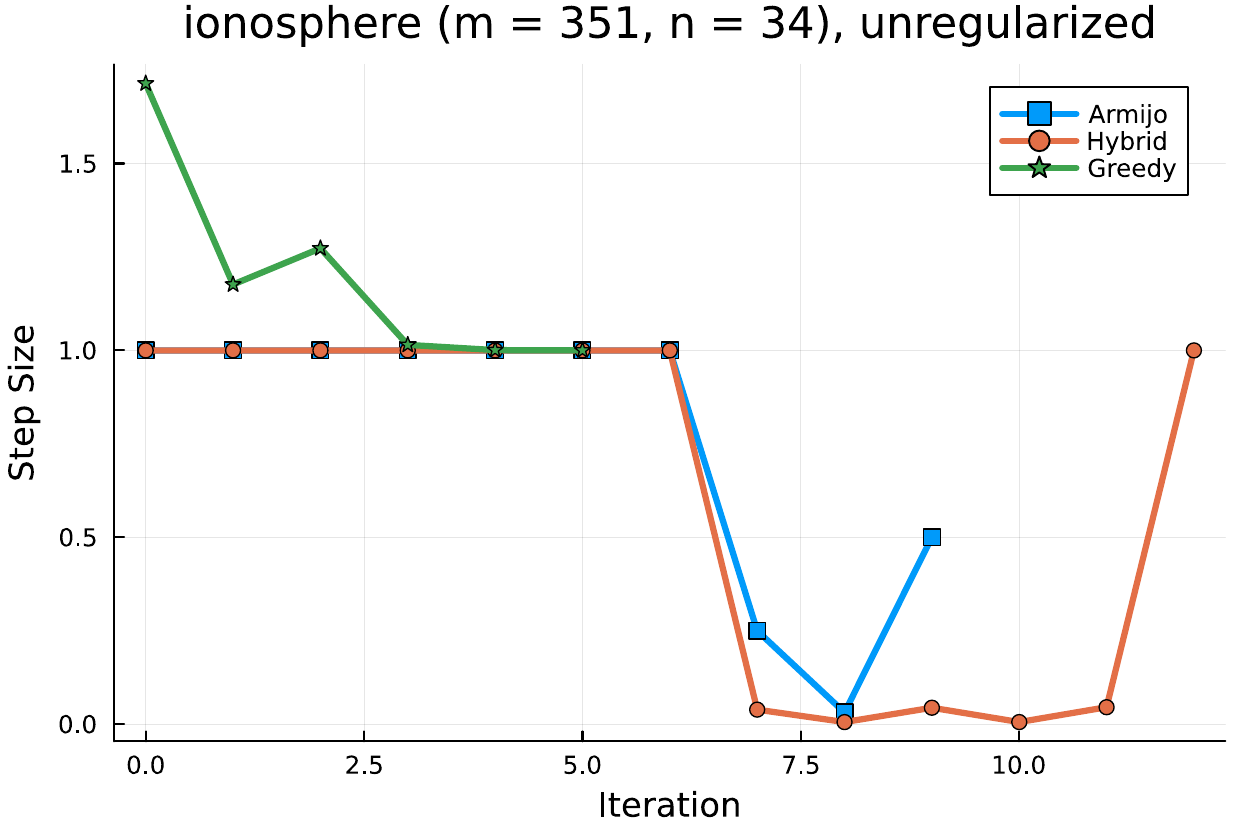}
	\includegraphics[width=.24\textwidth]{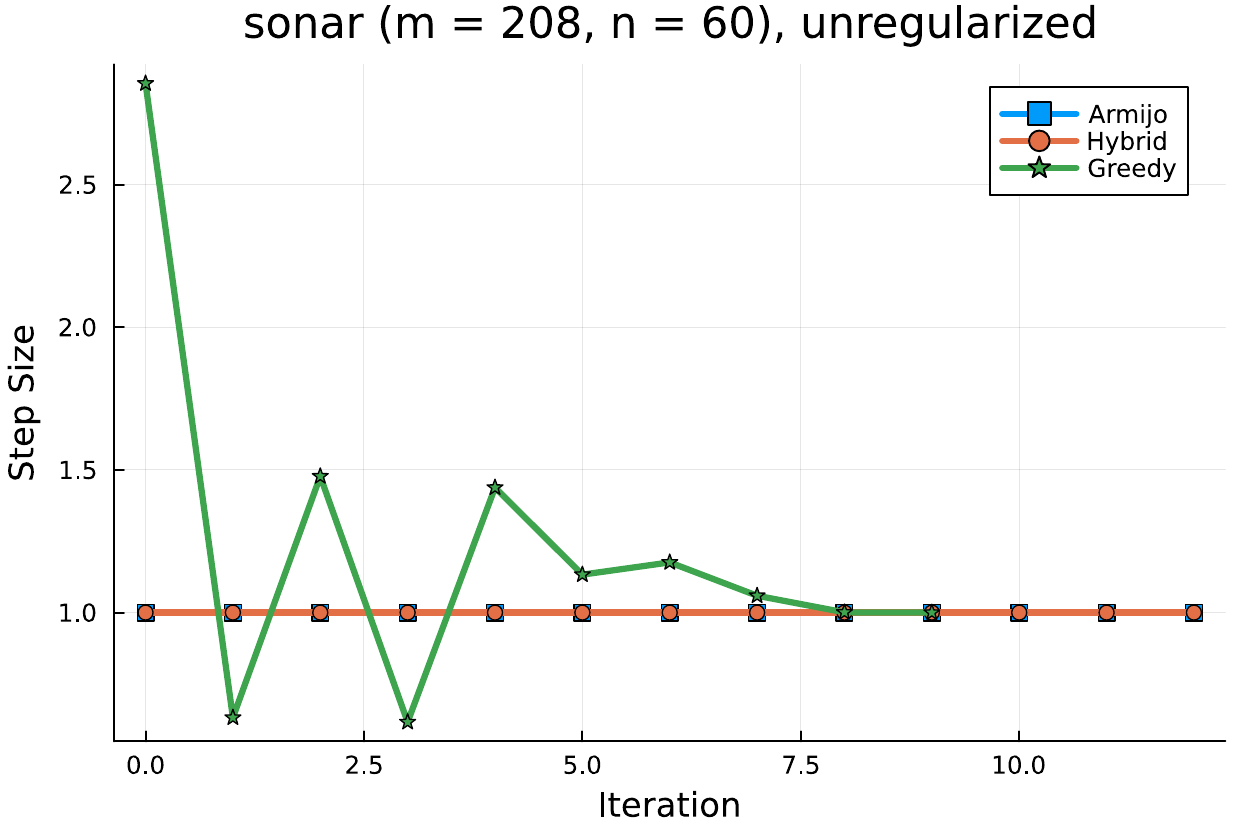}
	\includegraphics[width=.24\textwidth]{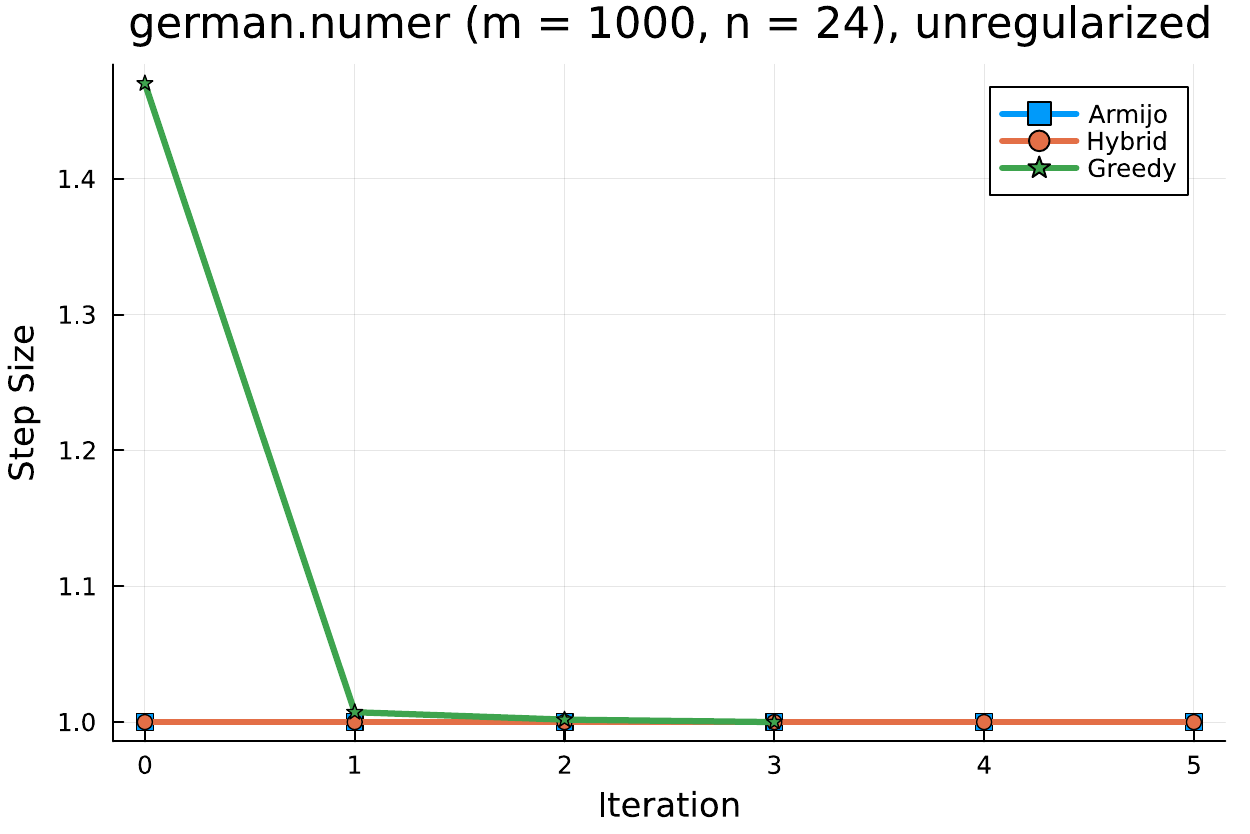}
	\includegraphics[width=.24\textwidth]{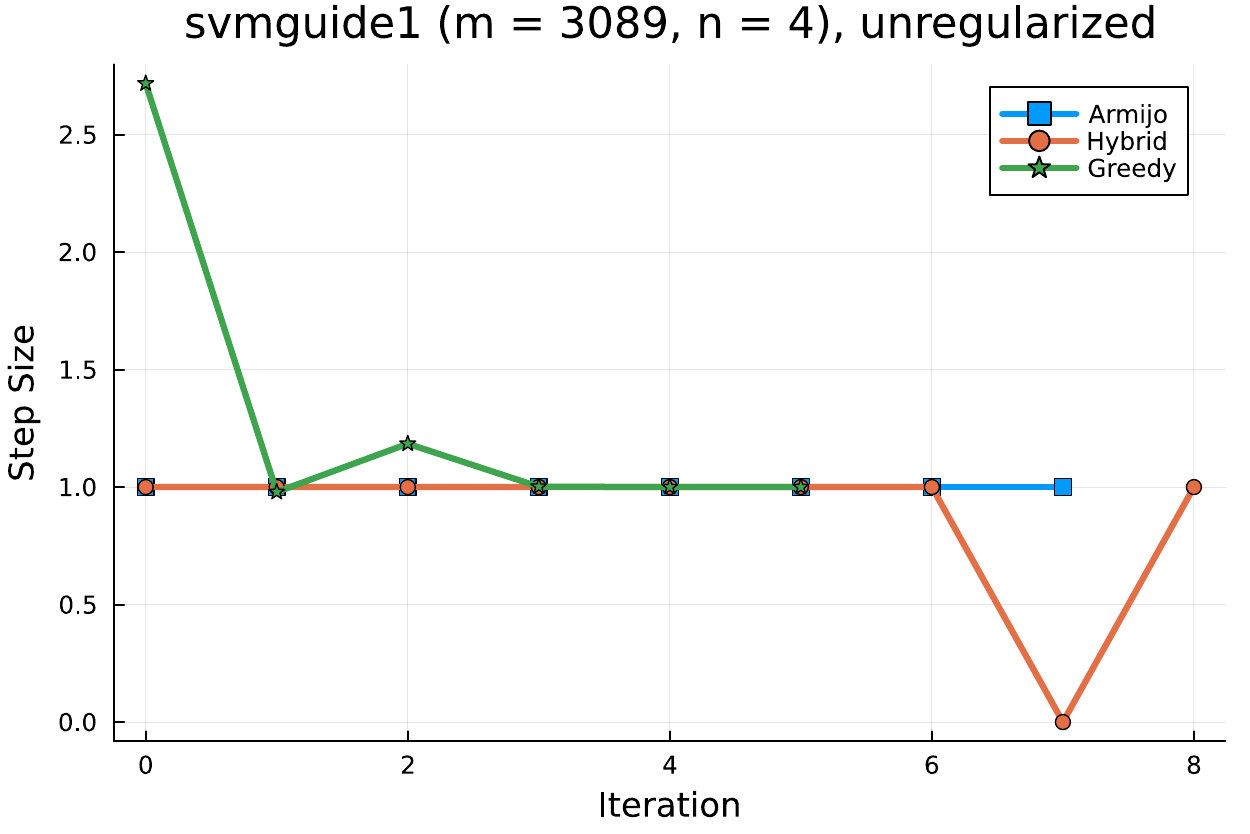}
	\includegraphics[width=.24\textwidth]{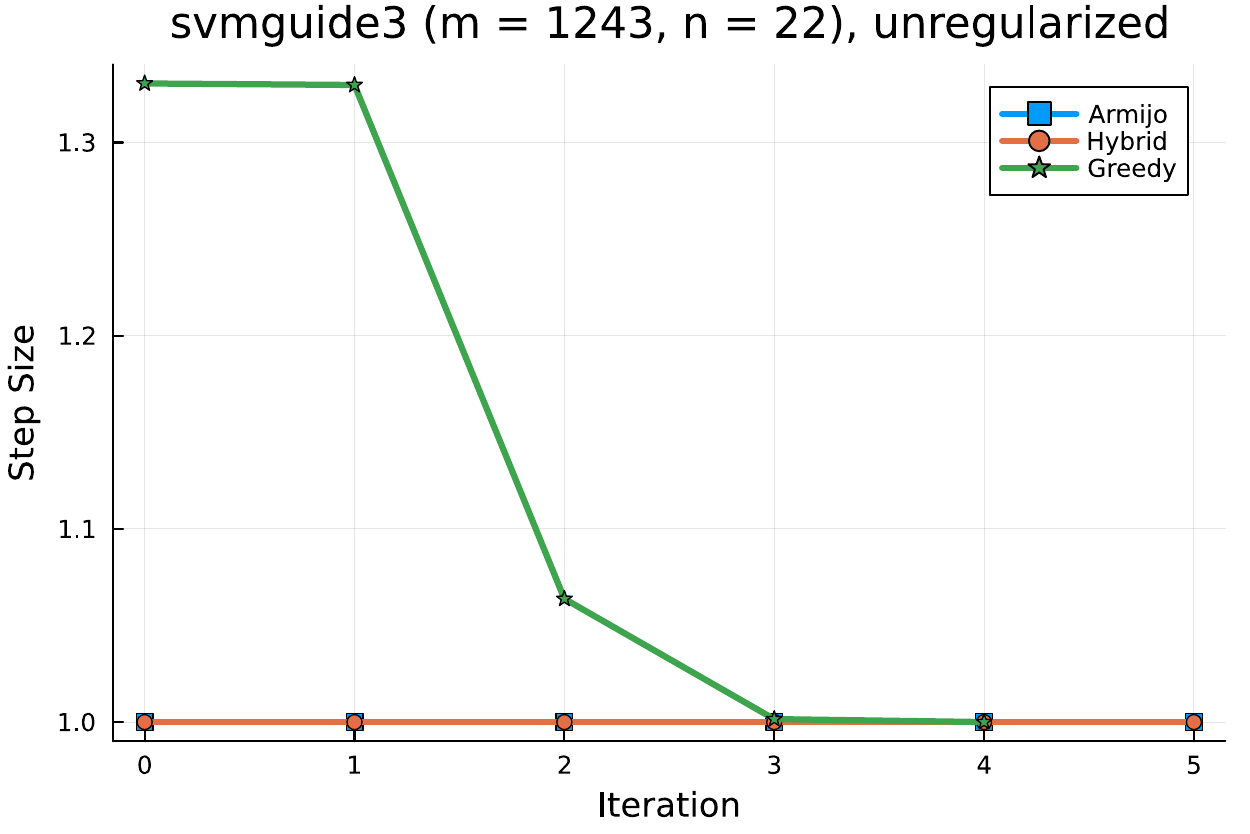}
	\includegraphics[width=.24\textwidth]{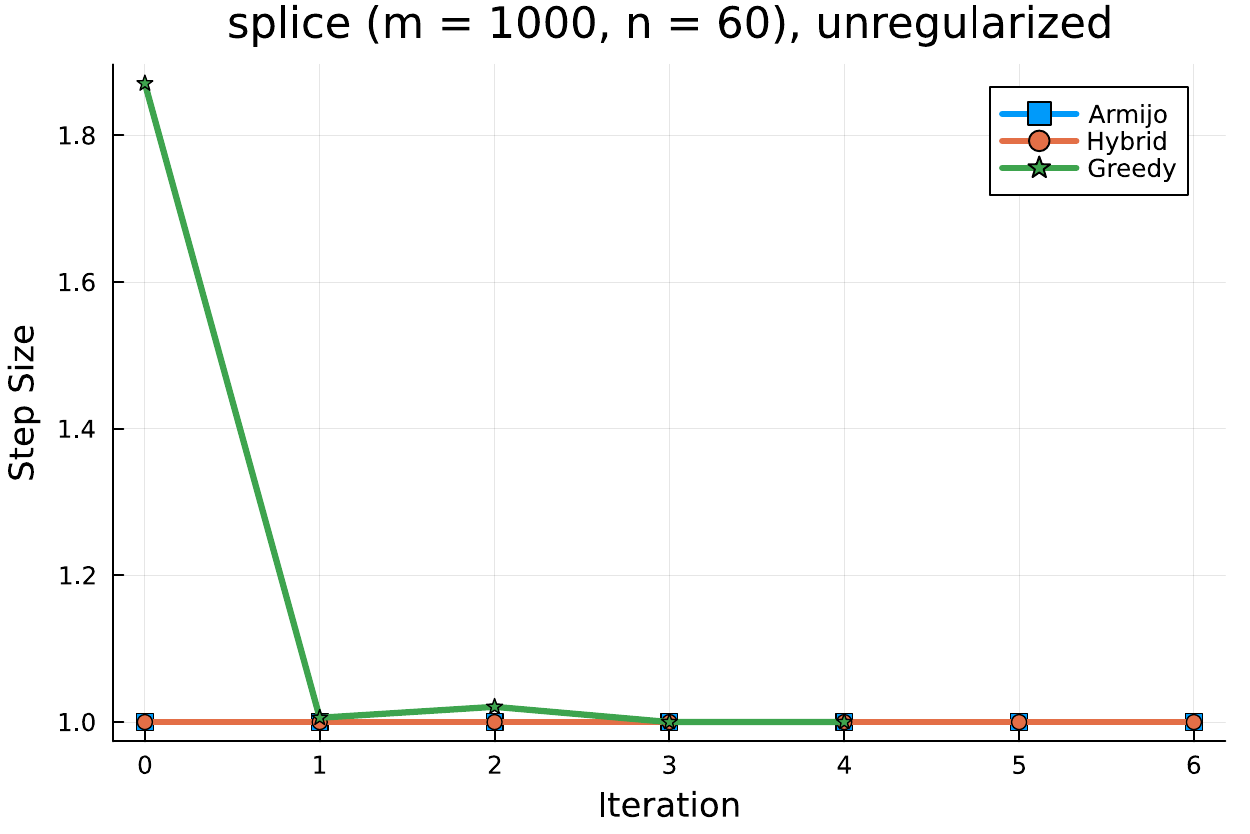}
	\includegraphics[width=.24\textwidth]{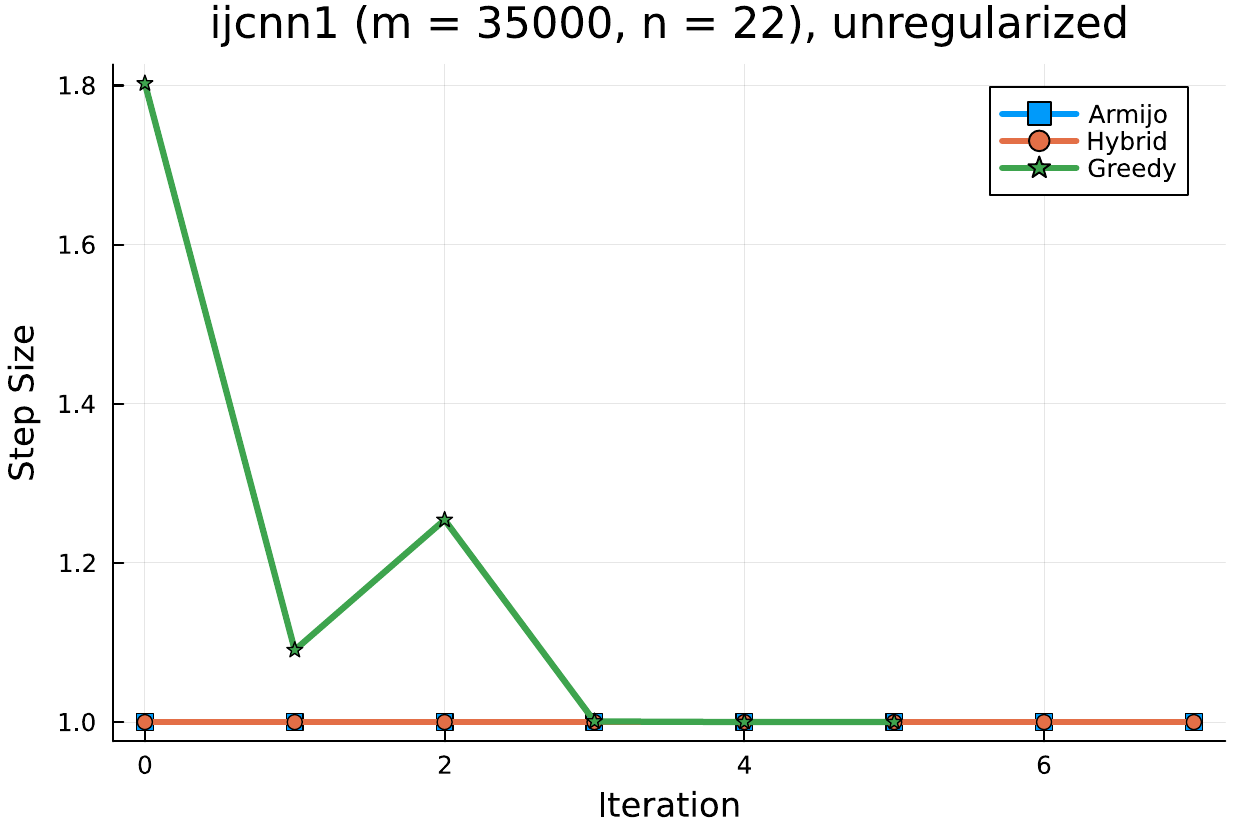}
	\caption{
		Step sizes of methods on real logistic regression datasets, in 8 cases where we observed typical performance.
	}
	\label{fig:logregRealt1}
\end{figure}

\begin{figure}
	\includegraphics[width=.24\textwidth]{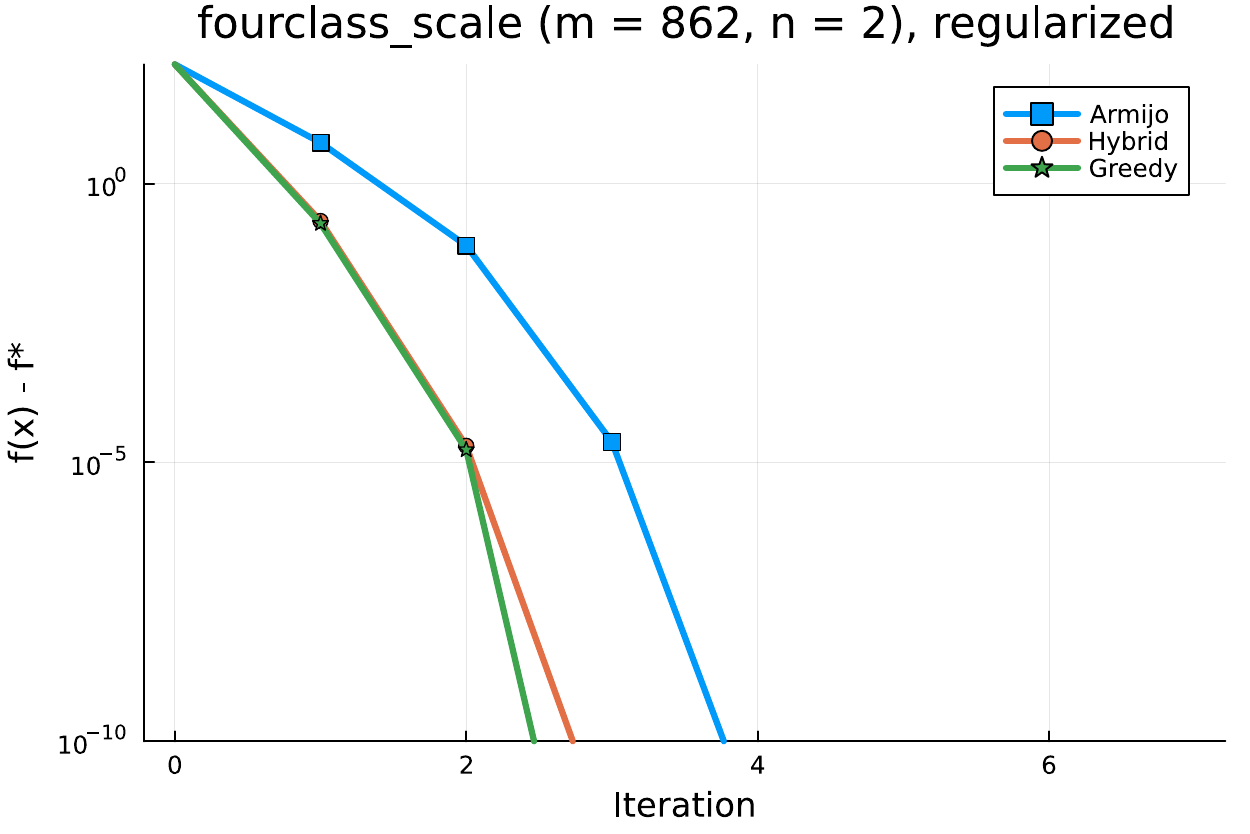}
	\includegraphics[width=.24\textwidth]{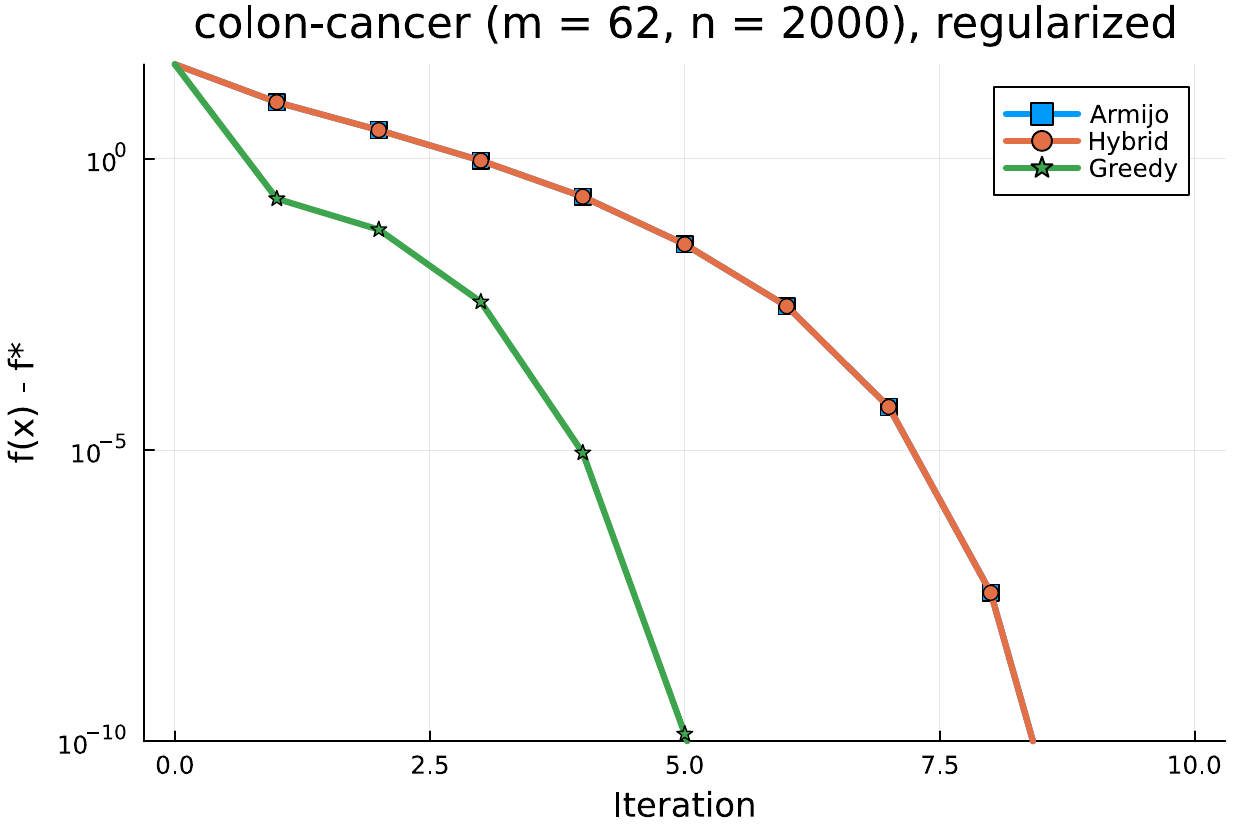}
	\includegraphics[width=.24\textwidth]{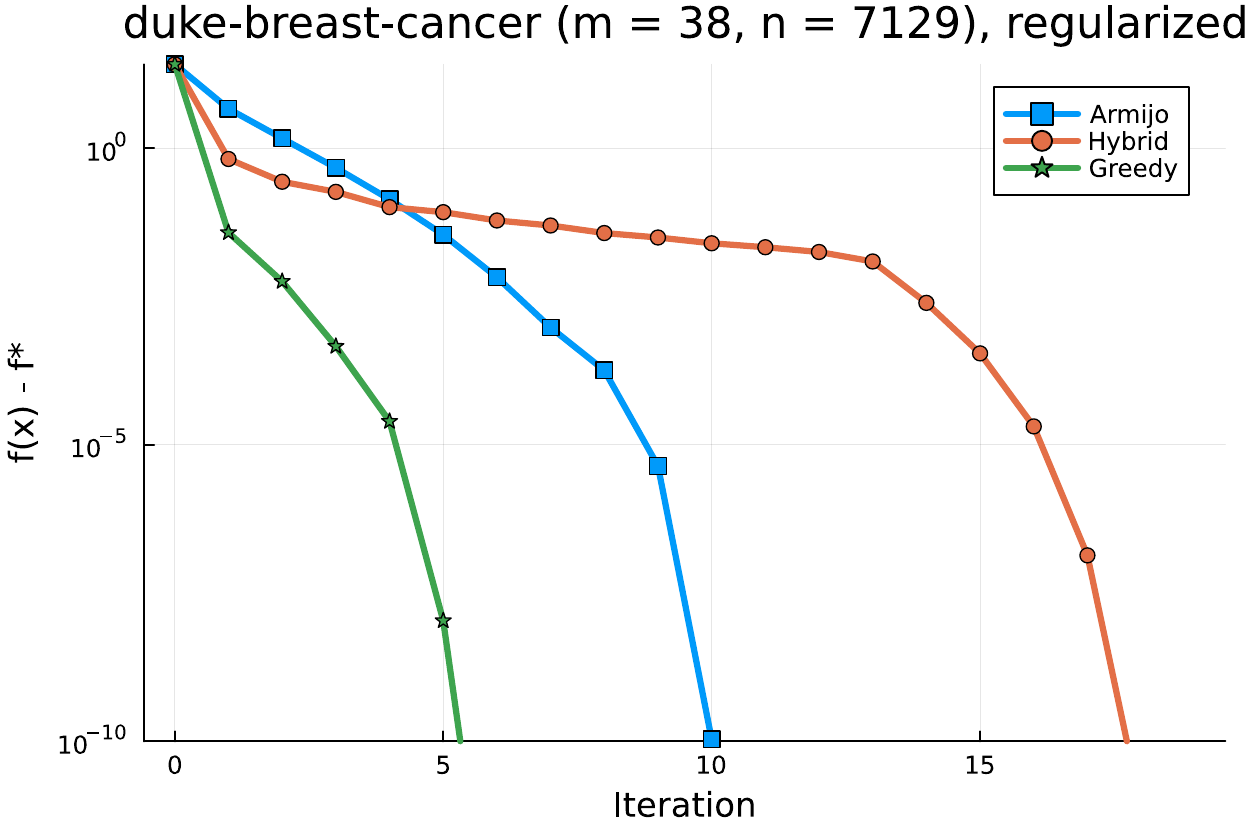}
	\includegraphics[width=.24\textwidth]{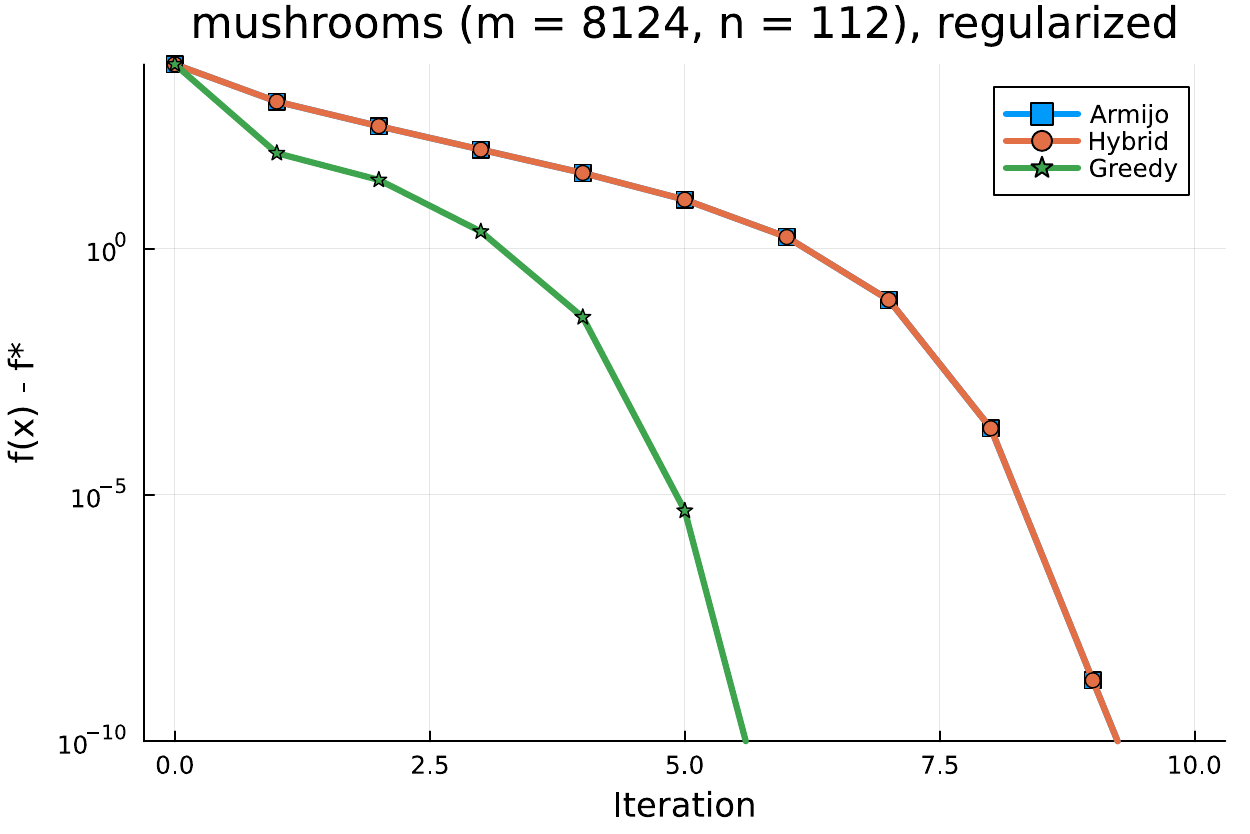}
	\includegraphics[width=.24\textwidth]{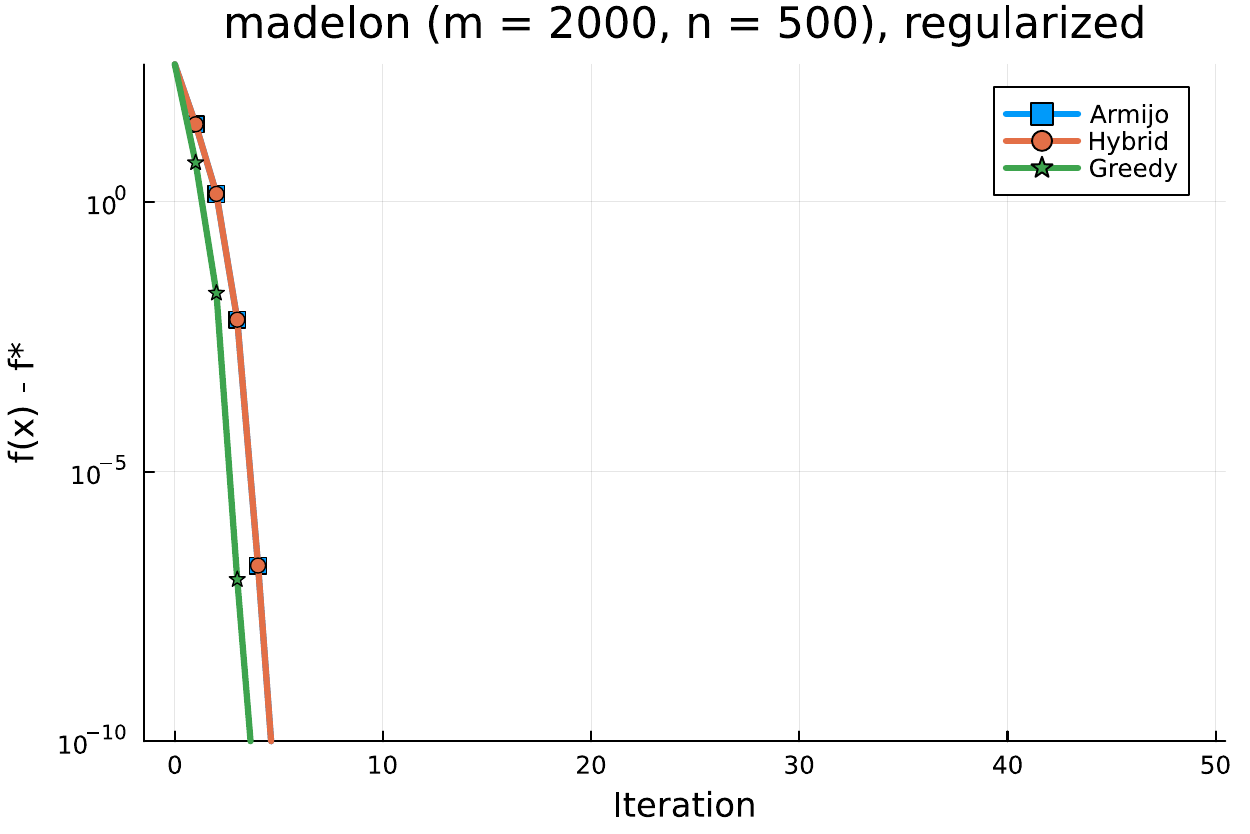}
	\includegraphics[width=.24\textwidth]{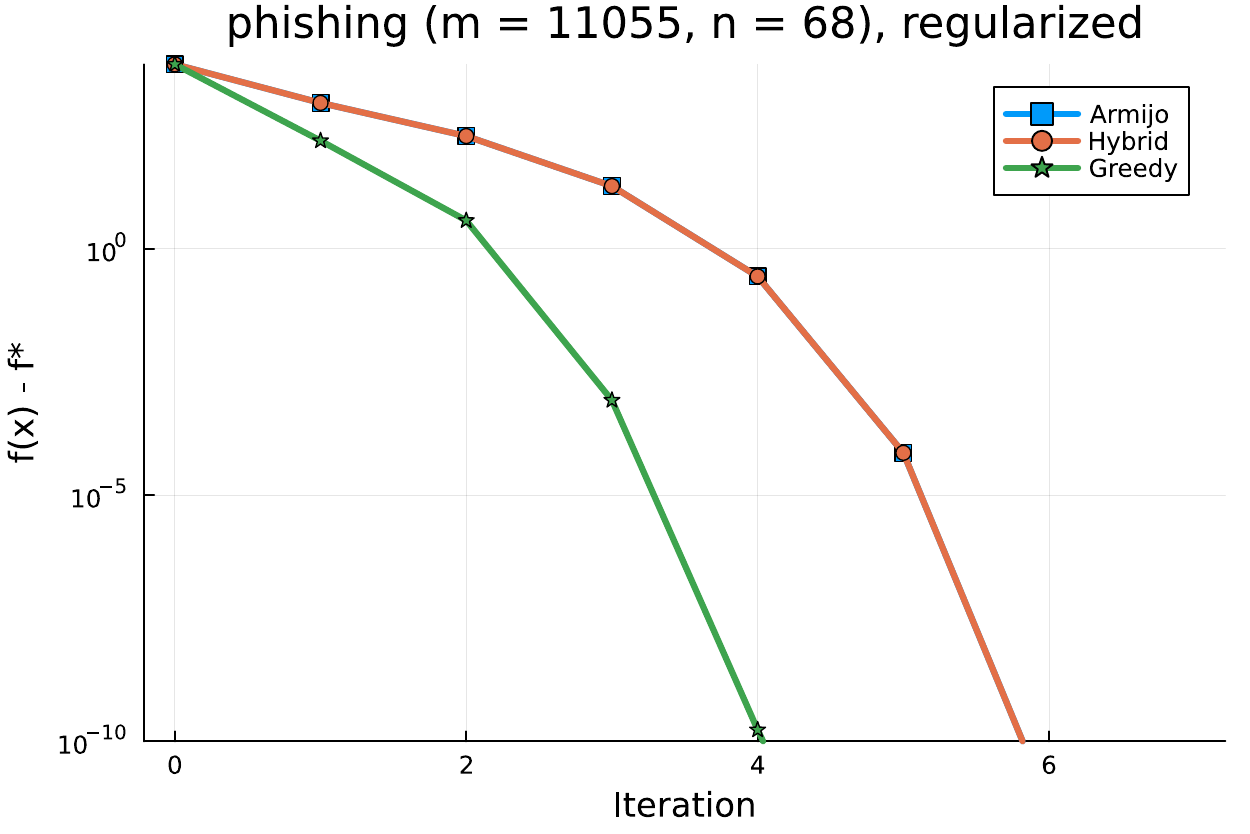}
	\includegraphics[width=.24\textwidth]{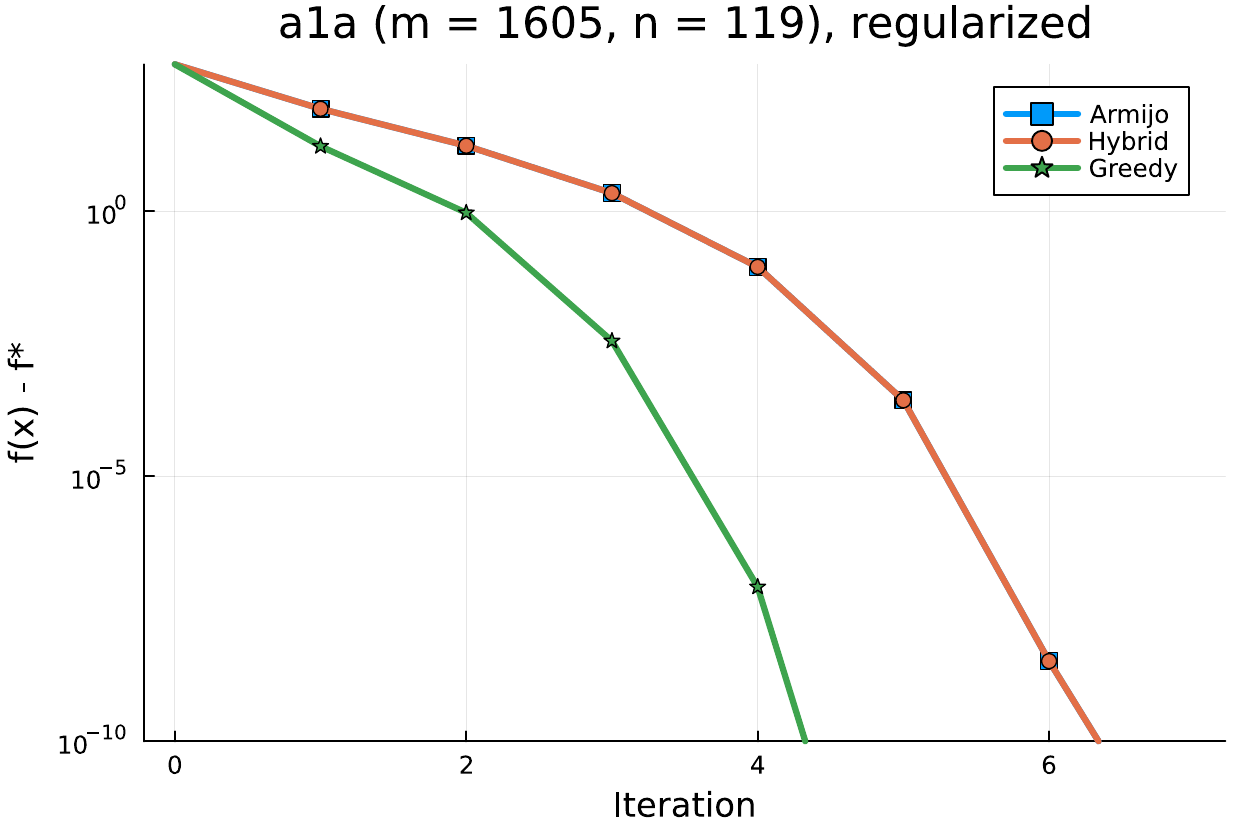}
	\includegraphics[width=.24\textwidth]{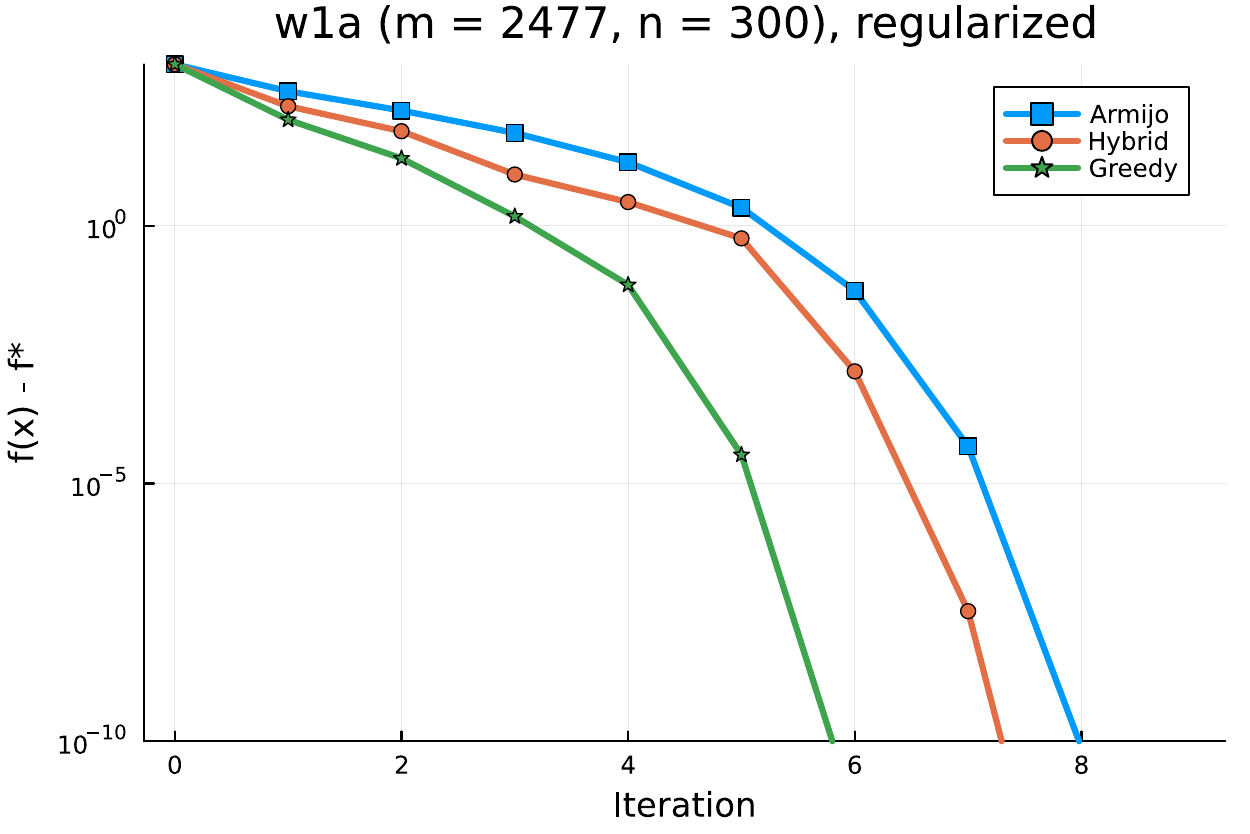}
	\includegraphics[width=.24\textwidth]{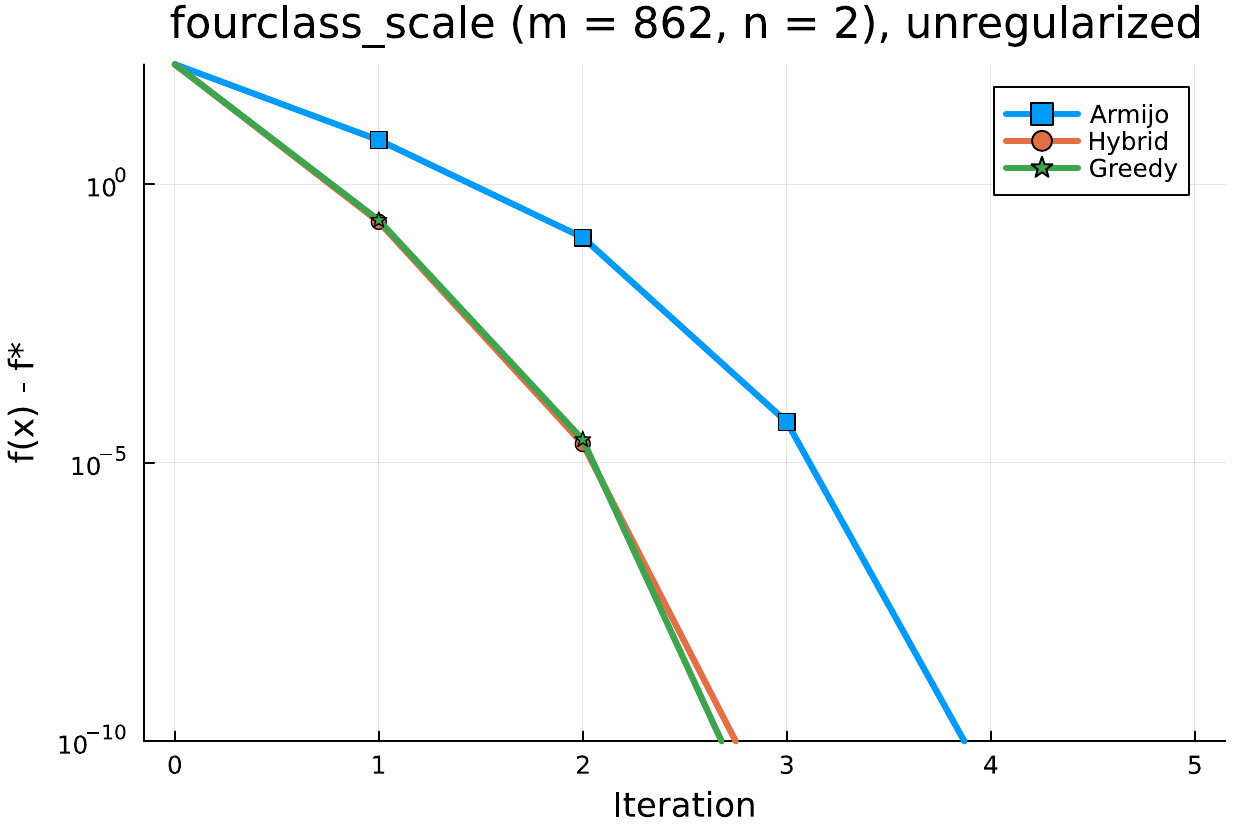}
	\includegraphics[width=.24\textwidth]{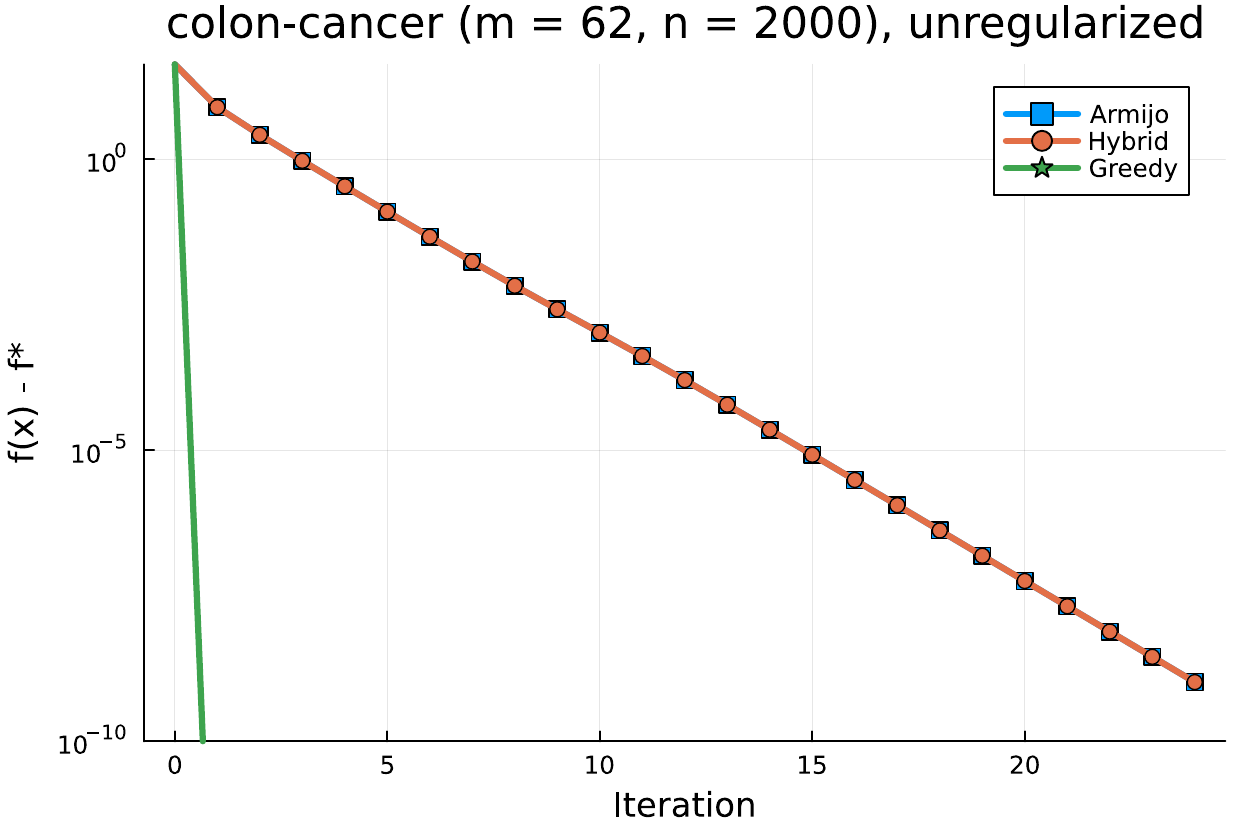}
	\includegraphics[width=.24\textwidth]{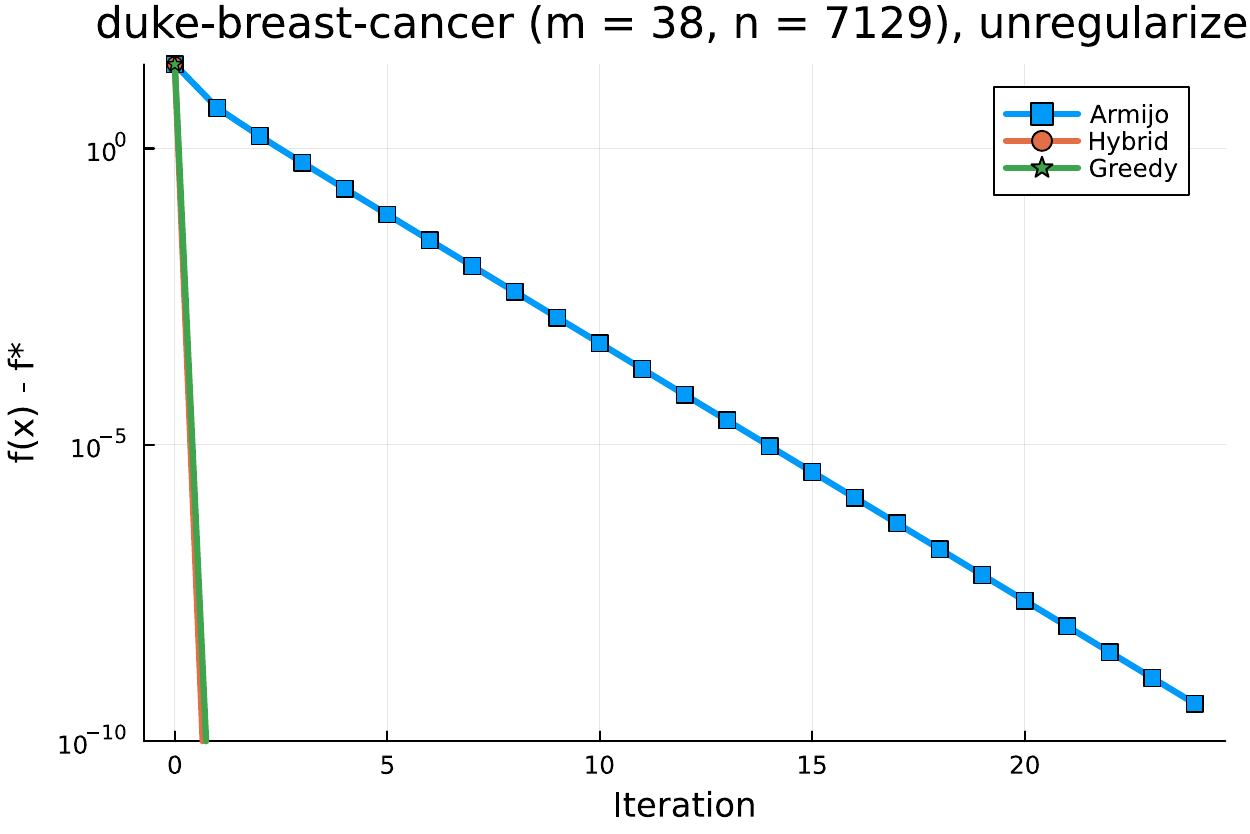}
	\includegraphics[width=.24\textwidth]{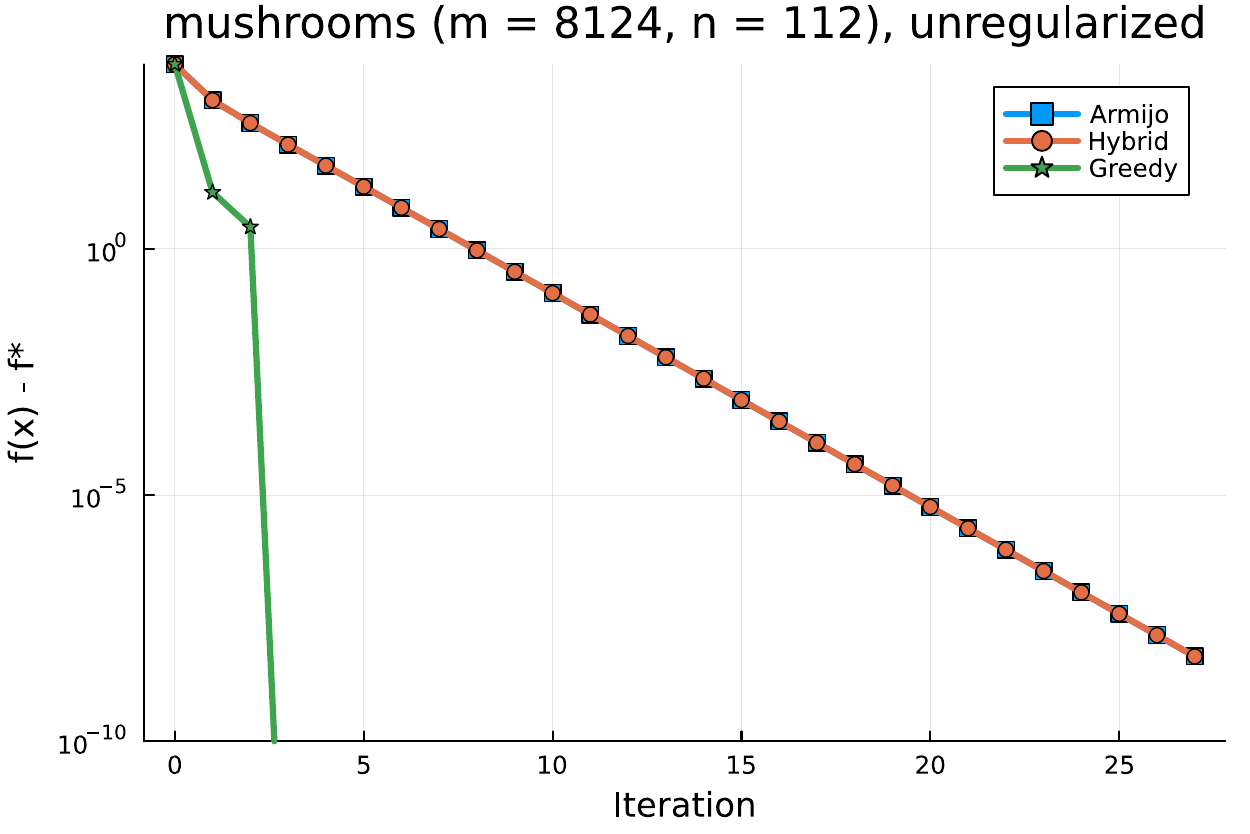}
	\includegraphics[width=.24\textwidth]{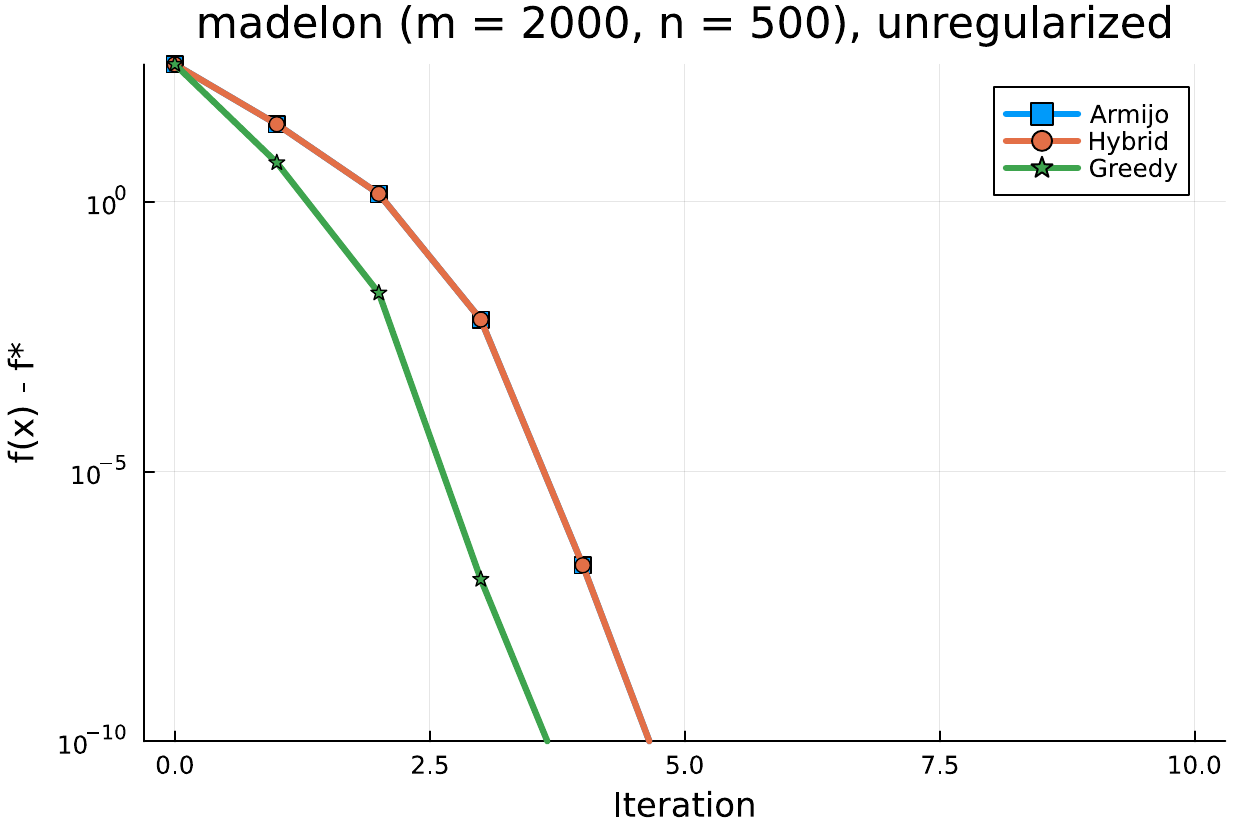}
	\includegraphics[width=.24\textwidth]{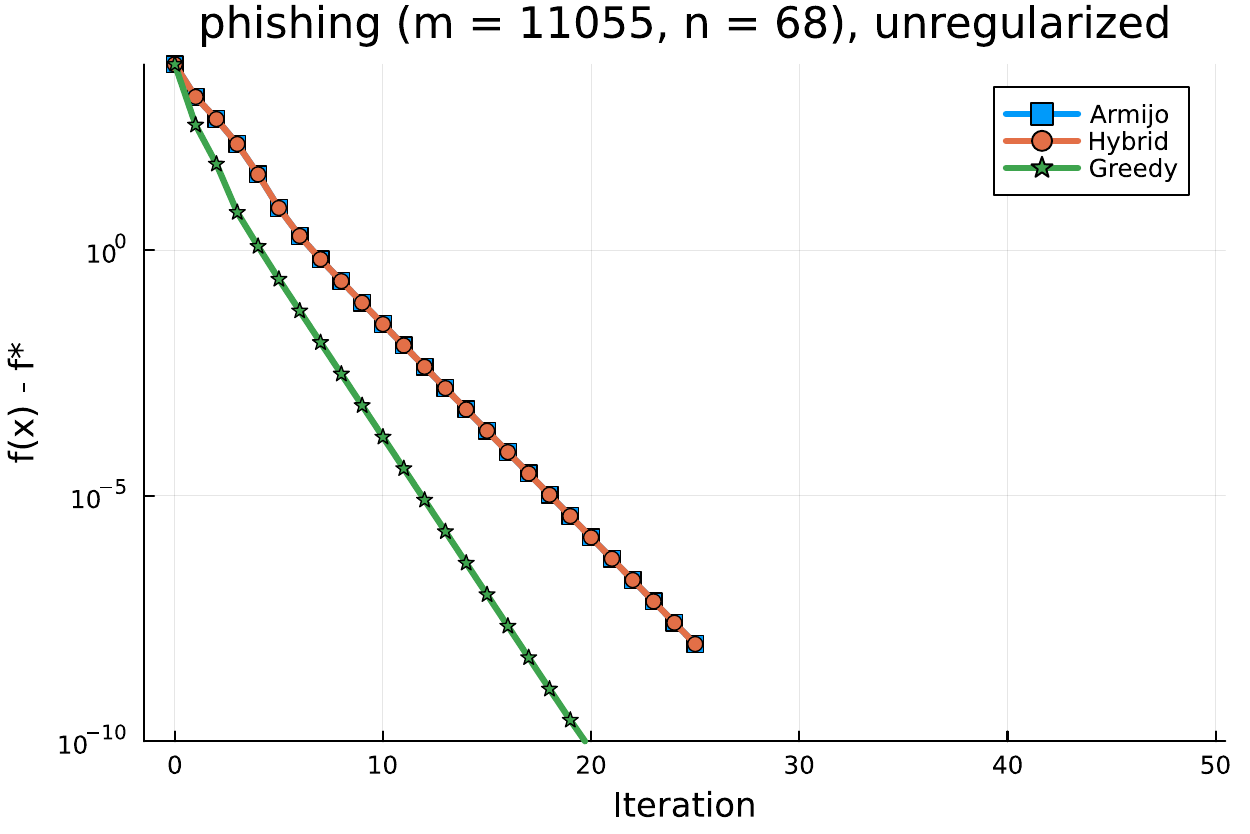}
	\includegraphics[width=.24\textwidth]{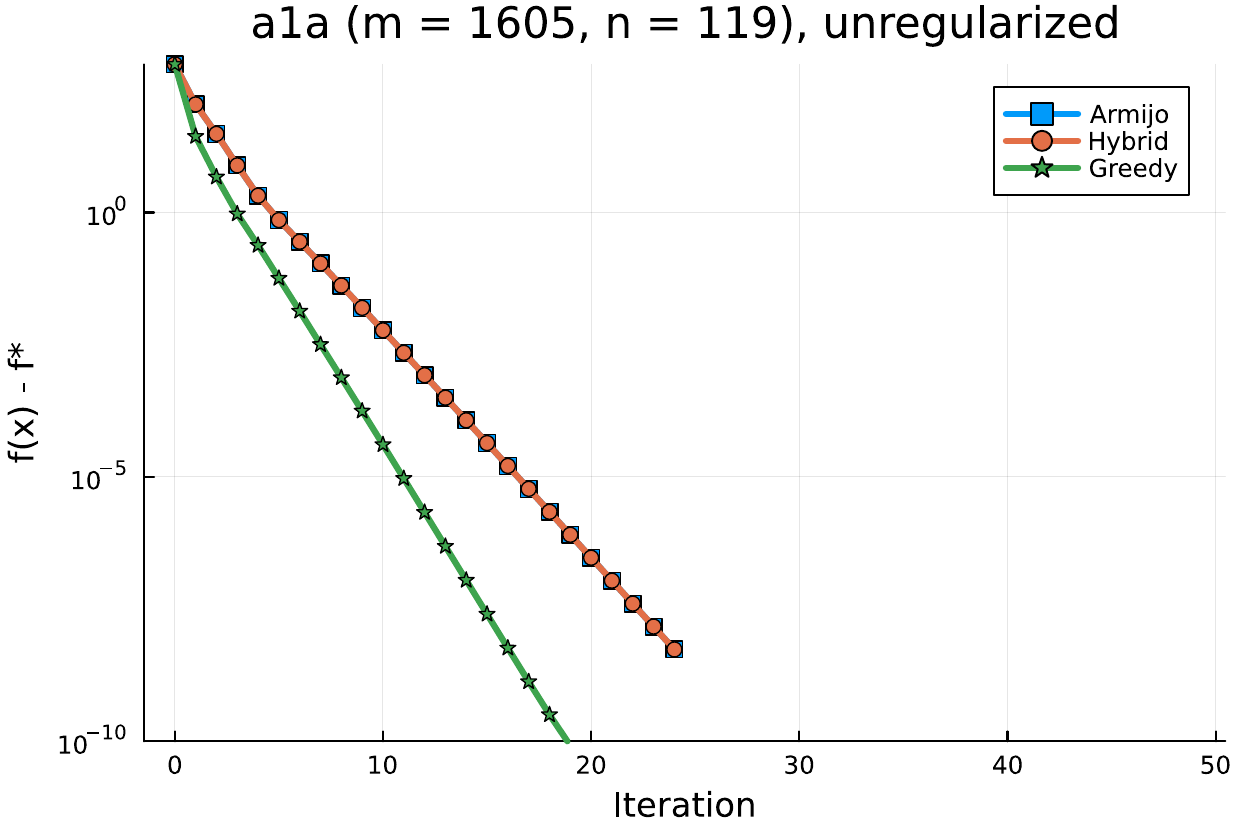}
	\includegraphics[width=.24\textwidth]{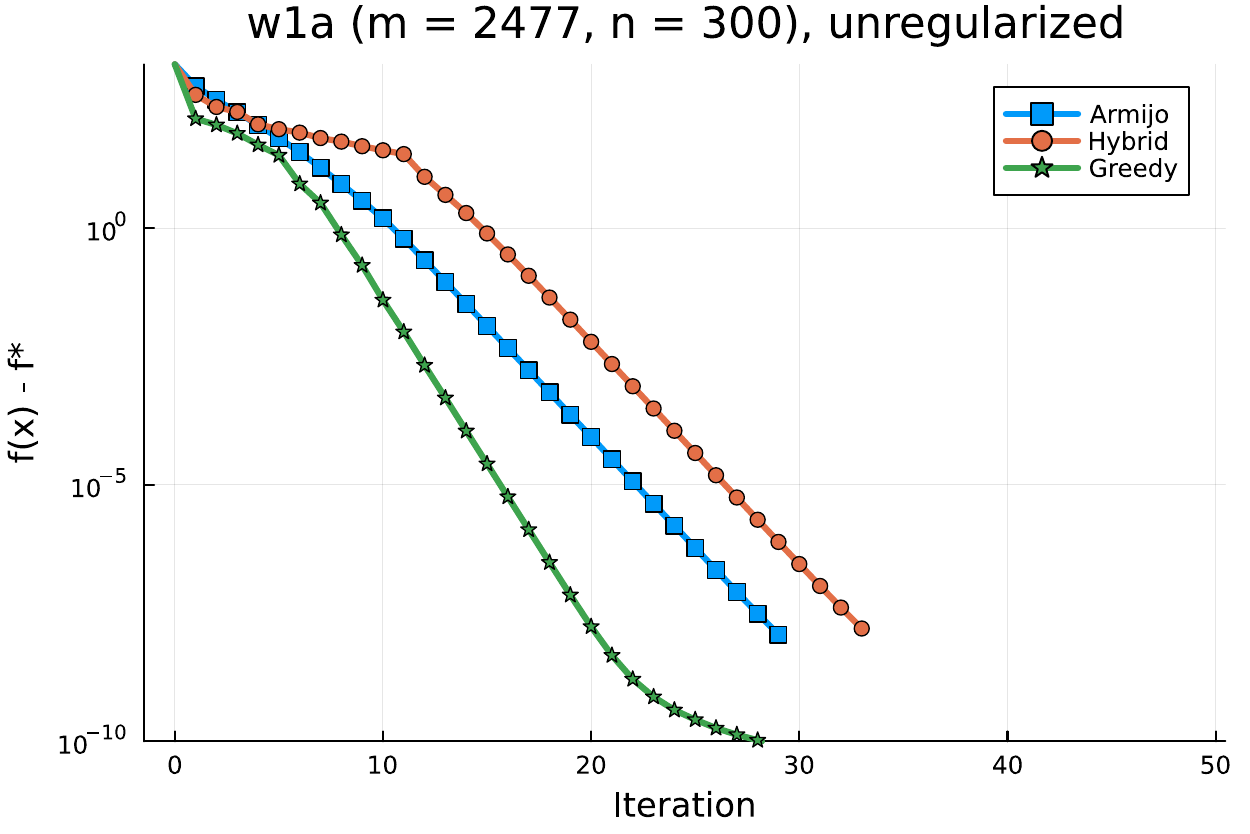}
	\caption{
		Comparison of methods on real logistic regression datasets, in 8 cases where we observed atypical performance in either the regularized or unregularized setting.
	}
	\label{fig:logregRealf2}
\end{figure}

\begin{figure}
	\includegraphics[width=.24\textwidth]{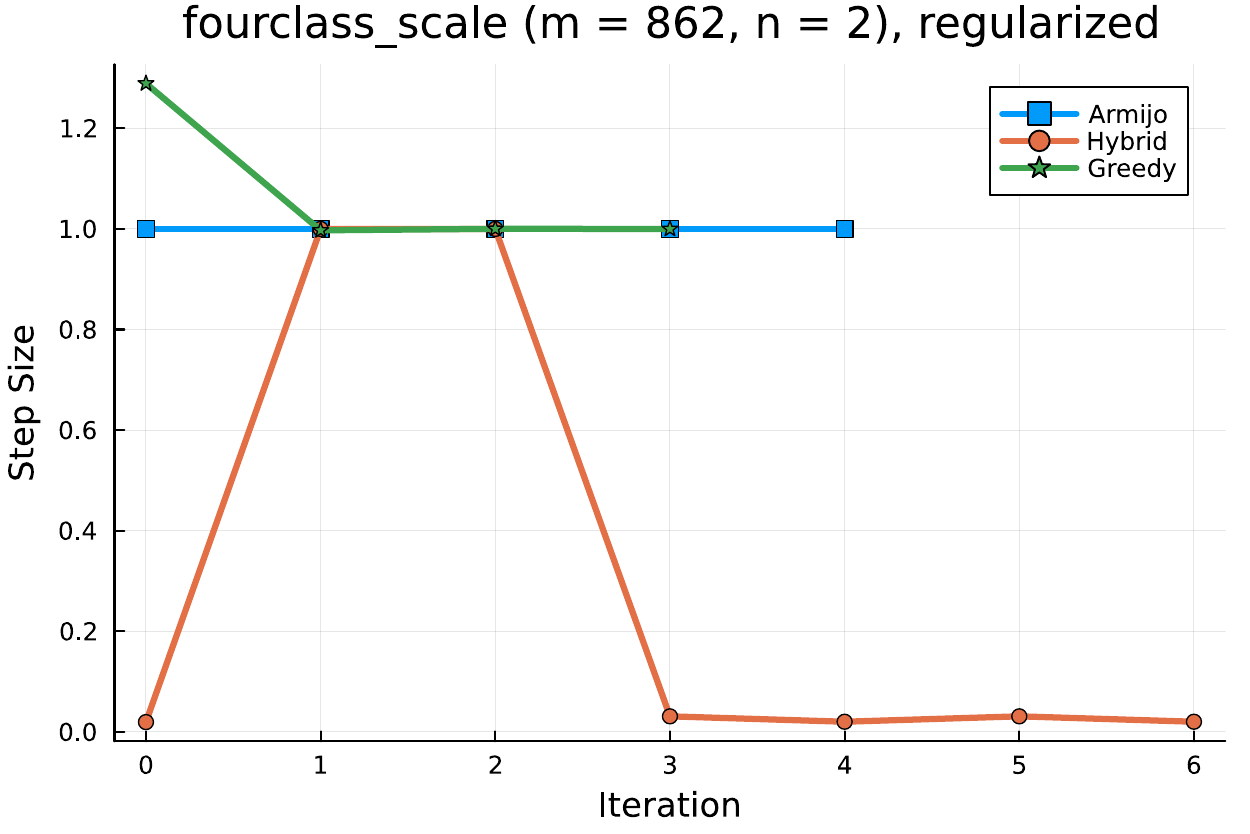}
	\includegraphics[width=.24\textwidth]{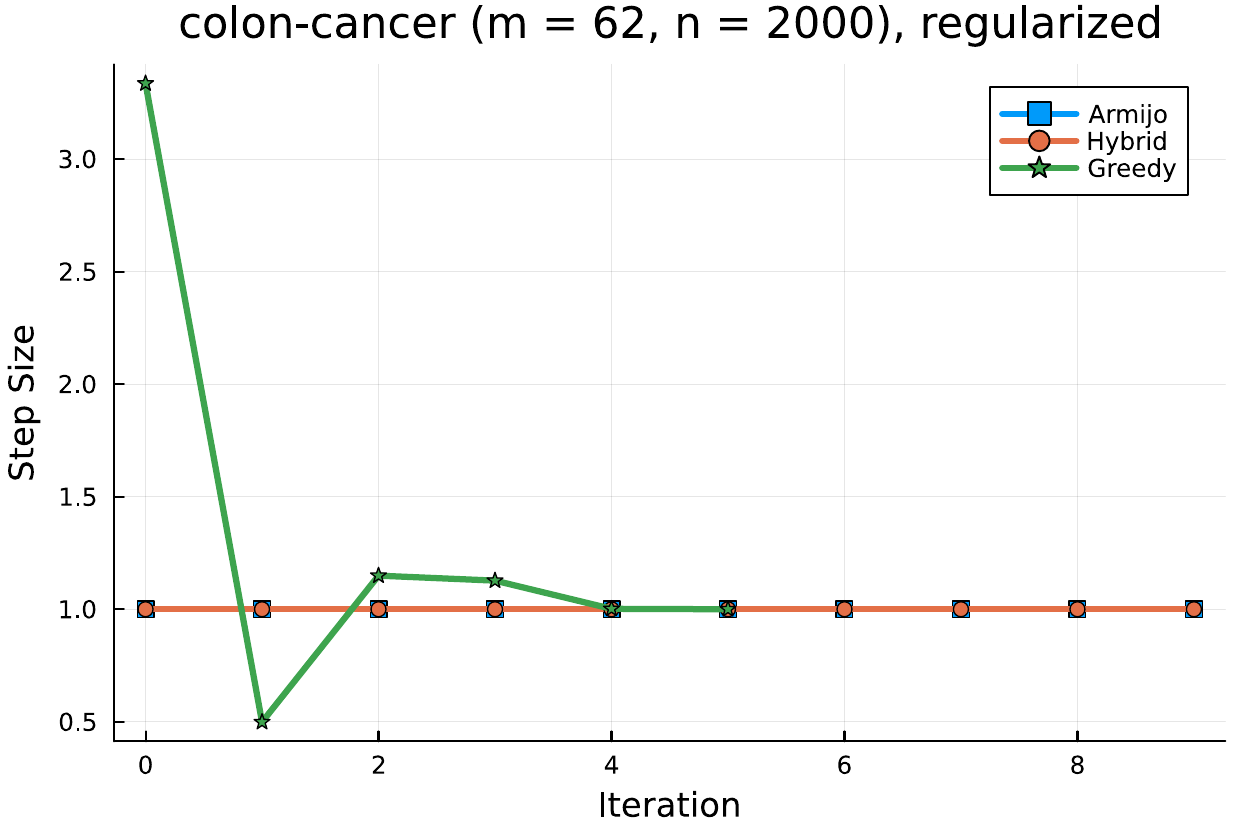}
	\includegraphics[width=.24\textwidth]{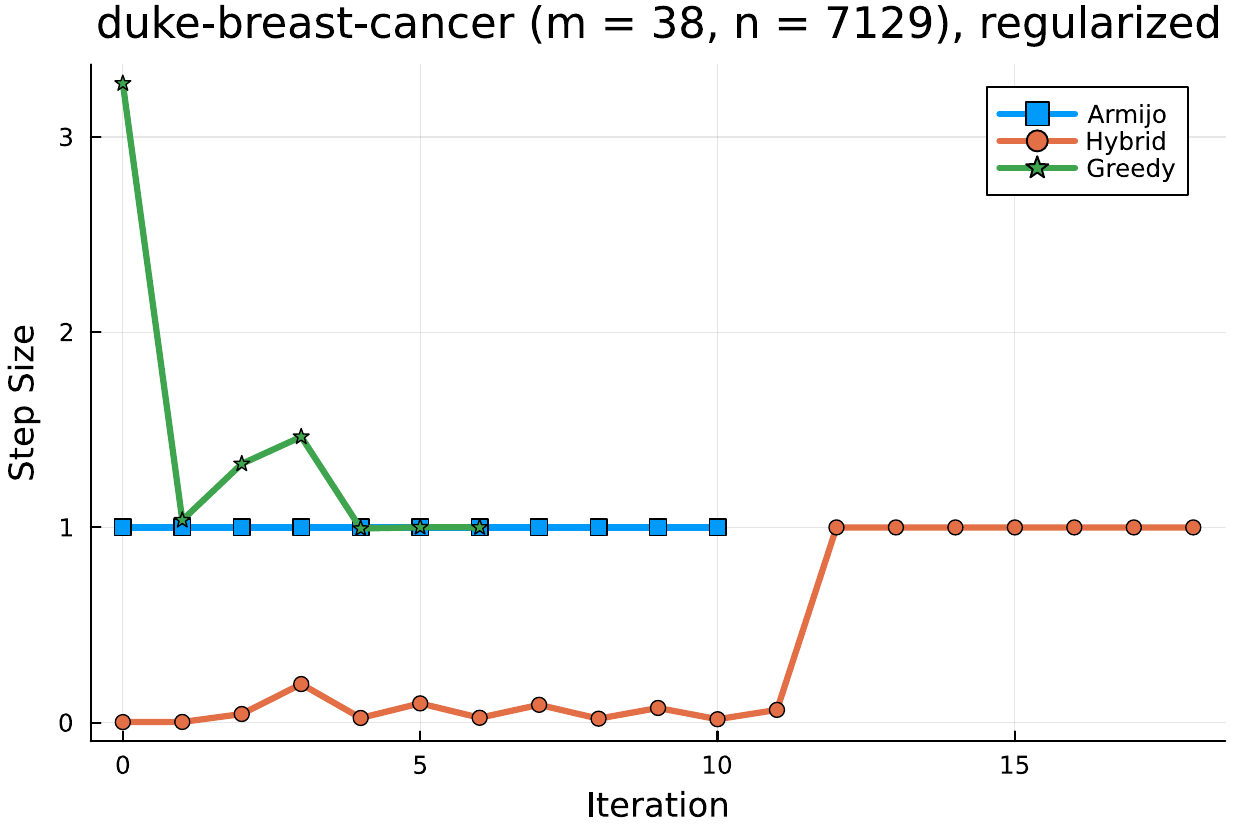}
	\includegraphics[width=.24\textwidth]{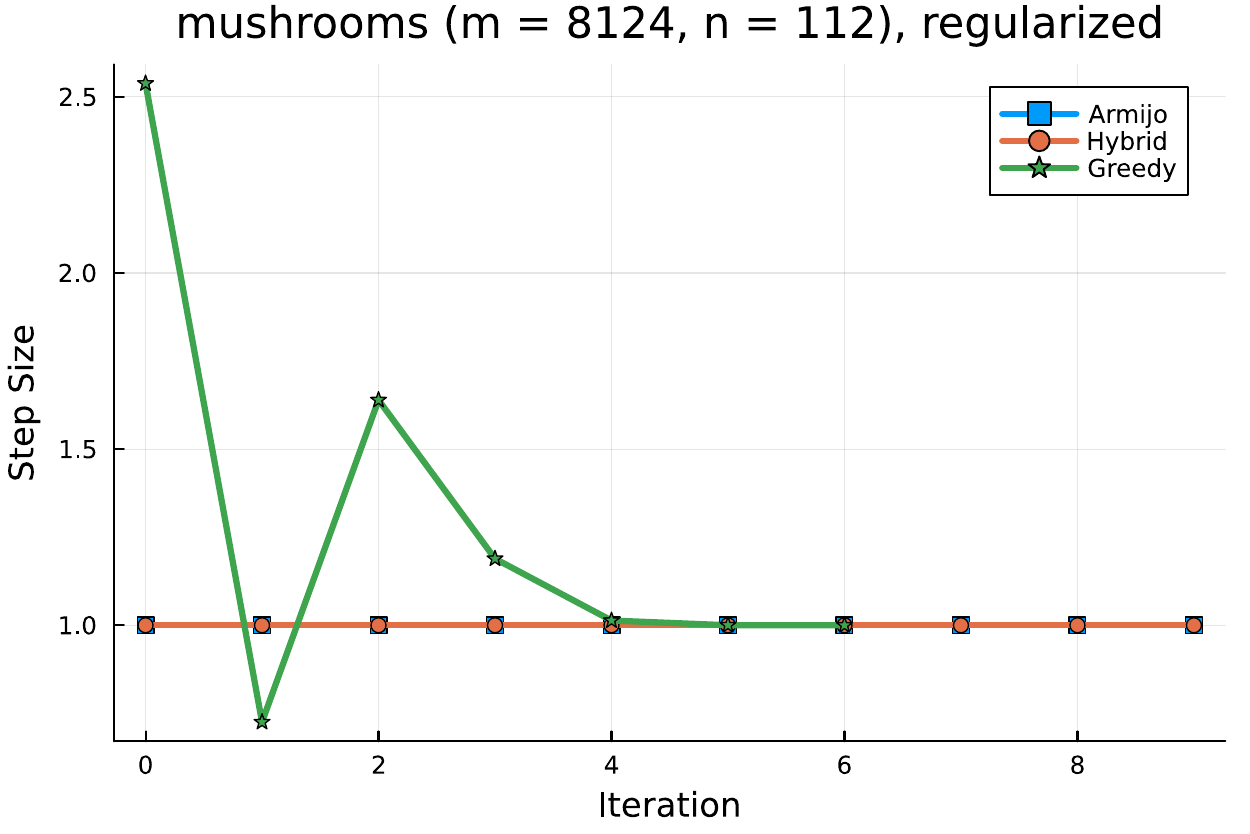}
	\includegraphics[width=.24\textwidth]{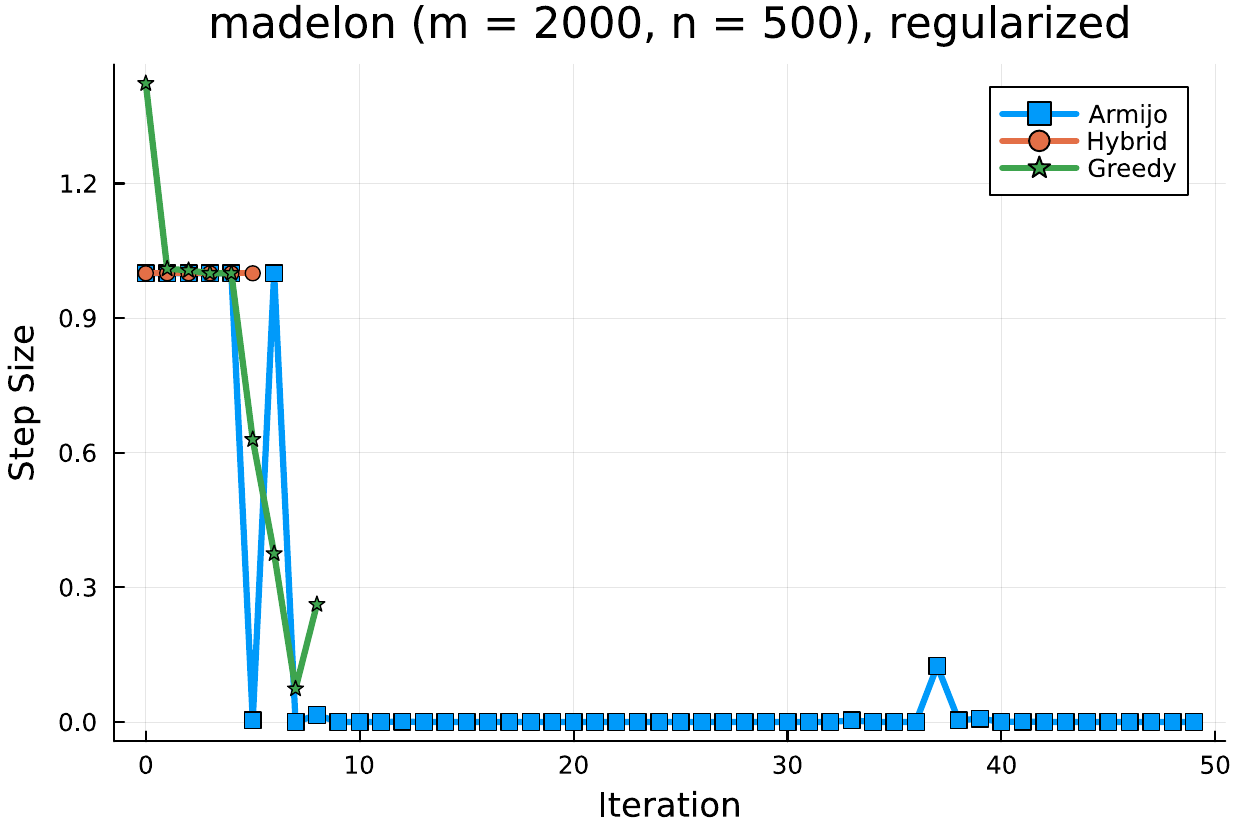}
	\includegraphics[width=.24\textwidth]{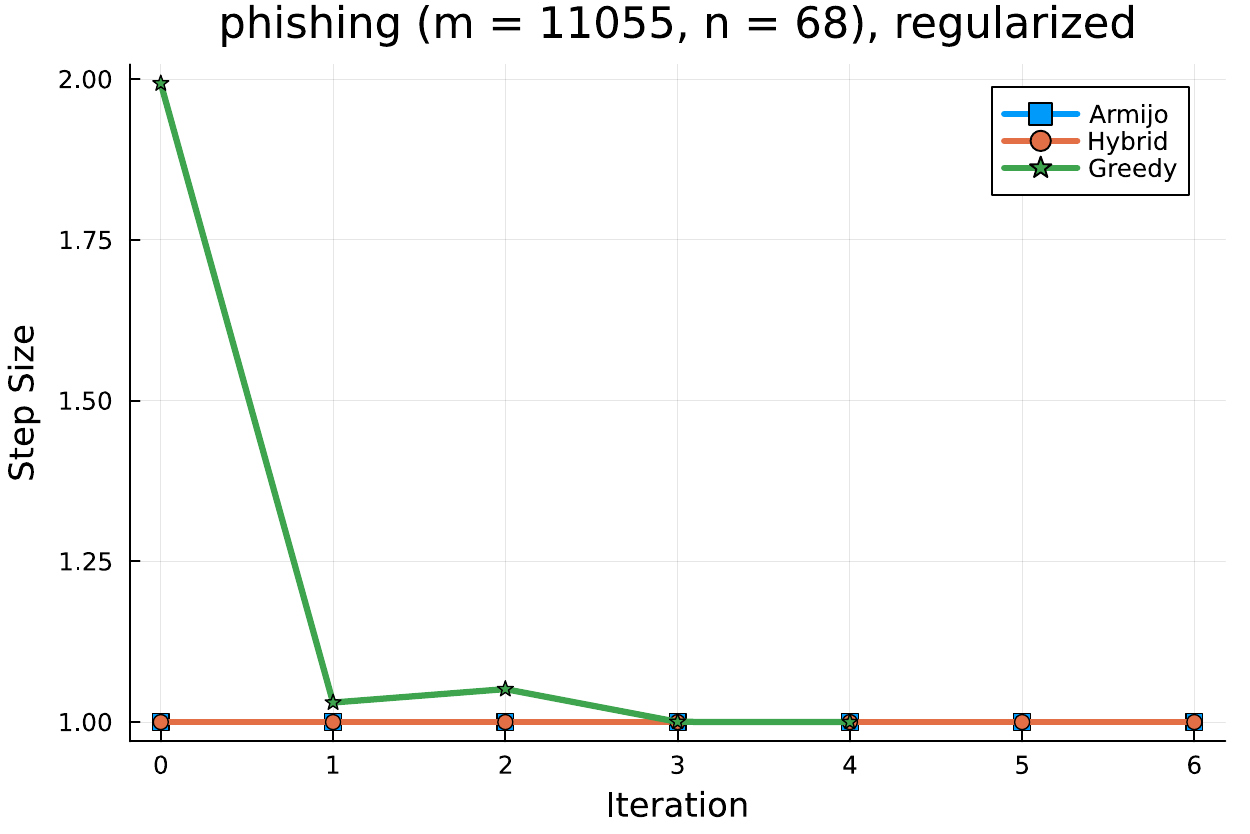}
	\includegraphics[width=.24\textwidth]{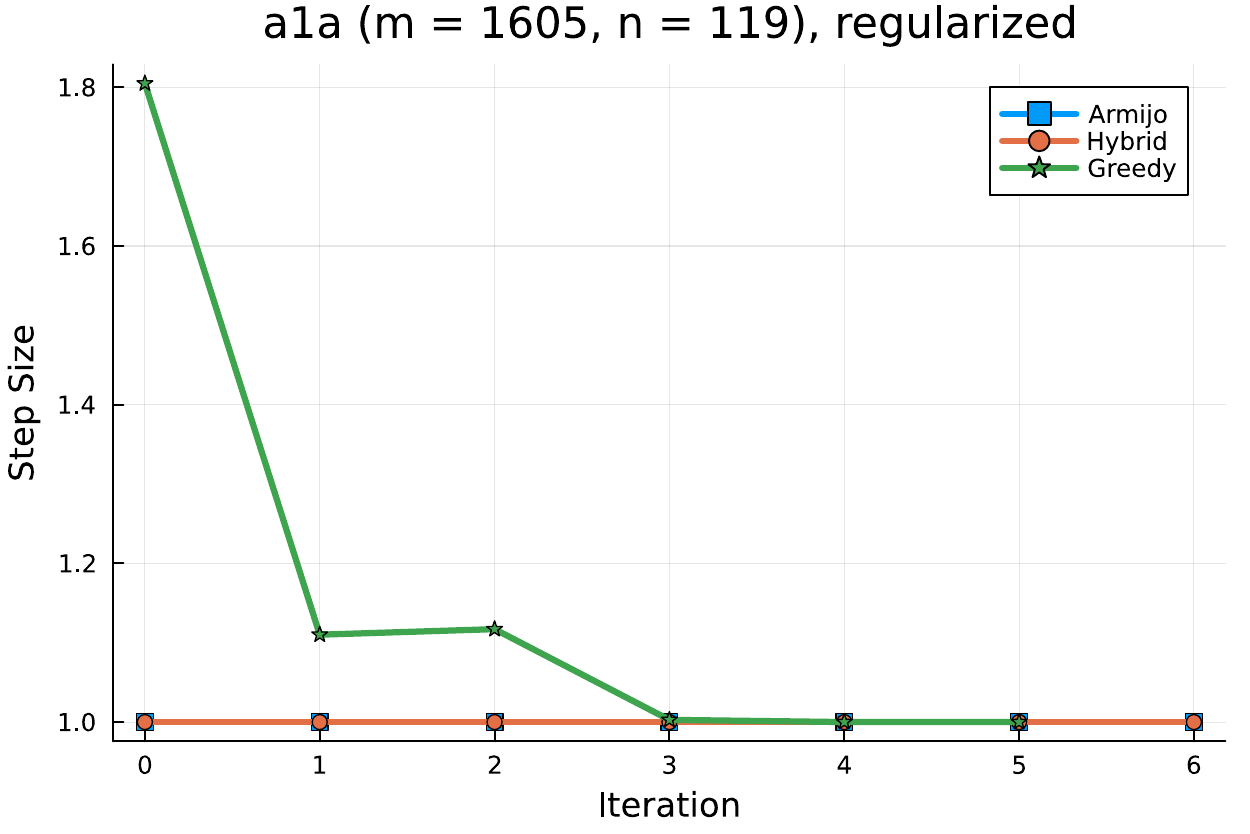}
	\includegraphics[width=.24\textwidth]{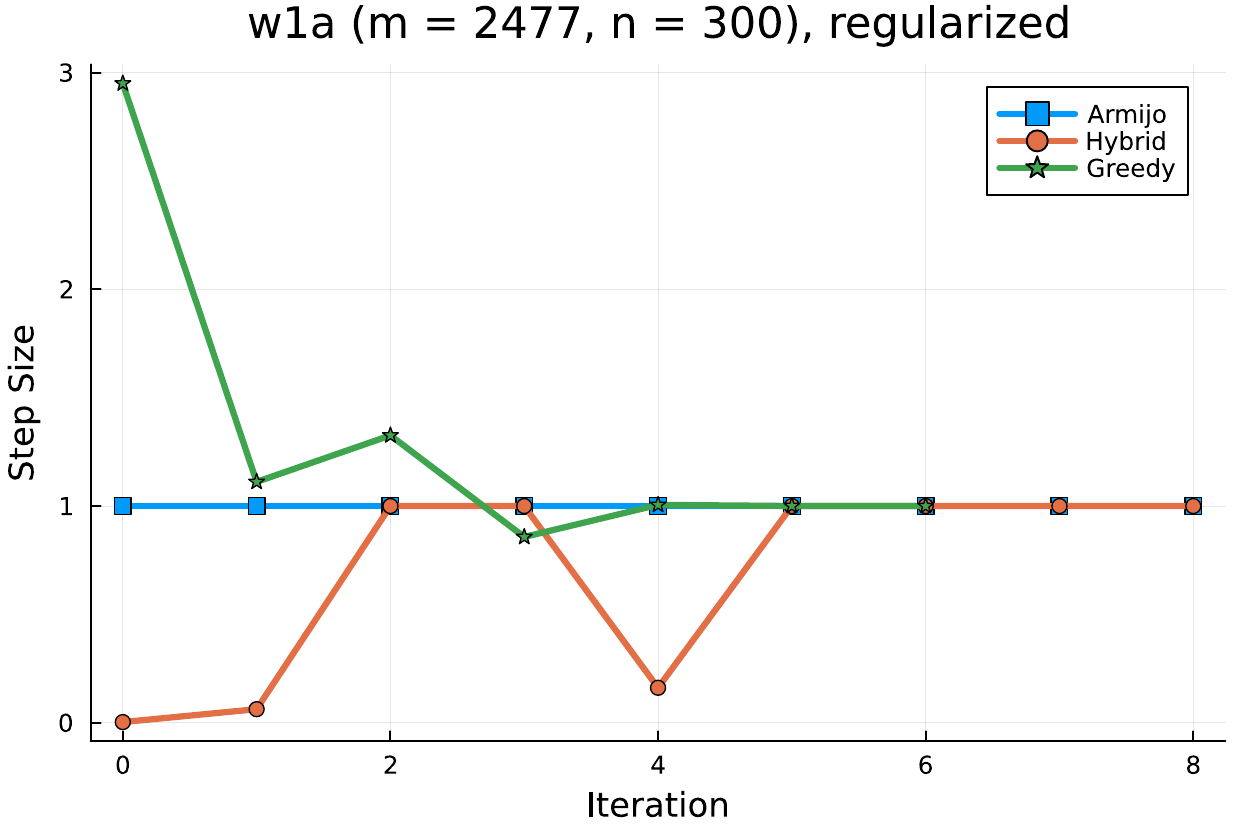}
	\includegraphics[width=.24\textwidth]{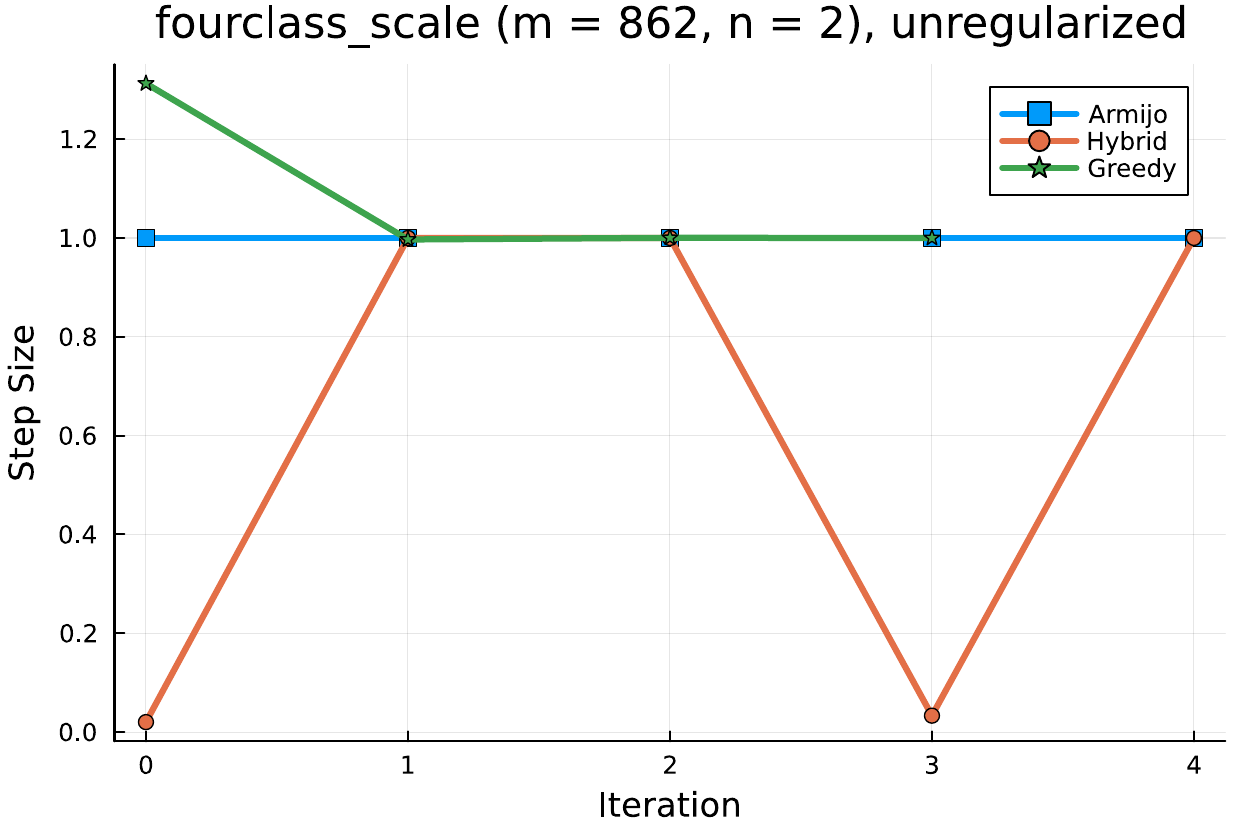}
	\includegraphics[width=.24\textwidth]{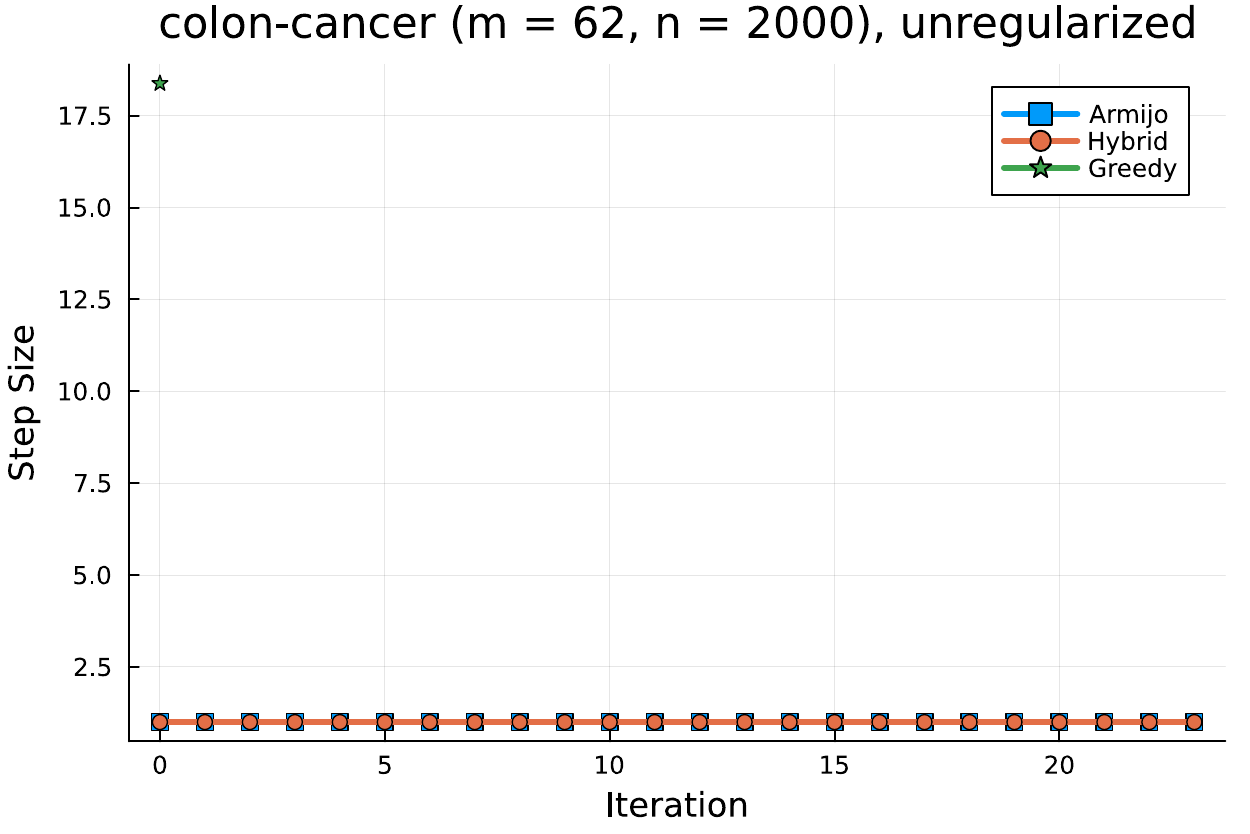}
	\includegraphics[width=.24\textwidth]{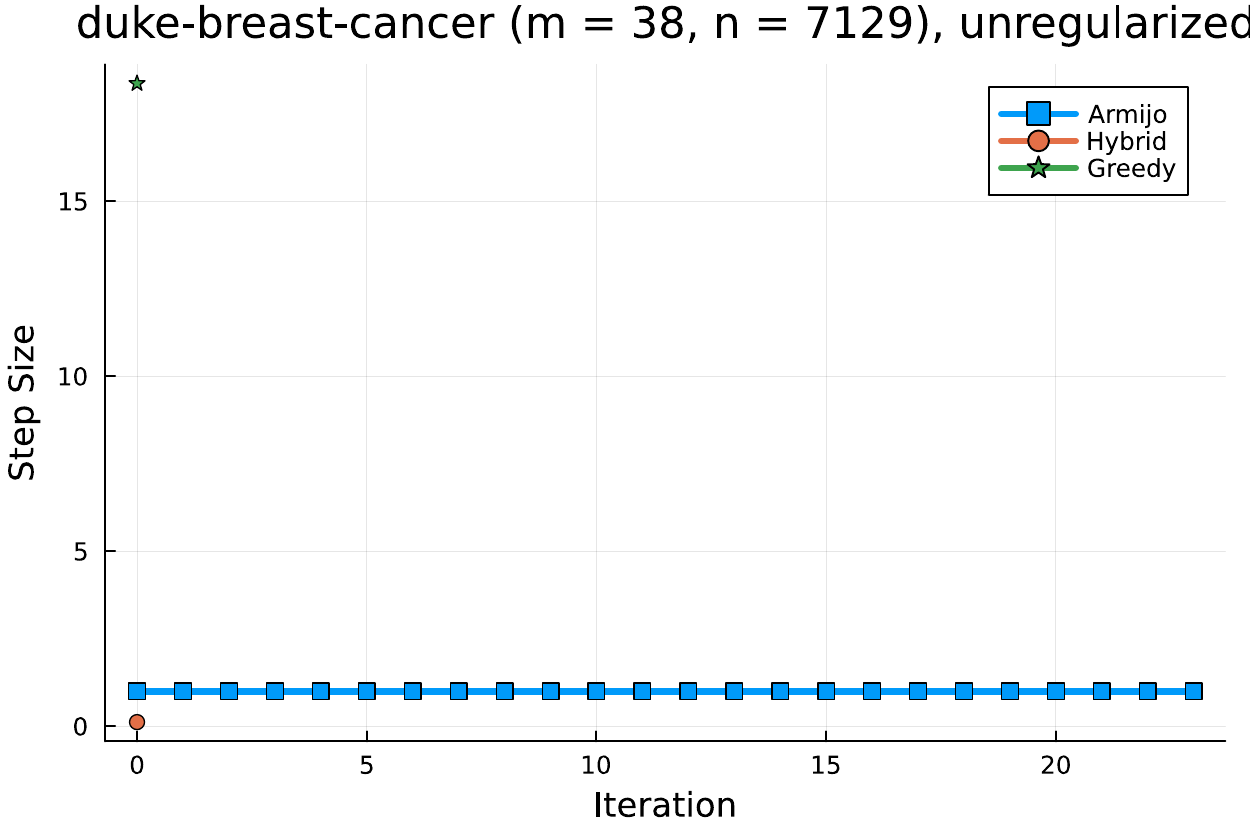}
	\includegraphics[width=.24\textwidth]{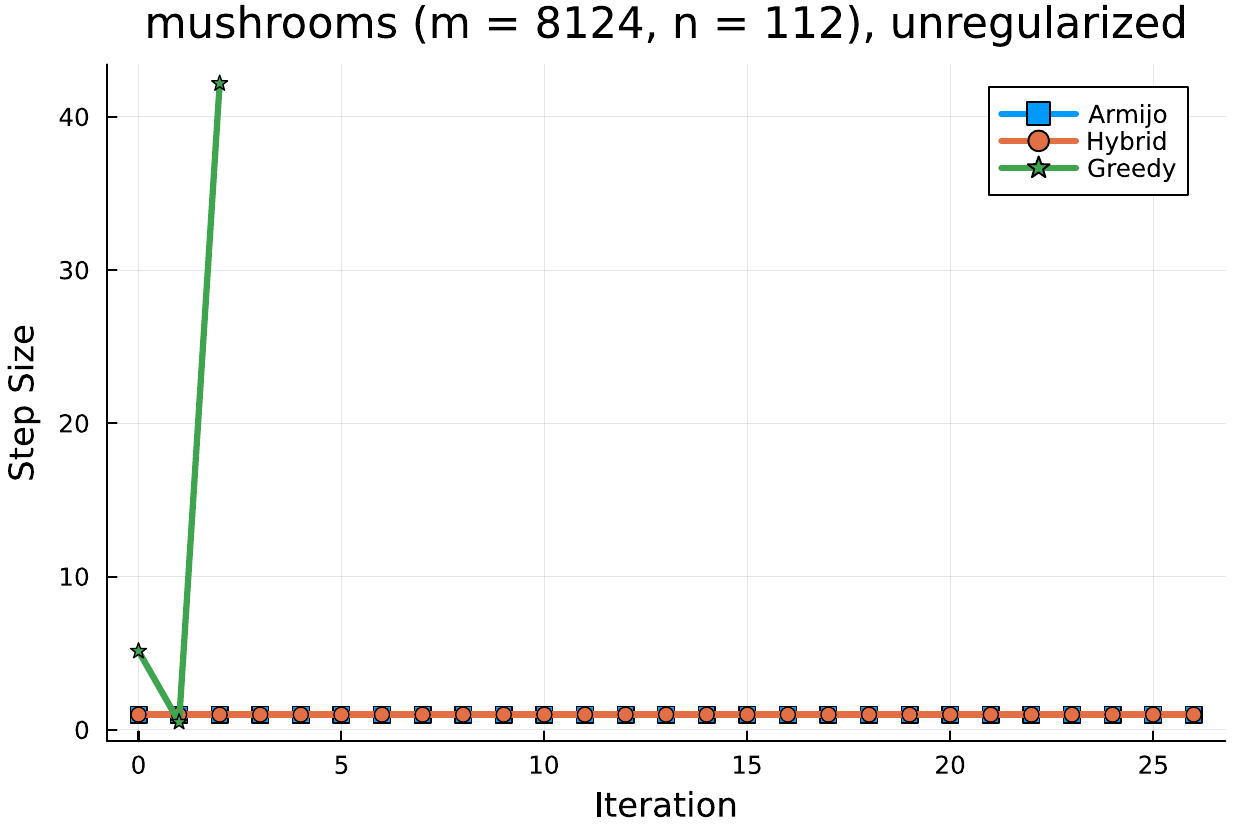}
	\includegraphics[width=.24\textwidth]{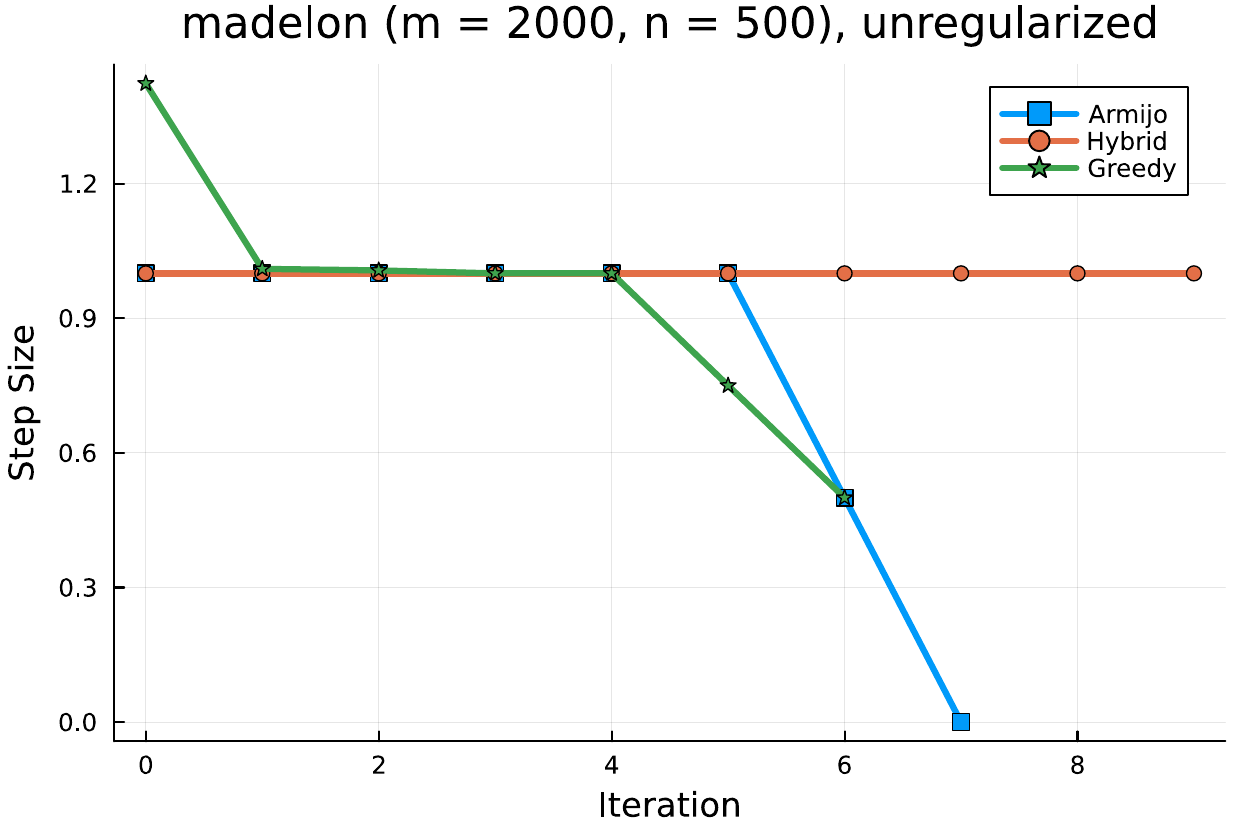}
	\includegraphics[width=.24\textwidth]{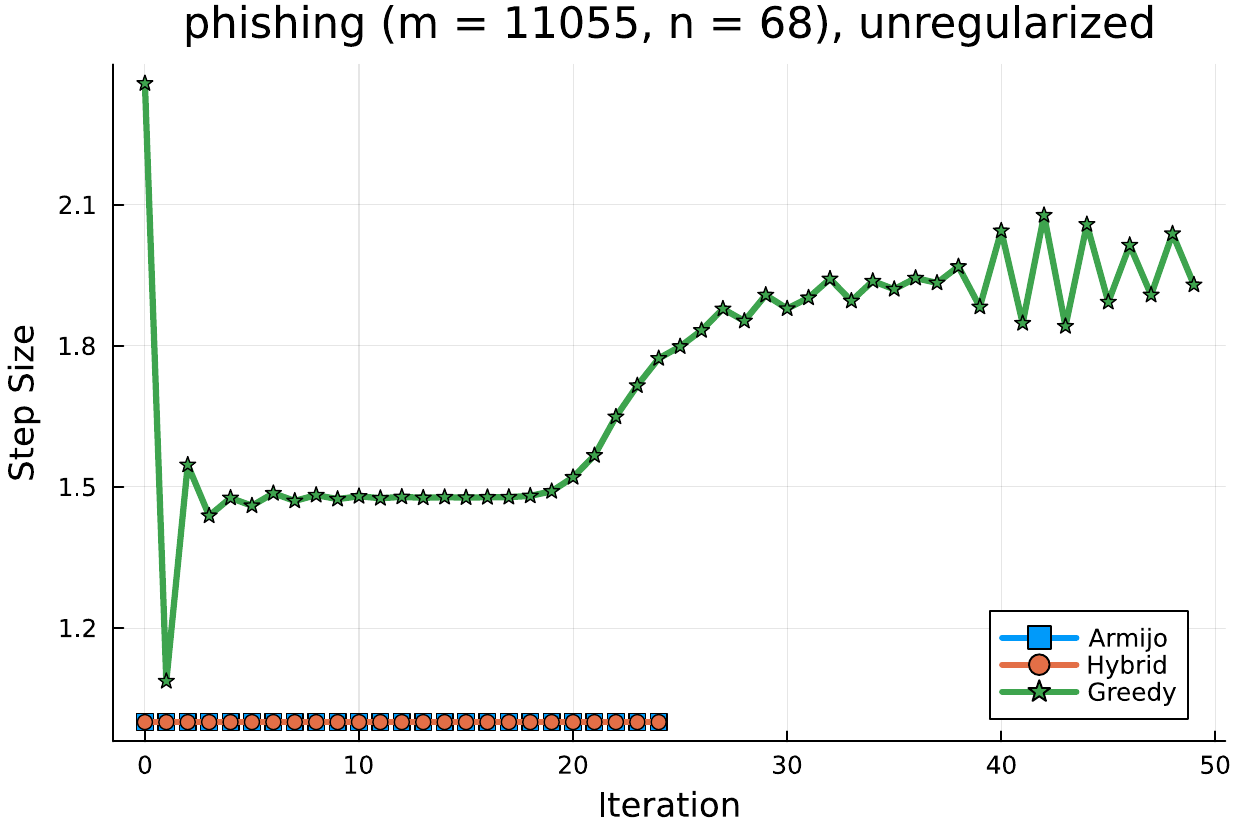}
	\includegraphics[width=.24\textwidth]{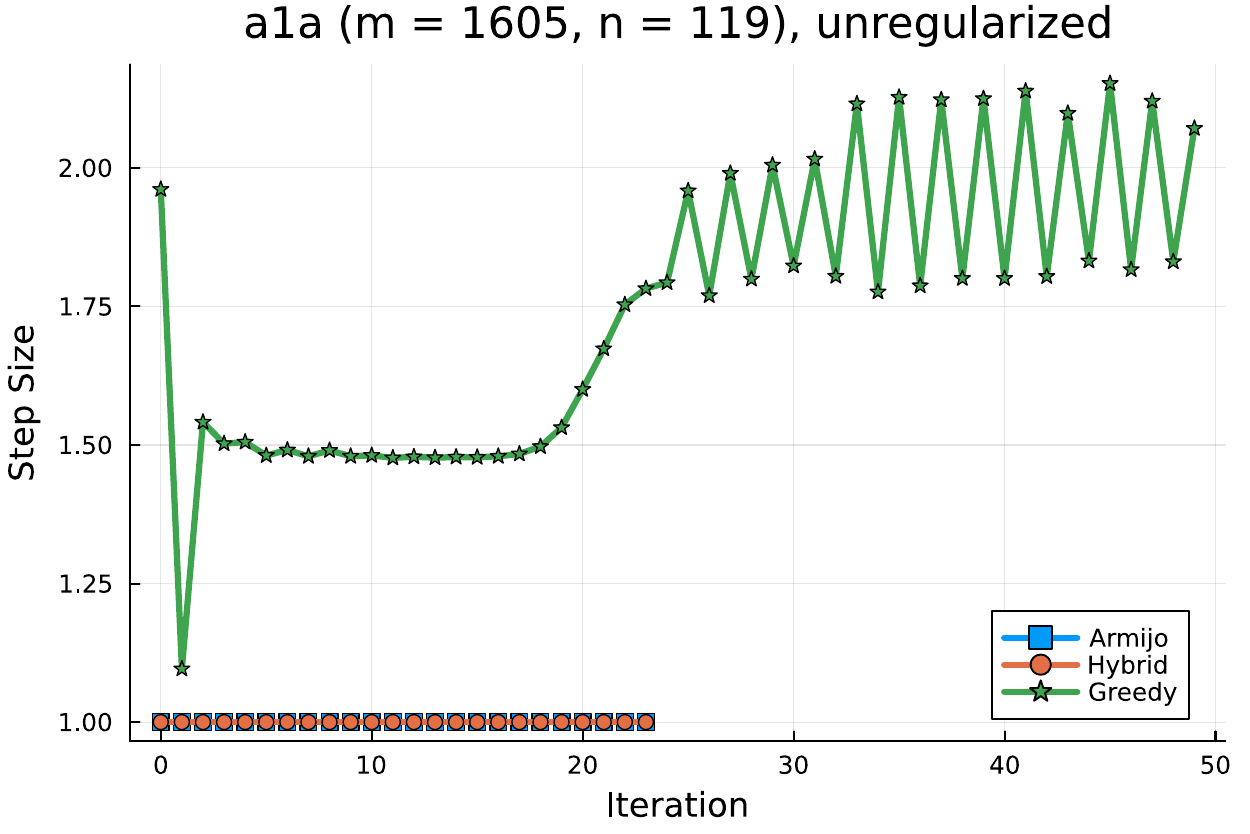}
	\includegraphics[width=.24\textwidth]{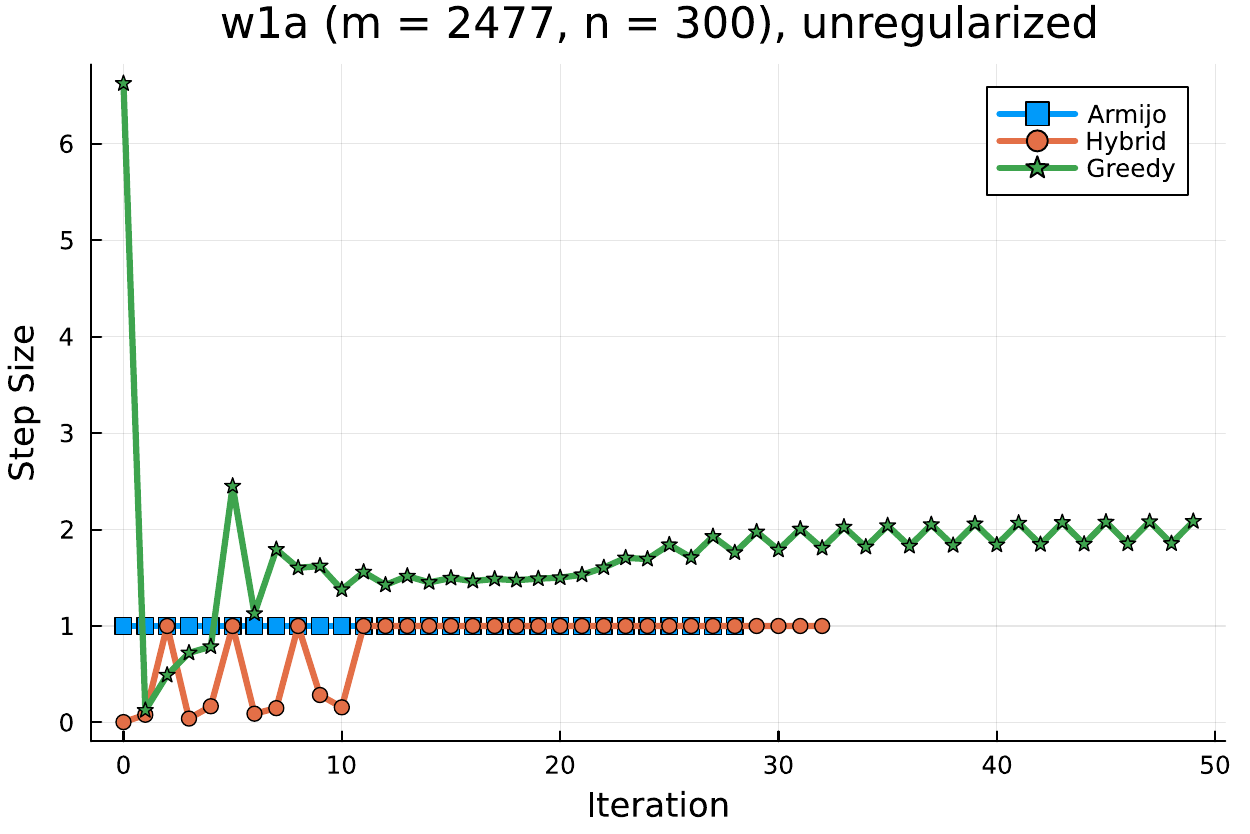}
	\caption{
		Step sizes of methods on real logistic regression datasets, in 8 cases where we observed atypical performance in either the regularized or unregularized setting.
	}
	\label{fig:logregRealt2}
\end{figure}

We performed logistic regression experiments on over 40 datasets obtained using the Dataset Downloader software (\url{https://github.com/fKunstner/dataset-downloader}). On the majority of these datasets, we observed the following trends whether we regularized or not:
\begin{itemize}
	\item The Armijo and hybrid methods produced identical iterations (no gradient descent steps were selected).
	\item The GN method outperformed the other methods.
	\item The GN method typically used a large initial step size but the step sizes quickly converged to 1.
\end{itemize}
In Figures~\ref{fig:logregRealf1} and~\ref{fig:logregRealt1} we plot the performance on 8 datasets where we observed this typical case.
On a smaller number of datasets, we observed different behaviours including:
\begin{itemize}
	\item Cases where the performance of the hybrid method was better or worse than the Armijo method.
	\item Cases where all methods converged extremely quickly.
	\item Cases where the hybrid method performed similar to greedy.
	\item Cases where the GN method converges in 1 step (and one case where the hybrid method did this).
	\item Cases where the GN method eventually begins to oscillate between two non-unit step sizes (this seemed to happen for problems with singular Hessians, so is likely due to the particular Hessian modification strategy we used). 
\end{itemize}
In Figures~\ref{fig:logregRealf2} and~\ref{fig:logregRealt2} we plot the performance on 8 datasets where we observed some of these atypical behaviours.
Despite these different behaviours, we note that GN performed the same or better than the other methods across all datasets.

\section{Runtime Experiments}
\label{app:time}

\begin{figure}
	\includegraphics[width=.24\textwidth]{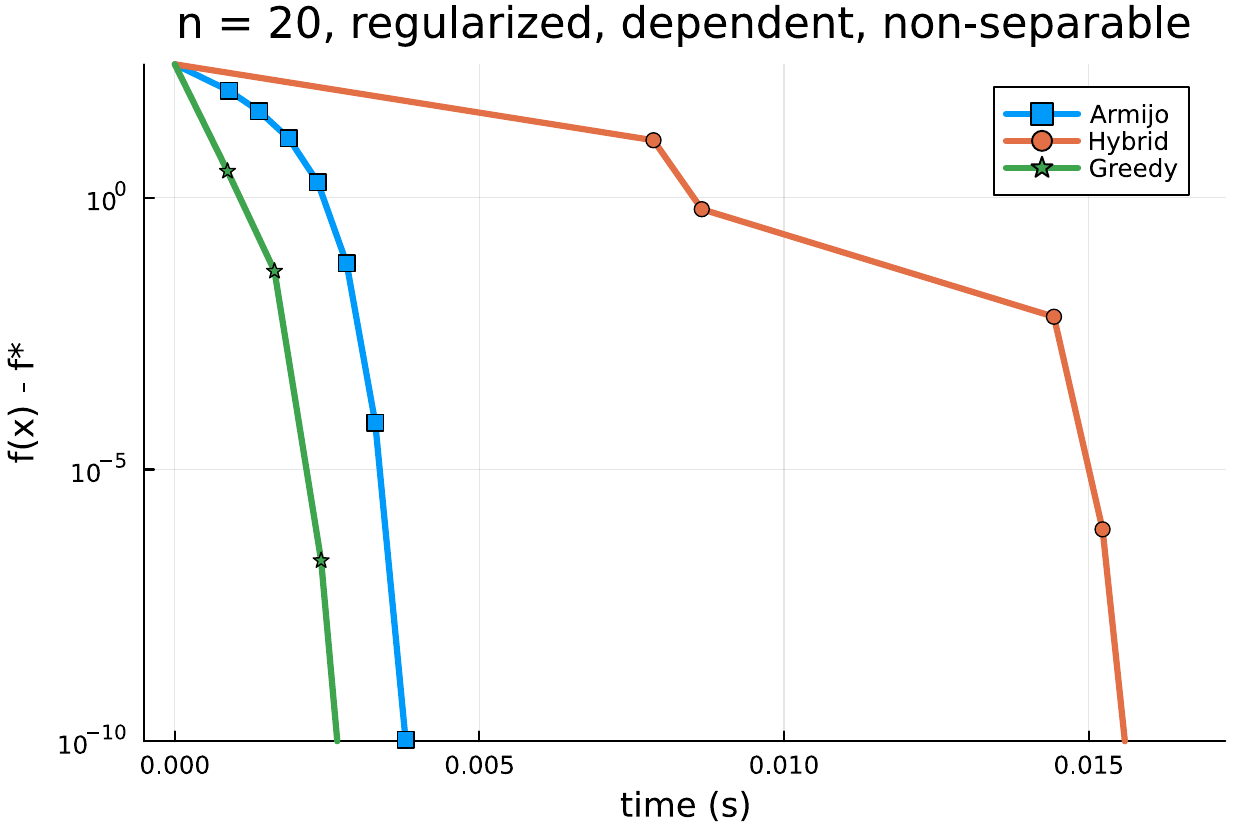}
	\includegraphics[width=.24\textwidth]{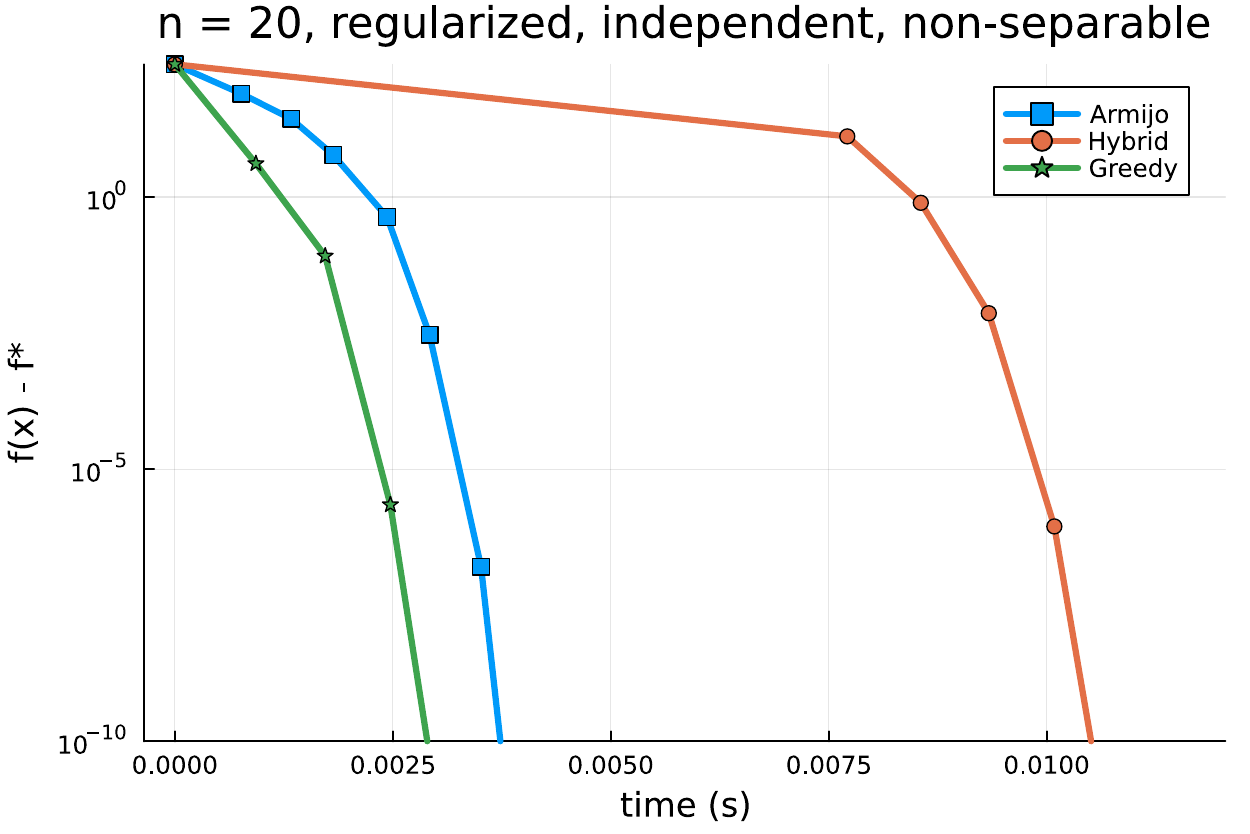}
	\includegraphics[width=.24\textwidth]{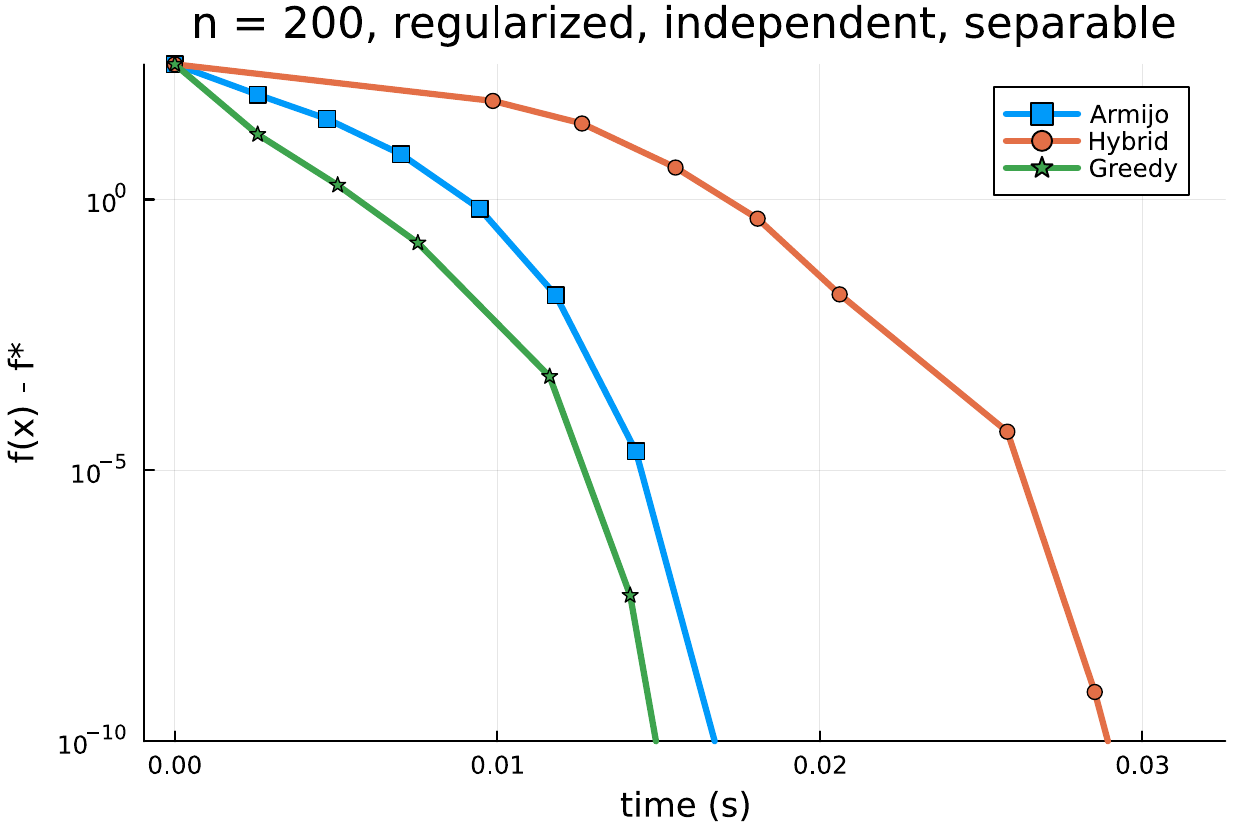}
	\includegraphics[width=.24\textwidth]{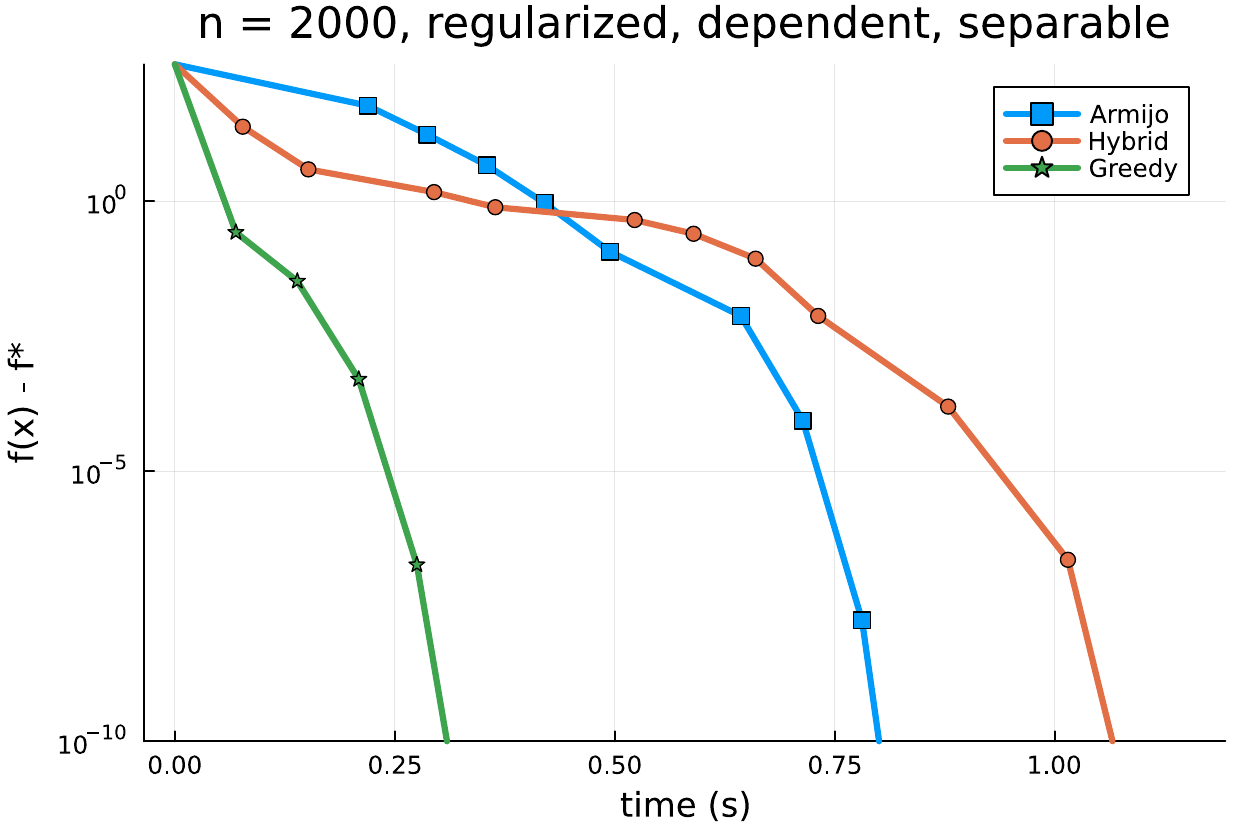}
	\includegraphics[width=.24\textwidth]{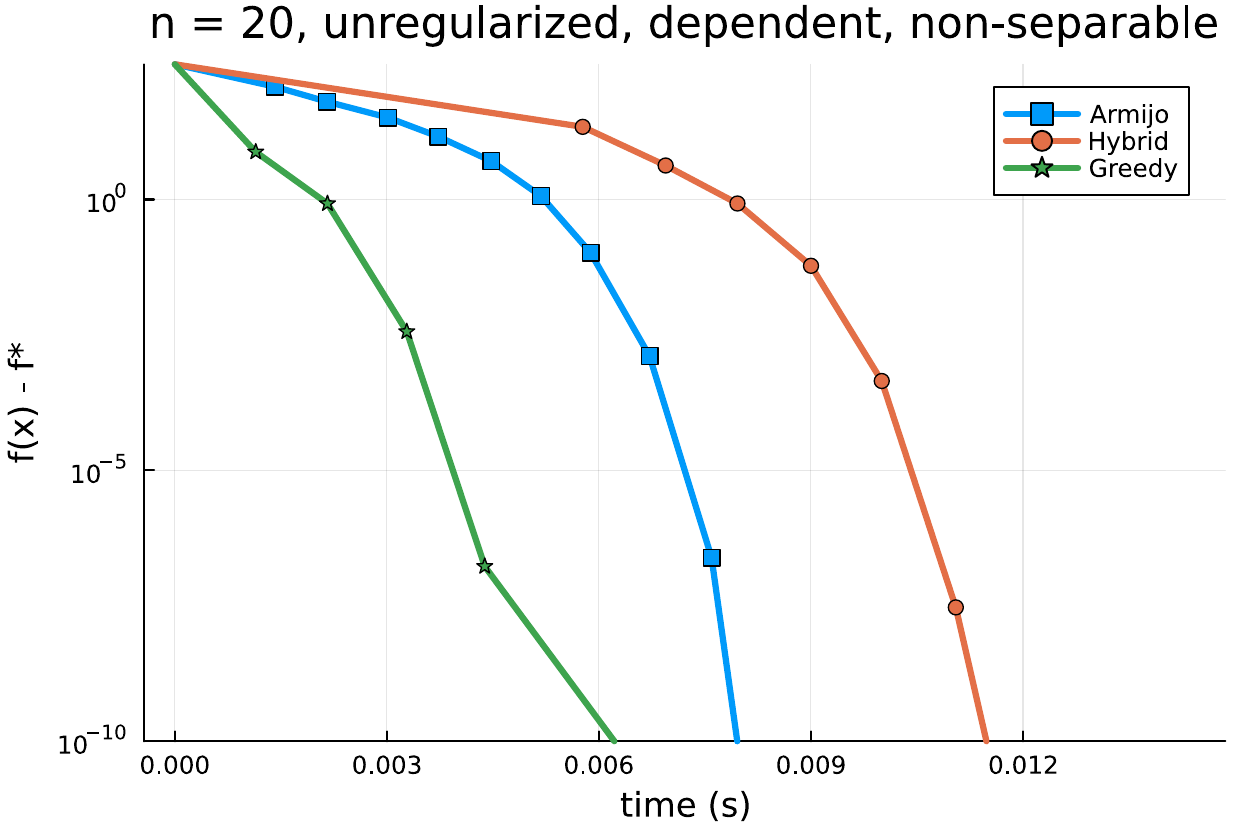}
	\includegraphics[width=.24\textwidth]{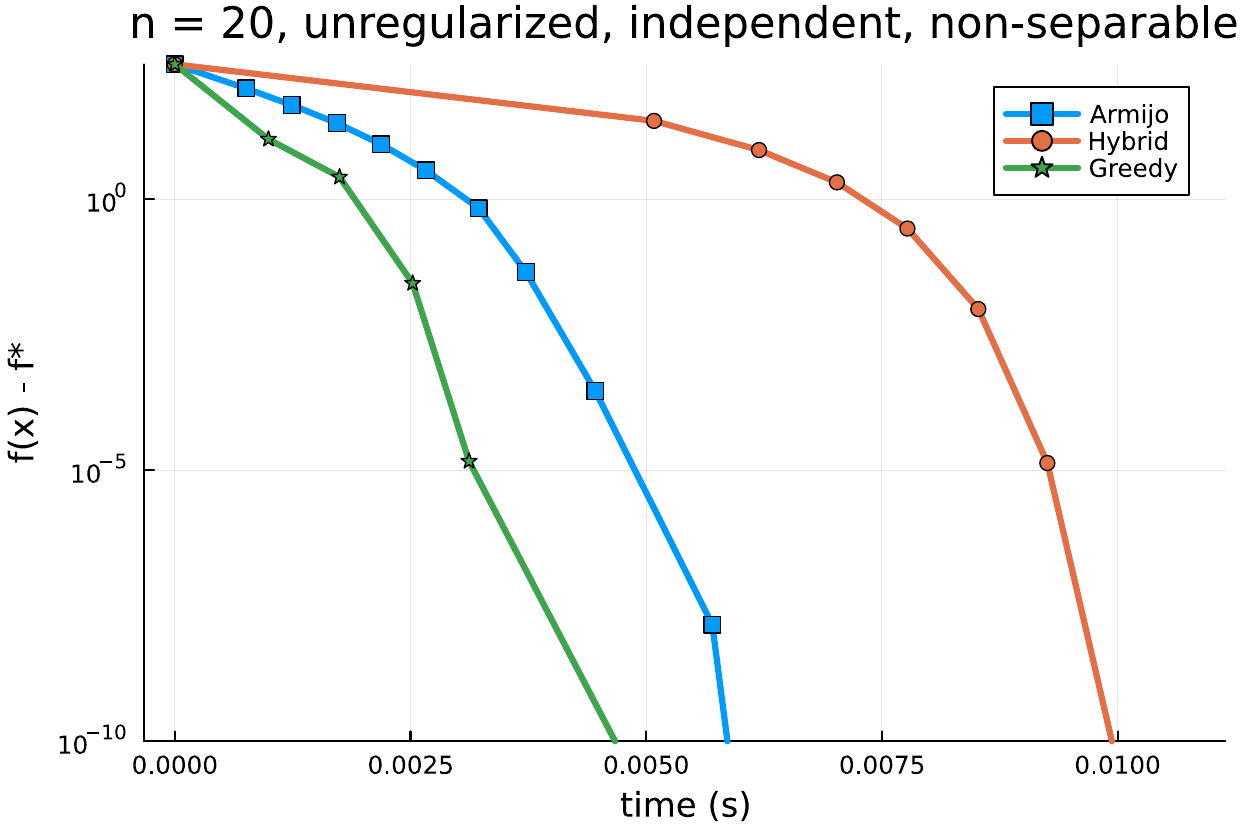}
	\includegraphics[width=.24\textwidth]{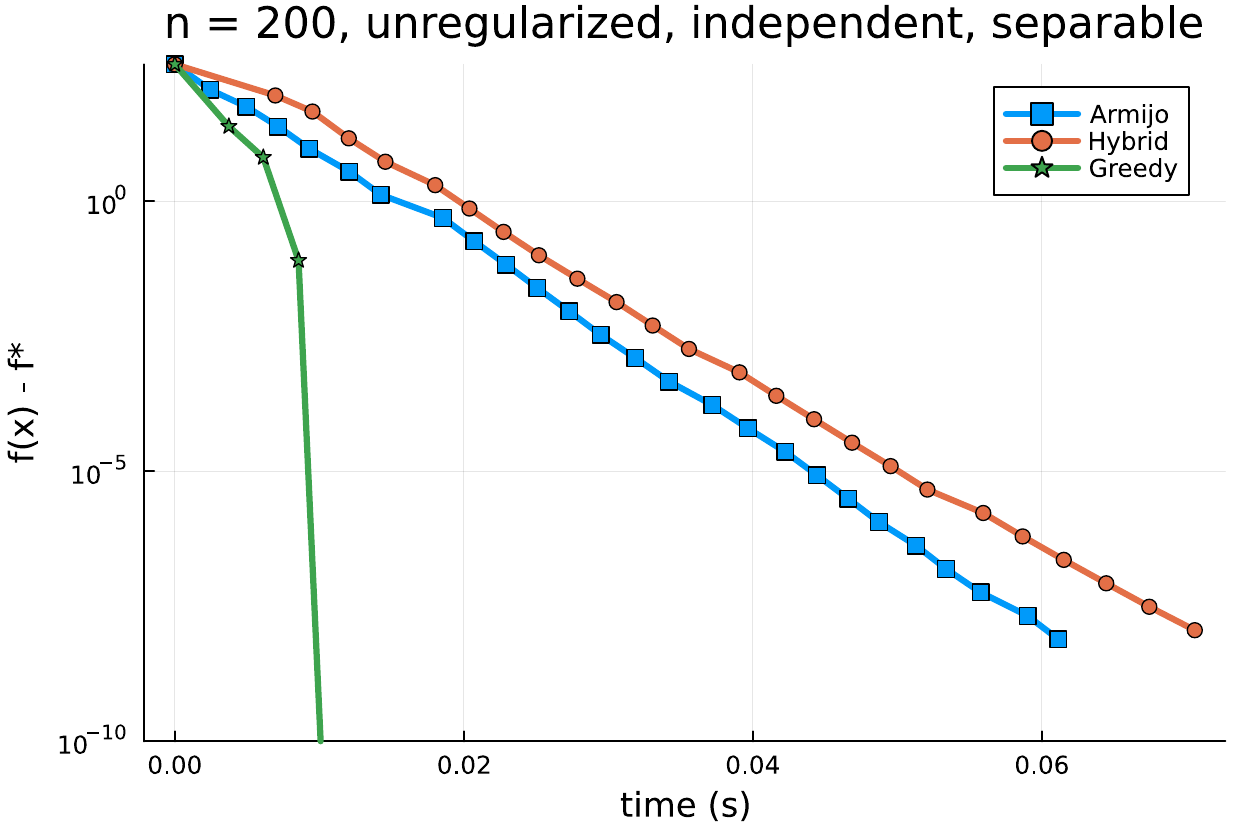}
	\includegraphics[width=.24\textwidth]{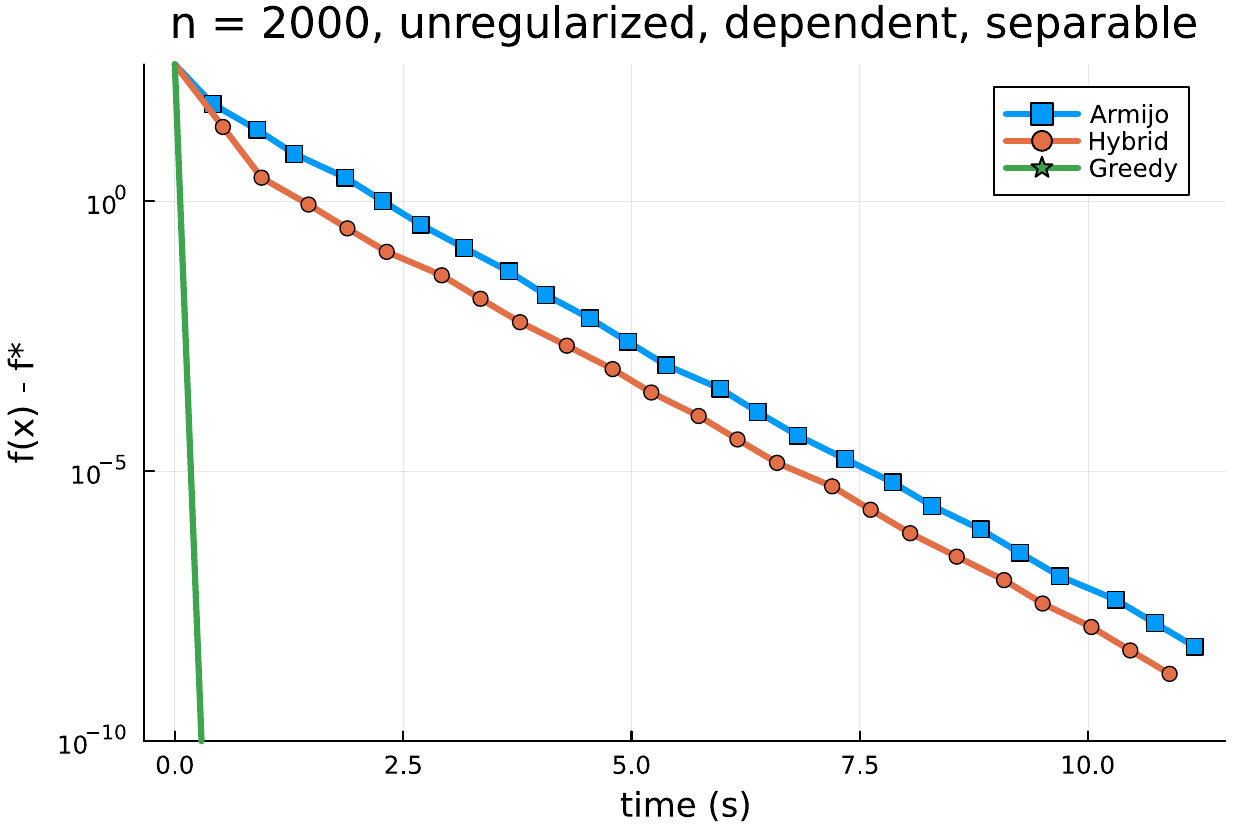}
\caption{Runtime comparison of Newton with Armijo backtracking, hybrid gradient-Newton, and greedy Newton on logistic regression problems with regularization (top row) and without regularization (bottom row).}
	\label{fig:time}
\end{figure}

In Figure~\ref{fig:time} we repeat the synthetic data experiments of Figure~\ref{fig:logregf} but measure performance in terms of runtime. The GN method still outperformed the other methods in terms of runtime. In some cases the hybrid method performed worse in terms of runtime, although it is not clear why (other than requiring an additional matrix-vector product per iteration compared to the other methods). We caution the reader that runtime comparisons are very sensitive to precise implementations, and that for larger datasets the iteration counts of Figure~\ref{fig:logregf} are likely to be a better predictor of performance than runtime performance on specific datasets with specific hardware.

\section{Armijo Backtracking with Larger Initialization}
\label{app:armijo}

\begin{figure}
	\includegraphics[width=.24\textwidth]{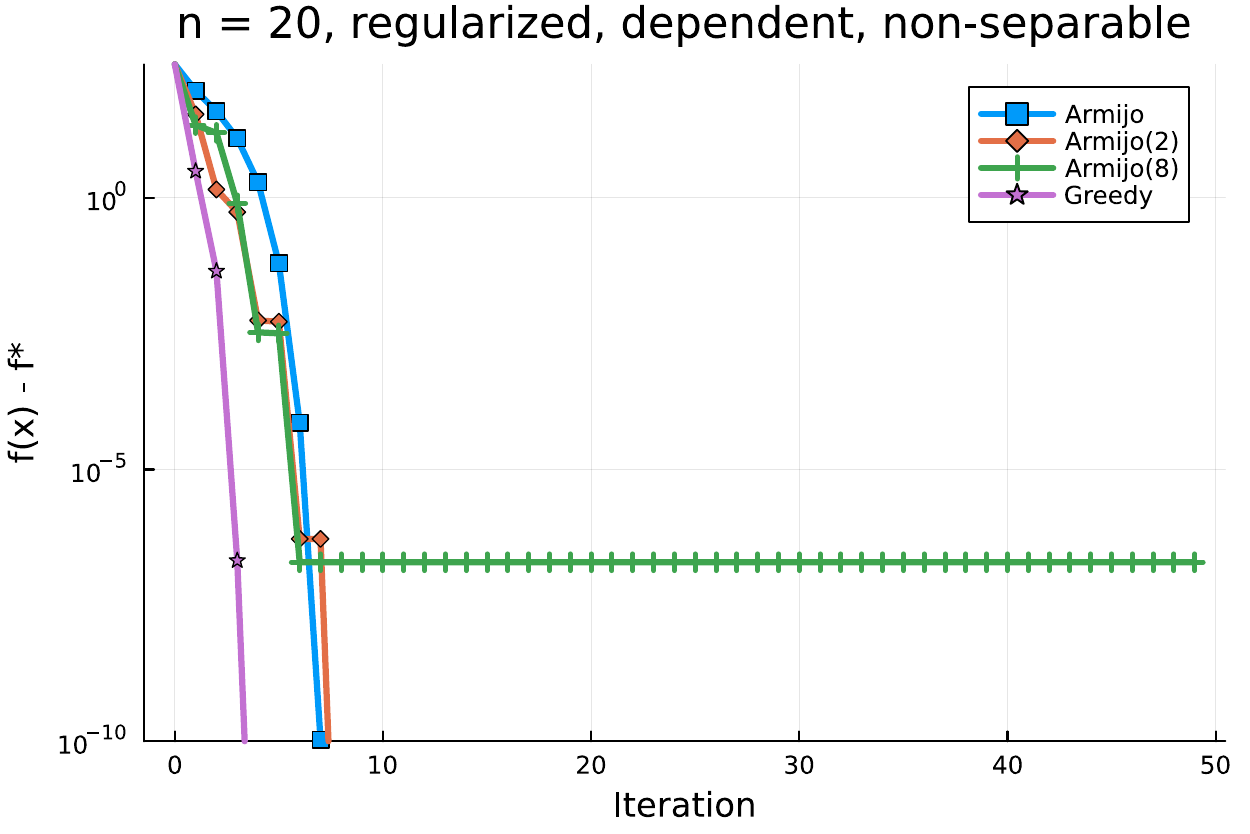}
	\includegraphics[width=.24\textwidth]{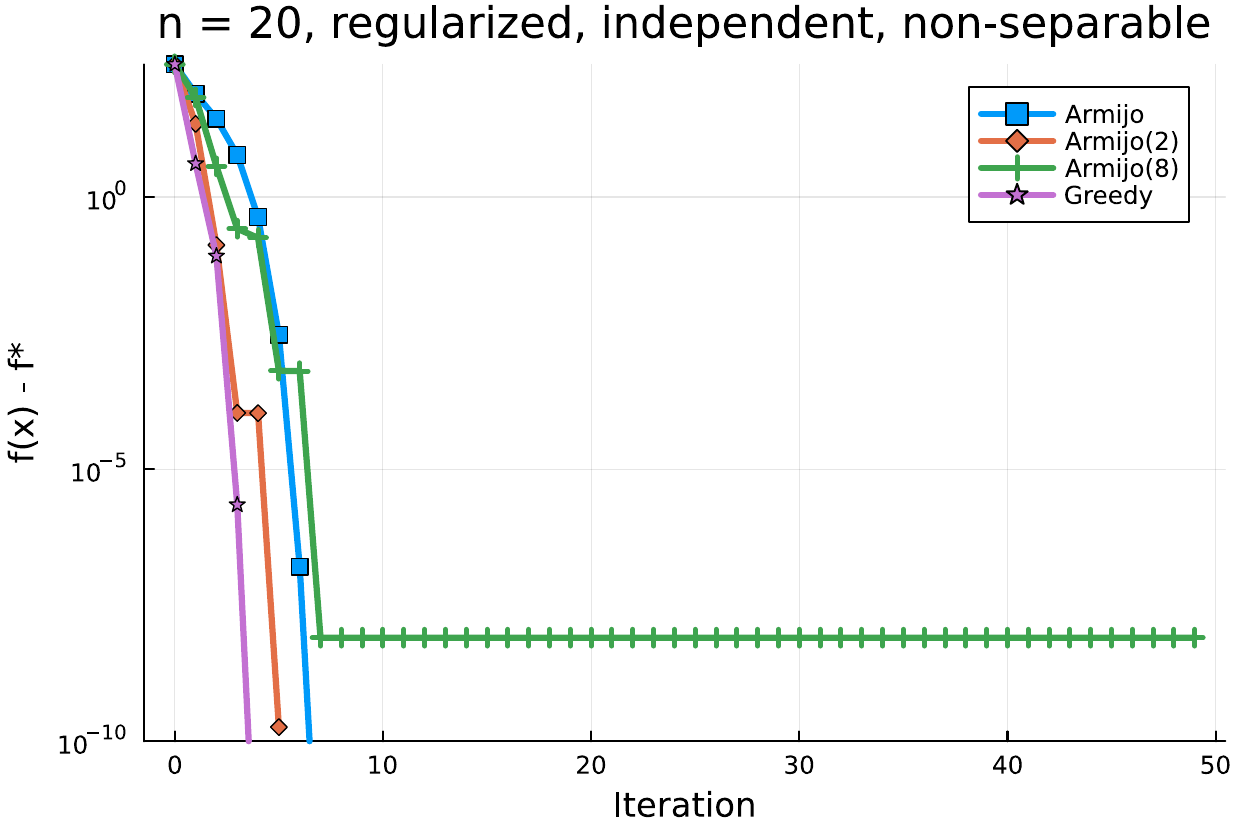}
	\includegraphics[width=.24\textwidth]{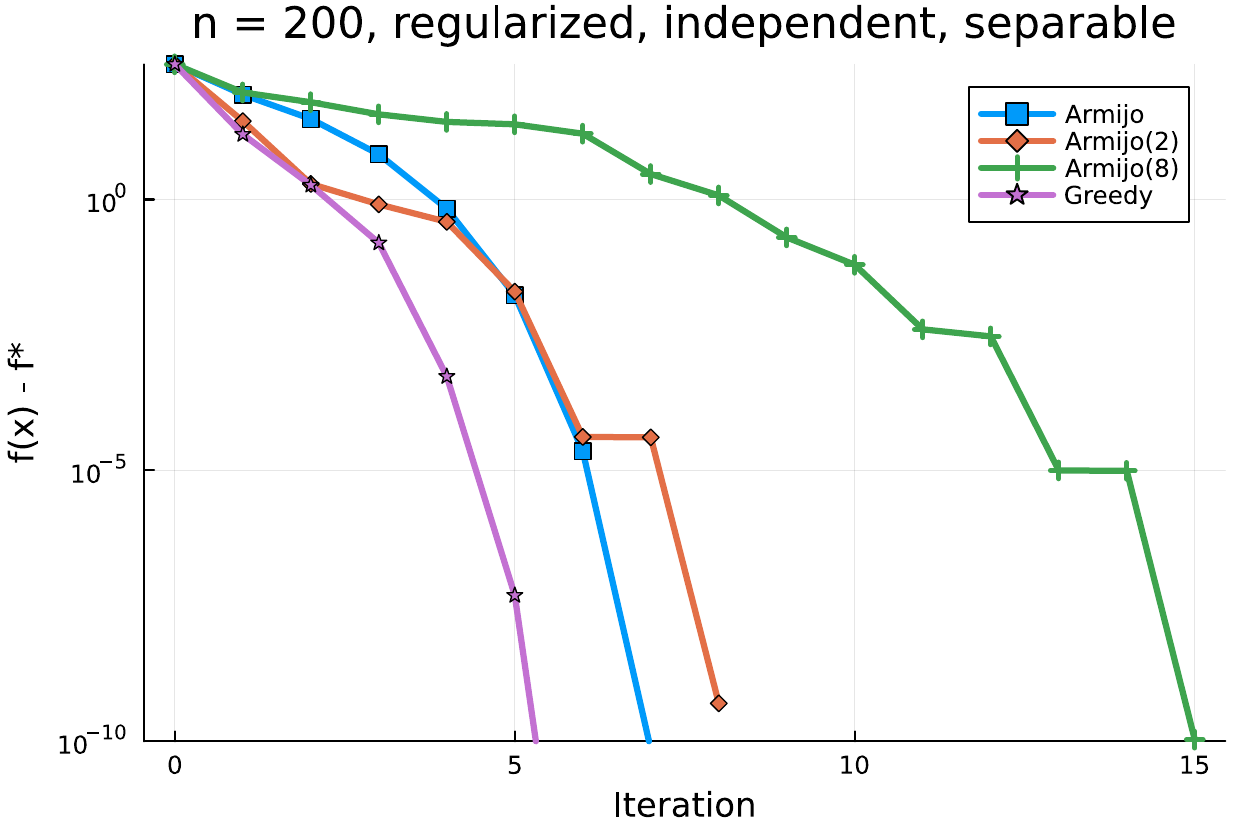}
	\includegraphics[width=.24\textwidth]{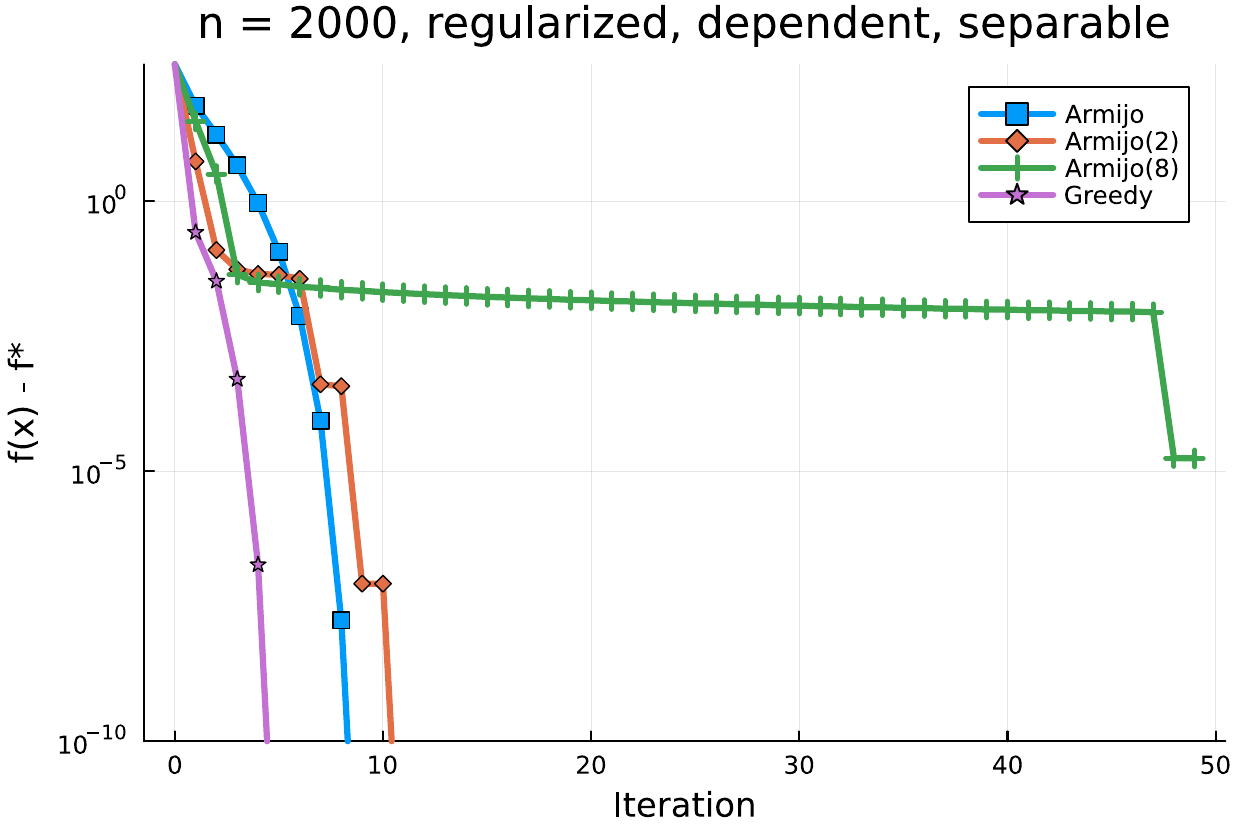}
	\includegraphics[width=.24\textwidth]{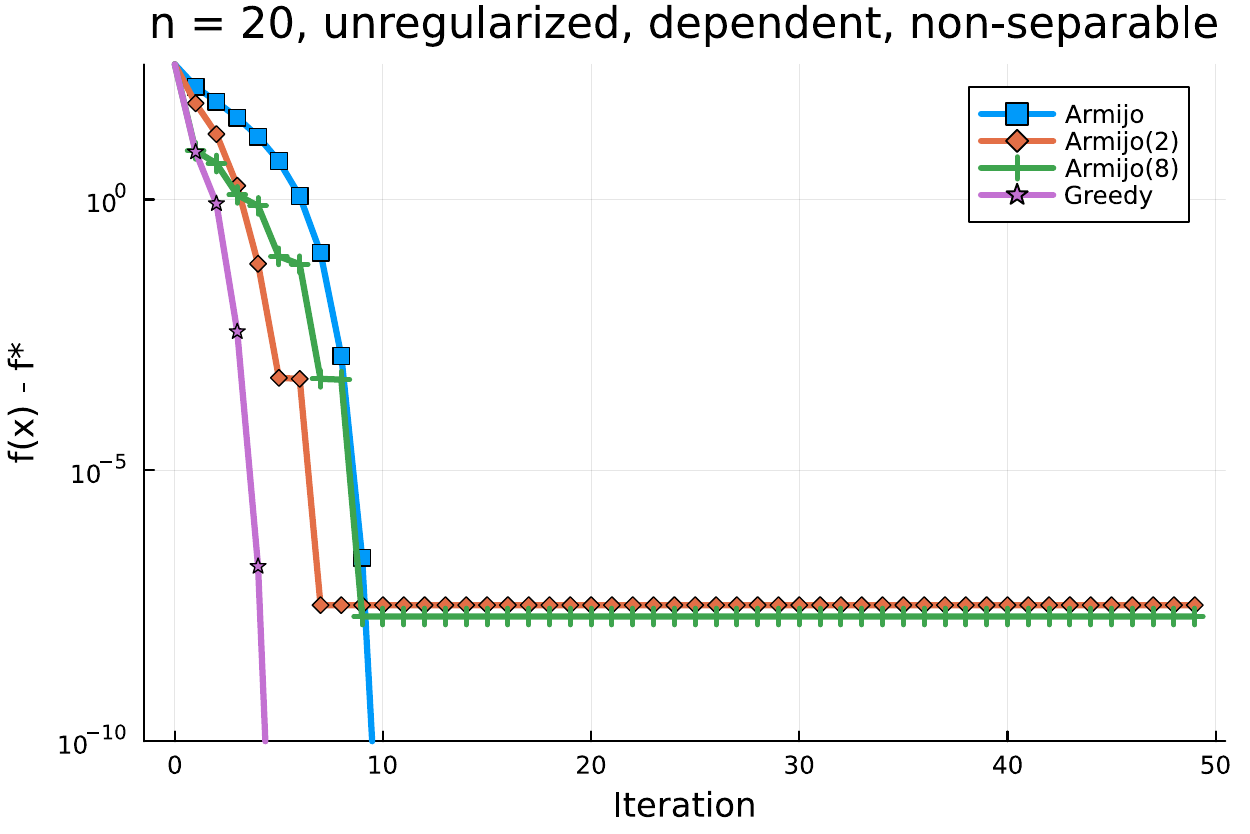}
	\includegraphics[width=.24\textwidth]{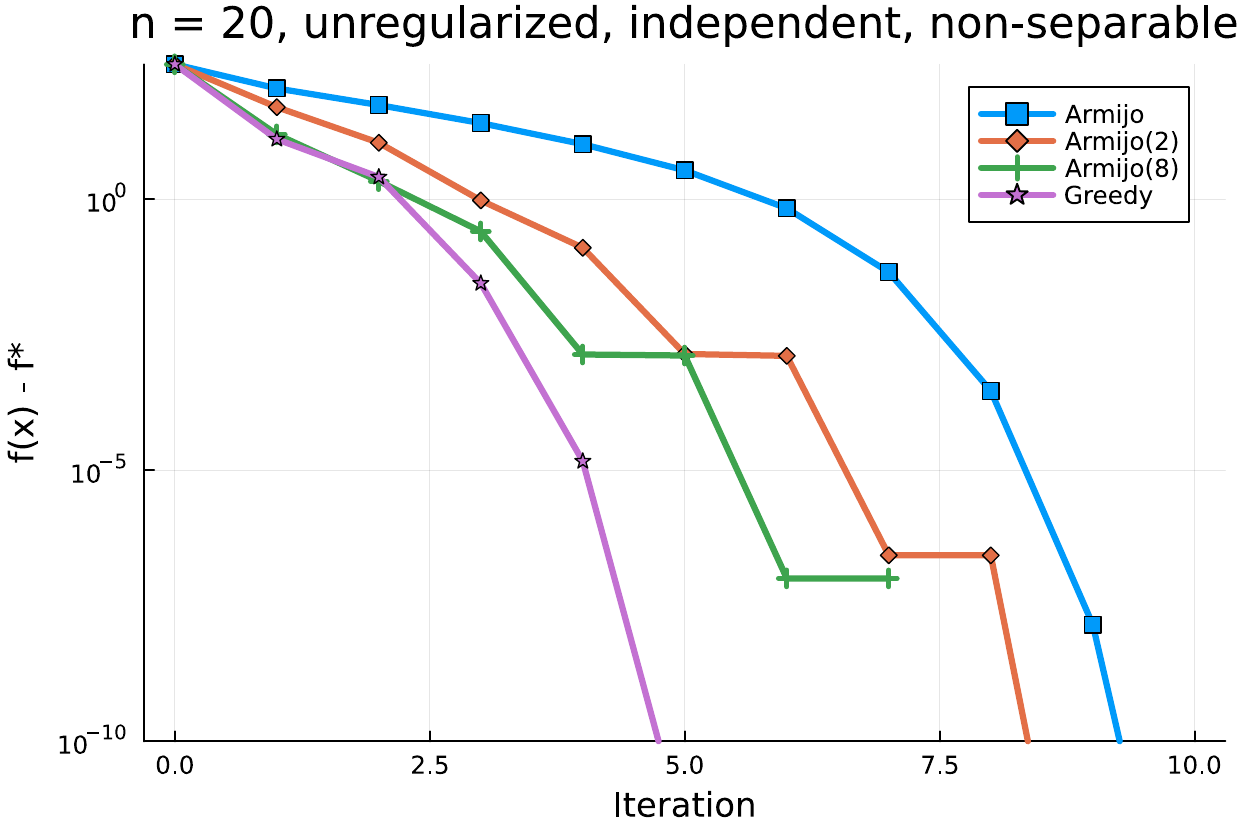}
	\includegraphics[width=.24\textwidth]{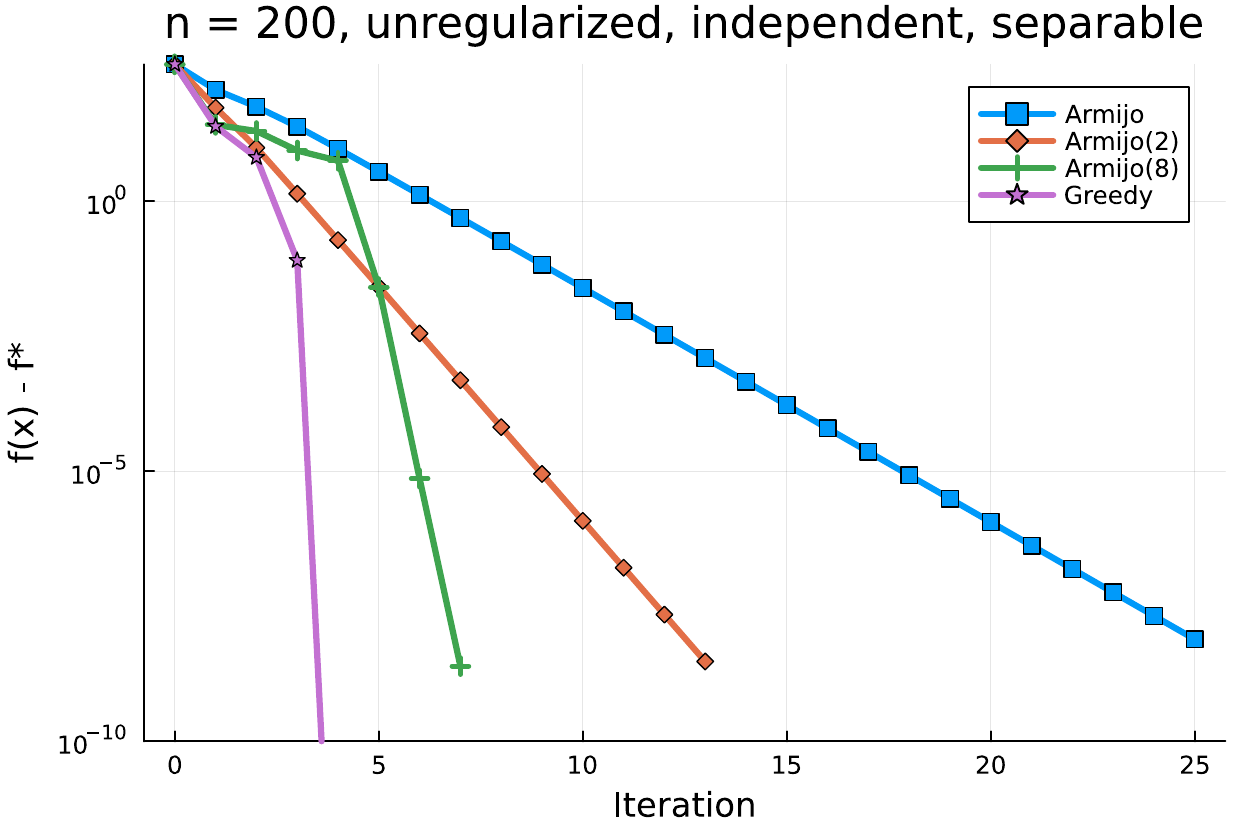}
	\includegraphics[width=.24\textwidth]{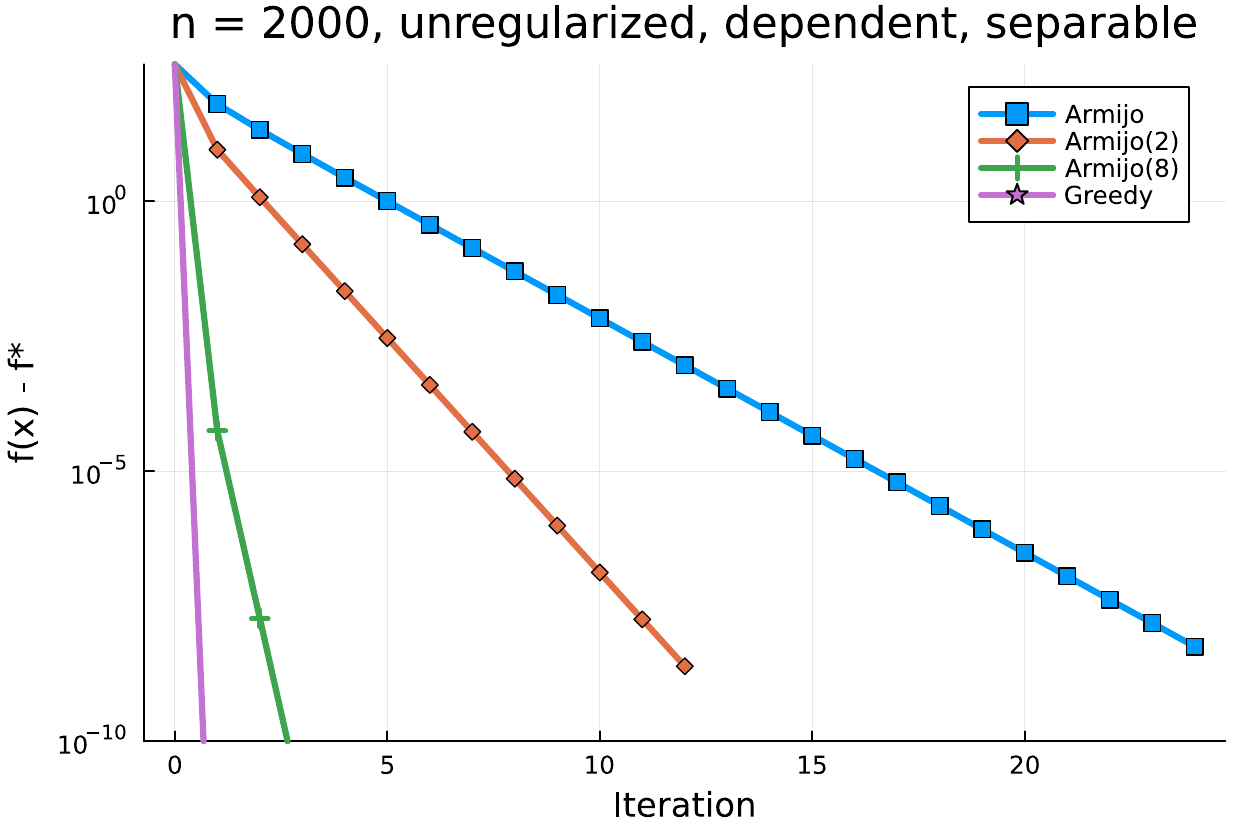}
	\caption{Comparison of greedy Newton to Newton with Armijo backtracking using three choices for the step size initialization (1, 2, and 8). We see that larger initializations can narrow or remove the gap with greedy Newton on early iterations, but can lead to poor performance on later iterations.
	}
	\label{fig:big}
\end{figure}

If Figure~\ref{fig:big} we repeat our synthetic experiments but explore different initializations of the Armijo backtracking line-search. In addition to the standard initialization of $\alpha=1$, we explored using a slightly larger initialization of $\alpha=2$ and a much larger initialization of $\alpha=8$. Note that $\alpha=8$ is larger than all the step sizes chosen by GN across iterations, except on the two unregularized seperable datasets. While these larger step sizes led to better performance on early iterations for many datasets, for some datasets they eventually led to much worse performance (due to using step sizes larger than 1 on later iterations).

\section{Greedy Newton with Cubic Regularization}
\label{app:cubic}

Given the Lipschitz constant of the Hessian $M$, Newton's method with cubic regularization~\cite{Nesterov2006} uses iterations of the form
\[
x_{k+1}^{C,M} \in \argmin_{y}\left\{f(x_k) + \nabla f(x_k)^T(y-x_k) + \frac{1}{2}(y-x_k)^T\nabla^2f(x_k)(y-x_k) + \frac{M}{6}\norm{y-x_k}^3\right\}.
\]
Unlike the classic Newton method, it is known that this method has a quadratic convergence rate directly in terms of function values. In particular, the method converges quadratically beginning from the first iteration where $f(x_k) - f(x_*) \leq \frac{M}{2\mu}$~\cite[see][Section~4.2.6]{Nesterov2018}.

There are several ways we could add a step size to this method:
\begin{enumerate}
	\item If $M$ is known, we could use iterations that take a step in the direction of a solution of the cubic sub-problem,
	\begin{equation}
		x_{k+1} = x_k +\alpha_k(x_{k+1}^{C,M} - x_k),
		\label{eq:cubic_search}
	\end{equation}
	where $\alpha_k$ is chosen to minimize the function value. The cubic step corresponds to choosing $\alpha_k=1$, but other values may decrease the objective function by a larger amount. Because we have $f(x_{k+1}) \leq f(x_k^C)$ with this method, it has the same radius of superlinear convergence as the basic Newton's method with cubic regularization.
	\item If $M$ is not known, it is common to use a backtracking procedure to set $M$~\cite{Nesterov2006}. In this case, we could alternate between backtracking to find an $M$ guaranteeing sufficient progress, and doing a line search on $\alpha_k$ in the direction of a solution of the cubic sub-problem~\eqref{eq:cubic_search} with the current approximation of $M$.
	\item If $M$ is not known, we could alternately search for the $M$ that minimizes the function value,
	\[
	M_k \in \argmin_{M}\left\{f(x_{k+1}^{C,M})\right\},
	\]
	and then use this $M_k$ in place of $M$ in the cubic update. A re-formulation of this update is using~\cite{Nesterov2006}
	\[
	x_{k+1} = x_k - (\nabla^2 f(x_k)^{-1} + \lambda I)^{-1}\nabla f(x_k),
	\]
	and choosing
	\[
	\lambda_k \in \argmin_{\lambda}f(x_k - (\nabla^2 f(x_k)^{-1} + \lambda I)^{-1}\nabla f(x_k)),
	\]
	which is a greedy version of a variation on the classic Levenberg-Marquardt update.
	This method decreases the function at least as much as any particular choice of $M$. Thus, it preserves the region of superlinear convergence without requiring us to know $M$ and would likely perform better in practice than using a fixed $M$. The disadvantage of this approach compared to line search methods is that it  involves additional operations with the Hessian. However, given the gradient and Hessian we could use a suitable factorization of the Hessian such as the Schur factorization or singular value decomposition~\cite{golub2013matrix} to implement this search in $O(n^3 + n^2\log(1/\epsilon))$ time plus the cost of evaluating $f$ $O(\log(1/\epsilon))$ times.
\end{enumerate}

\end{appendices}

\bigskip

\bmhead{Acknowledgements} We thank Frederik Kunstner and Nicolas Boumal for valuable discussions, and to the anonymous reviewers for suggestions that improved the paper. 

\bmhead{Funding} Betty Shea is funded by an NSERC Canada Graduate Scholarship. The work was partially supported by the Canada CIFAR AI Chair Program and NSERC Discovery Grant RGPIN-2022-036669.

\bmhead{Data availability}
The datasets used in this paper are publicly available at the minFunc package \citep{Schmidt2005} available at \url{https://www.cs.ubc.ca/~schmidtm/Software/minFunc.html}.

\bmhead{Conflict of interest} The authors have no financial or proprietary interests in any material discussed in this article.

\newpage
\bibliography{GN}

\end{document}